\documentclass[12pt,letterpaper]{amsart}
\usepackage{epic,eepic,latexsym, amssymb, amscd, amsfonts, xypic, floatflt}
 \newlength{\baseunit}               % the basic unit length
 \newcount{\numlines}                % depth of picture (in number of lines)
\newcommand{\point}{\vspace{3mm}\par \noindent \refstepcounter{subsection}{\bf \thesubsection.} }
\newcommand{\tpoint}[1]{\vspace{3mm}\par \noindent \refstepcounter{subsection}{\bf \thesubsection.} 
  {\em #1. ---} }
\newcommand{\epoint}[1]{\vspace{3mm}\par \noindent \refstepcounter{subsection}{\bf \thesubsection.} 
  {\em #1.} }
\newcommand{\bpoint}[1]{\vspace{3mm}\par \noindent \refstepcounter{subsection}{\bf \thesubsection.} 
  {\bf #1.} }
\newcommand{\bpf}{\noindent {\em Proof.  }}
\newcommand{\epf}{\qed \vspace{+10pt}}

\newcommand{\LR}{Geometric Littlewood-Richardson rule}

\newcommand{\A}{\mathbb{A}}

\newcommand{\proj}{\mathbb P}
\newcommand{\oh}{{\mathcal{O}}}

\newcommand{\al}{\alpha}

\newcommand{\be}{\beta}
\newcommand{\ga}{\gamma}
\newcommand{\de}{\delta}

\newcommand{\spam}{\operatorname{span}}

\newcommand{\Cl}{\operatorname{Cl}}

\newcommand{\init}{\rm{init}}
\newcommand{\final}{\rm{final}}

\newcommand{\poset}{\mathcal{P}}
\newcommand{\cq}{{\mathcal{Q}}}
\newcommand{\x}{\mathbf{x}}
\newcommand{\y}{\mathbf{y}}
\newcommand{\yp}{\mathbf{y'}}
\newcommand{\z}{\mathbf{z}}
\newcommand{\ff}{\mathbf{F}}  % red flag
\newcommand{\fm}{\mathbf{M}}  % black flag
\newcommand{\bn}{\bullet_{\text{next}}}
\newcommand{\bB}{\bullet \Box}
\newcommand{\cb}{\circ \bullet}
\newcommand{\cB}{\circ \Box}
\newcommand{\cbB} {\circ \bullet \Box}
\newcommand{\Xbb}{X_{\bullet} \cup X_{\bn}}

\newcommand{\Y}{\mathcal{Y}}

\newcommand{\pBox}{\poset_{\Box}}
\newcommand{\pcb}{p_{\circ \Box}}
\newcommand{\cs}{\circ_{\text{swap}}}
\newcommand{\cn}{\circ_{\text{stay}}}
\newcommand{\css}{\circ_{\text{sub}}}
\newcommand{\cbs}{\cs \bn}
\newcommand{\cbn}{\cn \bn}
\newcommand{\PF}{BS}  % was Q in earlier version!
\newcommand{\Qflag}{}  
\newcommand{\Sv}{\Omega}
\newcommand{\bu}{\cdot}  % This is a tiny bullet.
\newcommand{\ci}{{\sssize \circ}}
\newcommand{\sigkb}{Thm.~2.6}

\newcommand{\cited}{}

\newcommand{\secretnote}[1]{}
\newcommand{\notation}[1]{}
\newcommand{\lremind}[1]{{}}

\newcommand{\cut}[1]{}

\begin{document}
\pagestyle{plain}
\title{{\large {A geometric Littlewood-Richardson rule}}
%\footnote{1991 Mathematics Subject Classification:  Primary ?????, Secondary ?????}
}
\author{Ravi Vakil}
\address{Dept. of Mathematics, Stanford University, Stanford CA~94305--2125}
\email{vakil@math.stanford.edu}
\thanks{Partially supported by NSF Grant DMS--0228011, an AMS Centennial Fellowship, and an Alfred P. Sloan Research Fellowship.}
\date{Saturday, February 22, 2003.}
\begin{abstract}
We describe an explicit geometric
Littlewood-Richardson rule, interpreted as deforming the intersection of
two Schubert varieties so that they break into Schubert varieties.  
There are no restrictions on the base field, and all multiplicities
arising are 1; this is important
for applications.
This rule
should be seen as a generalization of Pieri's rule to arbitrary
Schubert classes, by way of explicit homotopies.  It has a
straightforward bijection to other Littlewood-Richardson rules, such
as tableaux and Knutson and Tao's puzzles.

%$K$-theory, suggested by Buch, which answers a question of Tao's on an
%extension of puzzles to $K$-theory.  The rule suggests a natural

This gives the first geometric proof and interpretation of the
Littlewood-Richardson rule. It has a host of geometric consequences,
described in \cite{si}.  The rule also has an interpretation in
$K$-theory, suggested by Buch, which gives an
extension of puzzles to $K$-theory.  The rule suggests a natural
approach to the open question of finding a Littlewood-Richardson rule
for the flag variety, leading to a conjecture, shown to be true up to
dimension 5.  Finally, the rule suggests approaches to similar open
problems, such as Littlewood-Richardson rules for the symplectic
Grassmannian and two-flag varieties.
\end{abstract}
\maketitle
\tableofcontents

{\parskip=12pt % closing bracket is just before the bibliography 

%{\em abstract for AG:}   Take the above abstract,
% and change [SI].  Is there anything wrong with having
% the abstract the 1st 2 pars of the paper?  Also,
% I've changed a couple of lines because of Terry Tao's response.
%
%We describe an explicit geometric
%Littlewood-Richardson rule, interpreted as deforming the intersection of
%two Schubert varieties so that they break into Schubert varieties.  
%There are no restrictions on the base field, and all multiplicities
%arising are 1; this is important
%for applications.
%This rule
%should be seen as a generalization of Pieri's rule to arbitrary
%Schubert classes, by way of explicit homotopies.  It has a
%straightforward bijection to other Littlewood-Richardson rules, such
%as tableaux and Knutson and Tao's puzzles.
%
%This gives the first geometric proof and interpretation of the
%Littlewood-Richardson rule. It has a host of geometric consequences,
%described in math.AG????.  The rule also has an interpretation in
%$K$-theory, suggested by Buch, which answers a question of Tao's on an
%extension of puzzles to $K$-theory.  The rule suggests a natural
%approach to the open question of finding a Littlewood-Richardson rule
%for the flag variety, leading to a conjecture, shown to be true up to
%dimension 5.  Finally, the rule suggests approaches to similar open
%problems, such as Littlewood-Richardson rules for the symplectic
%Grassmannian and two-flag varieties.

\section{Introduction}

The goal of this note is to describe an explicit geometric
Littlewood-Richardson rule, interpreted as deforming the intersection
of two Schubert varieties (with respect to transverse flags
$\ff_{\bu}$ and $\fm_{\bu}$) so that they break into Schubert
varieties.  There are no restrictions on the base field, and all
multiplicities arising are 1; this is important for applications.
This rule should be seen as a generalization of Pieri's rule to
arbitrary Schubert classes, by way of explicit homotopies.  It has a
straightforward bijection to other Littlewood-Richardson rules, such
as tableaux (Sect.~\ref{bijst}) and Knutson and Tao's puzzles
\cite{ktw, kt} (Sect.~\ref{checkerpuzzle}).

This gives the first geometric proof and interpretation of the
Littlewood-Richardson rule. It has a host of geometric consequences,
described in \cite{si}.  The rule also has an interpretation in
$K$-theory, suggested by Buch (Sect.~\ref{kcpt}), which suggests an
extension of puzzles to $K$-theory (Theorem~\ref{kglrp}), yielding a
triality of $K$-theory Littlewood-Richardson coefficients.  The rule
suggests a natural approach to the open question of finding a
Littlewood-Richardson rule for the flag variety, leading to a
conjecture, shown to be true up to dimension 5 (Sect.~\ref{flint}).
Finally, the rule suggests approaches to similar open problems, such
as Littlewood-Richardson rules for the symplectic Grassmannian and
two-flag varieties, and the quantum cohomology of the Grassmannian
(Sect.~\ref{open}).

The strategy is as follows.  We degenerate the ``{\bf M}oving flag'' 
 $\fm_{\bu}$ through a
series of codimension one degenerations (in $Fl(n)$) in a particular
way (the ``specialization order'') so that the cycle in $G(k,n)$
successively breaks into varieties we can easily understand
(``two-flag Schubert varieties'').  At each stage, the cycle breaks
into one or two pieces, and each piece appears with multiplicity one
(a key fact for applications).  At the end 
$\fm_{\bu}$ coincides
with the ``{\bf F}ixed flag'' $\ff_{\bu}$, and the limit cycle is a union of Schubert
varieties with respect to this flag. 
An explicit
example is shown in Figure~\ref{g24}.  ({\em Caution:} the theorem is
stated in terms of the affine description of the Grassmannian, but all
geometric descriptions are in terms of projective geometry.
Hence Figure~\ref{g24} shows a calculation in $G(2,4)$ using lines in
$\proj^3$.)  

We note that degeneration methods are a very old
technique.  See \cite{hilbert} for a historical discussion.
Sottile suggests that \cite{pieri} is an early example, proving
Pieri's formula using such methods; see also Hodge's proof \cite{hodge}.
More recent work by Sottile (especially \cite{lines}, dealing with 
$G(2,n)$, and \cite{explicitPieri}, concerning Pieri's formula in general)
provided inspiration for this work.
% I believe this degeneration
%approach is the same strategy originally used to prove both Pieri's
%formula and Monk's rule, but I will check to make sure.  

The degenerations can be described combinatorially, in terms of black
and white checkers on an $n \times n$ checkerboard.  Thus
Littlewood-Richardson coefficients count ``checker games''.  The input
is the data of two Schubert varieties, in the form of two $k$-subsets
of $\{1,\dots,n\}$; each move corresponds to moving black checkers
(encoding the position of $\fm_{\bu}$) in a certain way (the
``specialization order''), and determining then how the white checkers
(encoding the position of the $k$-plane) move.  The output is
interpreted as set of $k$-subsets.  Similarly, Schubert problems (i.e.
intersecting many Schubert classes) count ``checker tournaments''.
The checker rule is easy to use in practice, but somewhat awkward to
describe.  However, at all steps, (i) its geometric meaning is clear
at all steps, and (ii) there is a straightforward bijection to
partially completed tableaux and puzzles.  For many geometric
applications \cite{si}, the details of the combinatorial rule are
unnecessary.  R. Moriarty has written a program implementing this
rule.

\bpoint{Remarks on positive characteristic} We note that the only two
characteristic-dependent statements in the paper are invocations of
the Kleiman-Bertini theorem \cite{kleimanbertini}
(Sections~\ref{baltimore} and~\ref{wimbledon}).  Neither is used for
the proof of the main theorem (Theorem~\ref{glr3}).  They
will be replaced by a characteristic-free
generic smoothness theorem \cite[\sigkb]{si}
proved {\em using} the \LR.

\bpoint{Overview of paper}
Section~\ref{combdisc} is a description of
the rule in combinatorial terms; it contains no geometry.

Section~\ref{checkergeometry} describes the geometry related to
``two-flag Schubert varieties'' in the Grassmannian in the language of
checkers. Black checkers correspond to the relative position of the two flags; white checkers
correspond to the $k$-plane.  The main theorem (Theorem~\ref{glr3}) is stated
here.  Geometrically-minded readers may prefer to read some of
Section~\ref{checkergeometry} before Section~\ref{combdisc}.

The proof makes repeated use of smooth ``Bott-Samelson \Qflag varieties''
$\PF(\poset)$ corresponding to a ``planar poset'' $\poset$.  They
are simple objects, and some basic properties are described in
Section~\ref{bsv}.  Of particular importance is the Bott-Samelson \Qflag
variety $\PF(\pBox)$ parametrizing two hyperplanes and a codimension
two plane contained in both.  Each step in the degeneration is related to
$\PF(\pBox)$ via the key construction of Section~\ref{portland}.  

The proof of the \LR{} is given in Section~\ref{thmpf}; the strategy
is outlined in Section~\ref{strategy}.

Applications are discussed in Section~\ref{app2}.
Some are proved in detail, while others are merely sketched.
Further geometric applications are discussed in \cite{si}.

\epoint{Summary of notation} 
If $X \subset Y$, let $\Cl_Y X$ denote the closure in $Y$ of $X$.
Fix a base field $K$ (of any characteristic, not necessarily
algebraically closed), and non-negative integers $k \leq n$.  We work
in $G(k,n)$, the Grassmannian of dimension $k$ subspaces of 
$K^n$. \notation{$\Cl$, $K$}

We follow the notation of \cite{fulton}.  For example, let $\Sv_w$
denote the Schubert class in $H^*(G(k,n))$ or $A^*(G(k,n))$ (resp.
$H^*(Fl(n))$ or $A^* (Fl(k,n))$), where $w$ is a partition (resp. a
permutation).  Let $\Sv_w(F_{\bu})$ (resp. $\Sv_w^{\ci}(F_{\bu})$) be
the closed (resp. open) Schubert variety with respect to the flag
$F_{\bu}$.  In $Fl(n)$,
$$
\Sv_w(F_{\bu}) = \left\{ [V_{\bu}] \in Fl(n) : \dim ( V_p \cap F_q) \geq
\# \{ i \leq p : w(i) > n-q \} \text{ for all } p, q \right\}.$$
The Schur polynomials are denoted $s_{\bu}$.

Table~\ref{sumnot} is a summary of notation introduced in the article.

\begin{table}
\begin{tabular}{|l|l|} \hline
section  & notation \\ introduced & \\ \hline
\ref{firstblood}-- \ref{rulewhite} & configuration of checkers, 
  specialization order, 
  descending \\
& and rising checkers,
$\bullet$, $\bullet_{\init}$, $\bullet_{\final}$, $\bn$, 
 happy, $\circ$, $\cb$, $\circ_{A,B}$, critical \\
& row, critical diagonal, Phase 1, swap, blocker, Phase 2, 
$\cn, \cs$, $\circ_S$ \\
\ref{midsort} & mid-sort \\
\ref{domdefhere} & $\fm_{\bu}$, $\ff_{\bu}$, dominate, $X_{\bullet}$ \\
\ref{eugene} & 
open and closed two-flag
Schubert varieties $Y_{\cb}$ and $\overline{Y}_{\cb}$, \\
& universal two-flag Schubert variety,
$X_{\cb}$, $\overline{X}_{\cb}$ \\
\ref{hamilton} & $D_X$ \\
\ref{bsvs} & dimension, planar poset,  $\poset$, quadrilateral, southwest
and \\ & northeast border, Bott-Samelson variety, $\PF(\poset)$ \\
\ref{bsvss}--\ref{augusta} & strata of a Bott-Samelson variety, 
$\poset_{\bullet},
  \poset_{\circ}$,  $\PF(\poset \rightarrow \cq)$, $p_{\cb}$ \\
\ref{portland}--\ref{secdesc} & $\pBox$, $\x$, $\y$, $\yp$, $\z$, $V_{\x}$, 
$V_{\y}$, $V_{\yp}$, $V_{\z}$, $D_\Box$, 
$p_{\bullet \Box}$, $(r,c)$, $S_{\bu}$, $T_{\bu}$, $X_{\bB}$, $X_{\bB}$-position \\
\ref{strategy} & geometrically (ir)relevant \\
\ref{providence} & $Z$ \\
\ref{pit}--\ref{oakpark} & $\mathbf{a}$,$\mathbf{a'}$, ${\mathbf{a''}}$, 
left and right good quadrilaterals, columns, 
 $D_S^Z$, $D_S$, $Q_{\ell}$, $\poset'_{\circ}$ \\
\hline
\end{tabular}
\caption{Summary of notation \label{sumnot} \lremind{sumnot}}
\end{table}

\bpoint{Acknowledgments} The author is grateful to A. Buch and A.
Knutson for patiently explaining the combinatorial, geometric, and
representation-theoretic ideas behind this problem, and for comments
on earlier versions.  Any remaining misunderstandings are purely due to
the author.  The author also thanks L. Chen, W. Fulton, R. Moriarty,
and F. Sottile.

\section{The \LR:  Combinatorial
description} \label{combdisc}
\lremind{combdisc}

\point \label{firstblood} \lremind{firstblood}
Littlewood-Richardson coefficients $c_{\al \be}^{\ga}$ will be
counted in their guise of
structure constants of the cohomology (or Chow) ring
of the Grassmannian $G(k,n)$, which in turn will be 
interpreted in terms 
of {\em checker games}, involving $n$ black and $k$ white
checkers on an $n \times n$ board.  The rows and columns of
the board are numbered as in Figure~\ref{raleigh}; $(r,c)$ will
denote the square in row $r$ and column $c$.
A set of checkers on the board will be called a {\em configuration}
of checkers.\notation{configuration of checkers}

\begin{figure}[ht]
\begin{center} \setlength{\unitlength}{0.00083333in}
\begingroup\makeatletter\ifx\SetFigFont\undefined%
\gdef\SetFigFont#1#2#3#4#5{%
  \reset@font\fontsize{#1}{#2pt}%
  \fontfamily{#3}\fontseries{#4}\fontshape{#5}%
  \selectfont}%
\fi\endgroup%
{\renewcommand{\dashlinestretch}{30}
\begin{picture}(912,861)(0,-10)
\path(300,612)(900,612)(900,12)
	(300,12)(300,612)
\path(450,612)(450,12)
\path(600,612)(600,12)
\path(750,612)(750,12)
\path(300,462)(900,462)
\path(300,312)(900,312)
\path(300,162)(900,162)
\put(337,762){\makebox(0,0)[lb]{\smash{{{\SetFigFont{8}{9.6}{\rmdefault}{\mddefault}{\updefault}$1$}}}}}
\put(0,37){\makebox(0,0)[lb]{\smash{{{\SetFigFont{8}{9.6}{\rmdefault}{\mddefault}{\updefault}$n$}}}}}
\put(787,762){\makebox(0,0)[lb]{\smash{{{\SetFigFont{8}{9.6}{\rmdefault}{\mddefault}{\updefault}$n$}}}}}
\put(0,487){\makebox(0,0)[lb]{\smash{{{\SetFigFont{8}{9.6}{\rmdefault}{\mddefault}{\updefault}$1$}}}}}
\put(525,762){\makebox(0,0)[lb]{\smash{{{\SetFigFont{8}{9.6}{\rmdefault}{\mddefault}{\updefault}$\cdots$}}}}}
\put(20,237){\makebox(0,0)[lb]{\smash{{{\SetFigFont{8}{9.6}{\rmdefault}{\mddefault}{\updefault}$\vdots$}}}}}
\end{picture}
} \end{center}
\caption{\lremind{raleigh}}
\label{raleigh}
\end{figure}  

\bpoint{Black checkers}
\label{seattle} \lremind{seattle}
The moves of the
black checkers are prescribed in advance; 
the Littlewood-Richardson coefficients are counted
by the possible accompanying moves of the white checkers.

No two checkers of the same color are in the same row or column.
Hence the positions of the black checkers correspond to permutations
by the following bijection:   if the black checker in column $c$
is in row $r$, then the permutation sends $n+1-c$ to $r$.
For example, see Figure~\ref{butte}.

The initial position of the black checkers corresponds to the
identity permutation in $S_n$, and the final position corresponds
to the longest word $w_0$.  The intermediate positions in the 
checker game correspond to
partial factorizations from the right of $w_0$:
$$
w_0 = e_{n-1}  \cdots e_2 e_1 \quad \cdots  \quad  e_{n-1} e_{n-2} e_{n-3} \quad e_{n-1} e_{n-2}  \quad e_{n-1}
.$$
For example, in Figure~\ref{butte} shows the six moves of the
black checkers for $n=4$, along with the corresponding permutations.
(The geometric interpretation will be explained in Sect.~\ref{madison}.)
In essence this is a bubble-sort.
Call this the {\em specialization order} in the Bruhat order.
 Figure~\ref{boston} shows
a typical checker configuration in the specialization order.
Each  move will involve moving one checker one row down
(call this the {\em descending checker}), and another checker (call this the {\em rising checker})\notation{specialization order; descending and rising checker
(should descending be ``falling'')} 
one row up, as shown in the figure.

Note that this word neither begins nor ends with the corresponding word
for $n-1$, making induction impossible (see Sect.~\ref{cap}). 

We will use the notation $\bullet$ to represent a configuration
of black checkers.  
Denote the initial configuration by  $\bullet_{\init}$, and
the final configuration by $\bullet_{\final}$.
If $\bullet$ is one of the  configurations in the
specialization order, and $\bullet \neq \bullet_{\final}$,
then $\bn$ will denote the next configuration in 
the specialization order.\notation{$\bullet$, $\bullet_{\init}$, $\bullet_{\final}$, $\bn$}

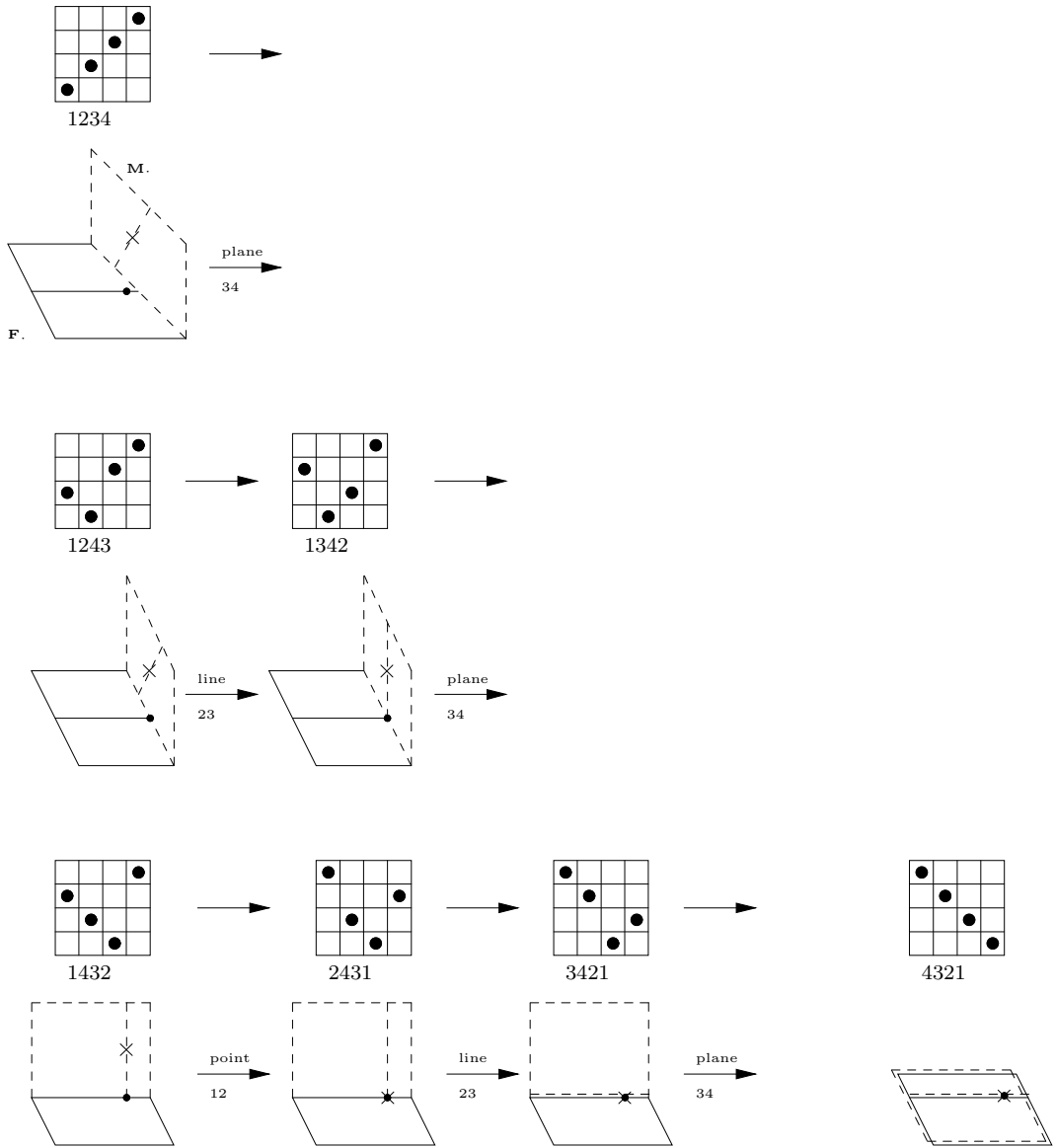
\begin{figure}[ht]
\begin{center} \setlength{\unitlength}{0.00083333in}
\begingroup\makeatletter\ifx\SetFigFont\undefined%
\gdef\SetFigFont#1#2#3#4#5{%
  \reset@font\fontsize{#1}{#2pt}%
  \fontfamily{#3}\fontseries{#4}\fontshape{#5}%
  \selectfont}%
\fi\endgroup%
{\renewcommand{\dashlinestretch}{30}
\begin{picture}(6624,7239)(0,-10)
\put(2412,312){\blacken\ellipse{40}{40}}
\put(2412,312){\ellipse{40}{40}}
\path(1812,312)(1962,12)(2712,12)
	(2562,312)(1812,312)
\dashline{60.000}(1812,312)(1812,912)(2562,912)(2562,312)
\dashline{60.000}(2412,912)(2412,312)
\put(3912,312){\blacken\ellipse{40}{40}}
\put(3912,312){\ellipse{40}{40}}
\path(3312,312)(3462,12)(4212,12)
	(4062,312)(3312,312)
\dashline{60.000}(3312,312)(3312,912)(4062,912)(4062,312)
\dashline{60.000}(3315,334)(4065,334)
\path(312,1812)(912,1812)(912,1212)
	(312,1212)(312,1812)
\path(612,1812)(612,1212)
\path(462,1812)(462,1212)
\path(762,1812)(762,1212)
\path(312,1512)(912,1512)
\path(312,1662)(912,1662)
\path(312,1362)(912,1362)
\path(1962,1812)(2562,1812)(2562,1212)
	(1962,1212)(1962,1812)
\path(2262,1812)(2262,1212)
\path(2112,1812)(2112,1212)
\path(2412,1812)(2412,1212)
\path(1962,1512)(2562,1512)
\path(1962,1662)(2562,1662)
\path(1962,1362)(2562,1362)
\path(3462,1812)(4062,1812)(4062,1212)
	(3462,1212)(3462,1812)
\path(3762,1812)(3762,1212)
\path(3612,1812)(3612,1212)
\path(3912,1812)(3912,1212)
\path(3462,1512)(4062,1512)
\path(3462,1662)(4062,1662)
\path(3462,1362)(4062,1362)
\path(312,7212)(912,7212)(912,6612)
	(312,6612)(312,7212)
\path(612,7212)(612,6612)
\path(462,7212)(462,6612)
\path(762,7212)(762,6612)
\path(312,6912)(912,6912)
\path(312,7062)(912,7062)
\path(312,6762)(912,6762)
\path(312,4512)(912,4512)(912,3912)
	(312,3912)(312,4512)
\path(612,4512)(612,3912)
\path(462,4512)(462,3912)
\path(762,4512)(762,3912)
\path(312,4212)(912,4212)
\path(312,4362)(912,4362)
\path(312,4062)(912,4062)
\path(1812,4512)(2412,4512)(2412,3912)
	(1812,3912)(1812,4512)
\path(2112,4512)(2112,3912)
\path(1962,4512)(1962,3912)
\path(2262,4512)(2262,3912)
\path(1812,4212)(2412,4212)
\path(1812,4362)(2412,4362)
\path(1812,4062)(2412,4062)
\put(6310,324){\blacken\ellipse{40}{40}}
\put(6310,324){\ellipse{40}{40}}
\path(5637,462)(5862,12)(6612,12)
	(6387,462)(5637,462)
\dashline{60.000}(5715,334)(6465,334)
\path(5712,312)(6462,312)
\dashline{60.000}(5598,487)(5823,37)(6573,37)
	(6348,487)(5598,487)
\path(5712,1812)(6312,1812)(6312,1212)
	(5712,1212)(5712,1812)
\path(6012,1812)(6012,1212)
\path(5862,1812)(5862,1212)
\path(6162,1812)(6162,1212)
\path(5712,1512)(6312,1512)
\path(5712,1662)(6312,1662)
\path(5712,1362)(6312,1362)
\put(6237,1287){\blacken\ellipse{74}{74}}
\put(6237,1287){\ellipse{74}{74}}
\put(6087,1437){\blacken\ellipse{74}{74}}
\put(6087,1437){\ellipse{74}{74}}
\put(5937,1587){\blacken\ellipse{74}{74}}
\put(5937,1587){\ellipse{74}{74}}
\put(5787,1737){\blacken\ellipse{74}{74}}
\put(5787,1737){\ellipse{74}{74}}
\put(5787,1062){\makebox(0,0)[lb]{\smash{{{\SetFigFont{8}{9.6}{\rmdefault}{\mddefault}{\updefault}$4321$}}}}}
\path(764,5788)(839,5713)
\path(839,5788)(764,5713)
\path(869,3050)(944,2975)
\path(944,3050)(869,2975)
\path(2370,3047)(2445,2972)
\path(2445,3047)(2370,2972)
\path(723,653)(798,578)
\path(798,653)(723,578)
\path(2376,349)(2451,274)
\path(2451,349)(2376,274)
\path(3874,349)(3949,274)
\path(3949,349)(3874,274)
\path(6272,361)(6347,286)
\path(6347,361)(6272,286)
\put(837,1737){\blacken\ellipse{74}{74}}
\put(837,1737){\ellipse{74}{74}}
\put(387,1587){\blacken\ellipse{74}{74}}
\put(387,1587){\ellipse{74}{74}}
\put(537,1437){\blacken\ellipse{74}{74}}
\put(537,1437){\ellipse{74}{74}}
\put(687,1287){\blacken\ellipse{74}{74}}
\put(687,1287){\ellipse{74}{74}}
\put(2487,1587){\blacken\ellipse{74}{74}}
\put(2487,1587){\ellipse{74}{74}}
\put(2037,1737){\blacken\ellipse{74}{74}}
\put(2037,1737){\ellipse{74}{74}}
\put(2187,1437){\blacken\ellipse{74}{74}}
\put(2187,1437){\ellipse{74}{74}}
\put(2337,1287){\blacken\ellipse{74}{74}}
\put(2337,1287){\ellipse{74}{74}}
\put(3987,1437){\blacken\ellipse{74}{74}}
\put(3987,1437){\ellipse{74}{74}}
\put(3837,1287){\blacken\ellipse{74}{74}}
\put(3837,1287){\ellipse{74}{74}}
\put(3687,1587){\blacken\ellipse{74}{74}}
\put(3687,1587){\ellipse{74}{74}}
\put(3537,1737){\blacken\ellipse{74}{74}}
\put(3537,1737){\ellipse{74}{74}}
\put(837,7137){\blacken\ellipse{74}{74}}
\put(837,7137){\ellipse{74}{74}}
\put(687,6987){\blacken\ellipse{74}{74}}
\put(687,6987){\ellipse{74}{74}}
\put(537,6837){\blacken\ellipse{74}{74}}
\put(537,6837){\ellipse{74}{74}}
\put(387,6687){\blacken\ellipse{74}{74}}
\put(387,6687){\ellipse{74}{74}}
\put(387,4137){\blacken\ellipse{74}{74}}
\put(387,4137){\ellipse{74}{74}}
\put(537,3987){\blacken\ellipse{74}{74}}
\put(537,3987){\ellipse{74}{74}}
\put(687,4287){\blacken\ellipse{74}{74}}
\put(687,4287){\ellipse{74}{74}}
\put(837,4437){\blacken\ellipse{74}{74}}
\put(837,4437){\ellipse{74}{74}}
\put(1887,4287){\blacken\ellipse{74}{74}}
\put(1887,4287){\ellipse{74}{74}}
\put(2187,4137){\blacken\ellipse{74}{74}}
\put(2187,4137){\ellipse{74}{74}}
\put(2037,3987){\blacken\ellipse{74}{74}}
\put(2037,3987){\ellipse{74}{74}}
\put(2337,4437){\blacken\ellipse{74}{74}}
\put(2337,4437){\ellipse{74}{74}}
\put(762,5412){\blacken\ellipse{40}{40}}
\put(762,5412){\ellipse{40}{40}}
\put(912,2712){\blacken\ellipse{40}{40}}
\put(912,2712){\ellipse{40}{40}}
\put(2412,2712){\blacken\ellipse{40}{40}}
\put(2412,2712){\ellipse{40}{40}}
\put(762,312){\blacken\ellipse{40}{40}}
\put(762,312){\ellipse{40}{40}}
\path(4287,1512)(4737,1512)
\blacken\path(4617.000,1482.000)(4737.000,1512.000)(4617.000,1542.000)(4617.000,1482.000)
\path(4287,462)(4737,462)
\blacken\path(4617.000,432.000)(4737.000,462.000)(4617.000,492.000)(4617.000,432.000)
\path(2787,462)(3237,462)
\blacken\path(3117.000,432.000)(3237.000,462.000)(3117.000,492.000)(3117.000,432.000)
\path(2787,1512)(3237,1512)
\blacken\path(3117.000,1482.000)(3237.000,1512.000)(3117.000,1542.000)(3117.000,1482.000)
\path(1212,1512)(1662,1512)
\blacken\path(1542.000,1482.000)(1662.000,1512.000)(1542.000,1542.000)(1542.000,1482.000)
\path(1212,462)(1662,462)
\blacken\path(1542.000,432.000)(1662.000,462.000)(1542.000,492.000)(1542.000,432.000)
\path(1287,6912)(1737,6912)
\blacken\path(1617.000,6882.000)(1737.000,6912.000)(1617.000,6942.000)(1617.000,6882.000)
\path(1287,5562)(1737,5562)
\blacken\path(1617.000,5532.000)(1737.000,5562.000)(1617.000,5592.000)(1617.000,5532.000)
\path(1137,4212)(1587,4212)
\blacken\path(1467.000,4182.000)(1587.000,4212.000)(1467.000,4242.000)(1467.000,4182.000)
\path(1137,2862)(1587,2862)
\blacken\path(1467.000,2832.000)(1587.000,2862.000)(1467.000,2892.000)(1467.000,2832.000)
\path(2712,4212)(3162,4212)
\blacken\path(3042.000,4182.000)(3162.000,4212.000)(3042.000,4242.000)(3042.000,4182.000)
\path(2712,2862)(3162,2862)
\blacken\path(3042.000,2832.000)(3162.000,2862.000)(3042.000,2892.000)(3042.000,2832.000)
\path(537,5712)(12,5712)(312,5112)(1137,5112)
\dashline{60.000}(912,5937)(687,5562)
\dashline{60.000}(537,5712)(1137,5112)(1137,5712)
	(537,6312)(537,5712)
\path(162,5412)(837,5412)
\path(762,3012)(162,3012)(462,2412)(1062,2412)
\dashline{60.000}(762,3012)(1062,2412)(1062,3012)
	(762,3612)(762,3012)
\path(312,2712)(912,2712)
\dashline{60.000}(987,3162)(837,2862)
\path(2262,3012)(1662,3012)(1962,2412)(2562,2412)
\dashline{60.000}(2262,3012)(2562,2412)(2562,3012)
	(2262,3612)(2262,3012)
\path(1812,2712)(2412,2712)
\dashline{60.000}(2412,3312)(2412,2712)
\path(162,312)(312,12)(1062,12)
	(912,312)(162,312)
\dashline{60.000}(162,312)(162,912)(912,912)(912,312)
\dashline{60.000}(762,912)(762,312)
\put(1287,537){\makebox(0,0)[lb]{\smash{{{\SetFigFont{5}{6.0}{\rmdefault}{\mddefault}{\updefault}point}}}}}
\put(2862,537){\makebox(0,0)[lb]{\smash{{{\SetFigFont{5}{6.0}{\rmdefault}{\mddefault}{\updefault}line}}}}}
\put(4362,537){\makebox(0,0)[lb]{\smash{{{\SetFigFont{5}{6.0}{\rmdefault}{\mddefault}{\updefault}plane}}}}}
\put(387,1062){\makebox(0,0)[lb]{\smash{{{\SetFigFont{8}{9.6}{\rmdefault}{\mddefault}{\updefault}$1432$}}}}}
\put(2037,1062){\makebox(0,0)[lb]{\smash{{{\SetFigFont{8}{9.6}{\rmdefault}{\mddefault}{\updefault}$2431$}}}}}
\put(3537,1062){\makebox(0,0)[lb]{\smash{{{\SetFigFont{8}{9.6}{\rmdefault}{\mddefault}{\updefault}$3421$}}}}}
\put(1287,312){\makebox(0,0)[lb]{\smash{{{\SetFigFont{5}{6.0}{\rmdefault}{\mddefault}{\updefault}$12$}}}}}
\put(2862,312){\makebox(0,0)[lb]{\smash{{{\SetFigFont{5}{6.0}{\rmdefault}{\mddefault}{\updefault}$23$}}}}}
\put(4362,312){\makebox(0,0)[lb]{\smash{{{\SetFigFont{5}{6.0}{\rmdefault}{\mddefault}{\updefault}$34$}}}}}
\put(762,6162){\makebox(0,0)[lb]{\smash{{{\SetFigFont{5}{6.0}{\rmdefault}{\mddefault}{\updefault}${\mathbf{M}}_{\cdot}$}}}}}
\put(387,6462){\makebox(0,0)[lb]{\smash{{{\SetFigFont{8}{9.6}{\rmdefault}{\mddefault}{\updefault}$1234$}}}}}
\put(12,5112){\makebox(0,0)[lb]{\smash{{{\SetFigFont{5}{6.0}{\rmdefault}{\mddefault}{\updefault}${\mathbf{F}}_{\cdot}$}}}}}
\put(1362,5637){\makebox(0,0)[lb]{\smash{{{\SetFigFont{5}{6.0}{\rmdefault}{\mddefault}{\updefault}plane}}}}}
\put(1362,5412){\makebox(0,0)[lb]{\smash{{{\SetFigFont{5}{6.0}{\rmdefault}{\mddefault}{\updefault}$34$}}}}}
\put(1212,2937){\makebox(0,0)[lb]{\smash{{{\SetFigFont{5}{6.0}{\rmdefault}{\mddefault}{\updefault}line}}}}}
\put(2787,2937){\makebox(0,0)[lb]{\smash{{{\SetFigFont{5}{6.0}{\rmdefault}{\mddefault}{\updefault}plane}}}}}
\put(387,3762){\makebox(0,0)[lb]{\smash{{{\SetFigFont{8}{9.6}{\rmdefault}{\mddefault}{\updefault}$1243$}}}}}
\put(1887,3762){\makebox(0,0)[lb]{\smash{{{\SetFigFont{8}{9.6}{\rmdefault}{\mddefault}{\updefault}$1342$}}}}}
\put(1212,2712){\makebox(0,0)[lb]{\smash{{{\SetFigFont{5}{6.0}{\rmdefault}{\mddefault}{\updefault}$23$}}}}}
\put(2787,2712){\makebox(0,0)[lb]{\smash{{{\SetFigFont{5}{6.0}{\rmdefault}{\mddefault}{\updefault}$34$}}}}}
\end{picture}
} \end{center}
\caption{The {\em specialization order} for $n=4$, in terms
of checkers and permutations, with the
geometric interpretation of Section~\ref{checkergeometry} \lremind{butte}}
\label{butte}
\end{figure} 

\bpoint{White checkers} 
\label{rulewhite} \lremind{rulewhite}
At each stage of the game, each white checker
has the following property: there is a black checker in the same
square or in a square above it, and there is a black checker in the
same square or in a square to the left of it.  A white
checker satisfying this property is said to be {\em happy}.
(In particular, a white checker is happy if it is  in the same square as a
black checker.)
A configuration of white checkers will often be denoted $\circ$,
and a configuration of  white  and black checkers
will often be denoted $\cb$.\notation{happy, $\circ$, $\cb$}

The initial position of the white checkers depends on two subsets
$A = \{ a_1, \dots, a_k \}$ and $B = \{ b_1, \dots, b_k \}$
of $\{ 1, \dots, n \}$, where $a_1 < \dots < a_k$ and $b_1 < \dots < b_k$.
The $k$ white checkers are in the squares $(a_1, b_k)$, $(a_2, b_{k-1})$, 
\dots, $(a_k, b_1)$.
We denote this configuration by $\circ_{A,B}$.\notation{$\circ_{A,B}$}
If (and only if) any of these white checkers are not {happy} (i.e. if $a_i + b_{k+1-i} \leq n$ for some $i$), then
there are no checker games corresponding to $(A,B)$.  (This
will mean that  all corresponding  Littlewood-Richardson coefficients are 0.
Geometrically, this means that the two Schubert varieties do not intersect.)
This happens, for example, if $n=2$ and $A=B= \{ 1 \}$, computing
the intersection of two general points in $\proj^1$.

For each move of the black checkers, there will be one or two 
possible moves of the white checkers, which we describe next.
Define the {\em critical row} and
the {\em critical diagonal} as in Figure~\ref{boston}.\notation{critical row, critical diagonal}
The movement of the white checkers takes place in two phases. 

\begin{figure}[ht]
\begin{center} \setlength{\unitlength}{0.00083333in}
\begingroup\makeatletter\ifx\SetFigFont\undefined%
\gdef\SetFigFont#1#2#3#4#5{%
  \reset@font\fontsize{#1}{#2pt}%
  \fontfamily{#3}\fontseries{#4}\fontshape{#5}%
  \selectfont}%
\fi\endgroup%
{\renewcommand{\dashlinestretch}{30}
\begin{picture}(1513,1558)(0,-10)
\put(562.500,368.500){\arc{706.701}{2.5822}{4.1530}}
\blacken\path(336.396,619.681)(375.000,668.000)(318.199,643.533)(336.396,619.681)
\put(1334.606,796.019){\arc{910.169}{3.4427}{5.0962}}
\blacken\path(907.318,992.412)(900.000,931.000)(935.357,981.744)(907.318,992.412)
\put(1200,1456){\blacken\ellipse{76}{76}}
\put(1200,1456){\ellipse{76}{76}}
\put(1050,1306){\blacken\ellipse{76}{76}}
\put(1050,1306){\ellipse{76}{76}}
\put(150,1156){\blacken\ellipse{76}{76}}
\put(150,1156){\ellipse{76}{76}}
\put(300,1006){\blacken\ellipse{76}{76}}
\put(300,1006){\ellipse{76}{76}}
\put(900,856){\blacken\ellipse{76}{76}}
\put(900,856){\ellipse{76}{76}}
\put(450,706){\blacken\ellipse{76}{76}}
\put(450,706){\ellipse{76}{76}}
\put(600,556){\blacken\ellipse{76}{76}}
\put(600,556){\ellipse{76}{76}}
\put(750,406){\blacken\ellipse{76}{76}}
\put(750,406){\ellipse{76}{76}}
\path(75,1531)(1275,1531)(1275,331)
	(75,331)(75,1531)
\path(225,1531)(225,331)
\path(375,1531)(375,331)
\path(525,1531)(525,331)
\path(675,1531)(675,331)
\path(825,1531)(825,331)
\path(975,1531)(975,331)
\path(1125,1531)(1125,331)
\path(75,1381)(1275,1381)
\path(75,1231)(1275,1231)
\path(75,1081)(1275,1081)
\path(75,931)(1275,931)
\path(75,781)(1275,781)
\path(75,631)(1275,631)
\path(75,481)(1275,481)
\dashline{60.000}(788,968)(1313,968)(1313,743)
	(788,743)(788,968)
\dashline{60.000}(338,818)(413,818)(863,368)
	(863,293)(788,293)(338,743)(338,818)
\path(450,743)(450,893)
\blacken\path(465.000,833.000)(450.000,893.000)(435.000,833.000)(465.000,833.000)
\path(900,818)(900,668)
\blacken\path(885.000,728.000)(900.000,668.000)(915.000,728.000)(885.000,728.000)
\put(0,31){\makebox(0,0)[lb]{\smash{{{\SetFigFont{8}{9.6}{\rmdefault}{\mddefault}{\updefault}rising checker}}}}}
\put(750,181){\makebox(0,0)[lb]{\smash{{{\SetFigFont{8}{9.6}{\rmdefault}{\mddefault}{\updefault}critical diagonal}}}}}
\put(1500,1081){\makebox(0,0)[lb]{\smash{{{\SetFigFont{8}{9.6}{\rmdefault}{\mddefault}{\updefault}descending checker}}}}}
\put(1425,806){\makebox(0,0)[lb]{\smash{{{\SetFigFont{8}{9.6}{\rmdefault}{\mddefault}{\updefault}critical row}}}}}
\end{picture}
} \end{center}
\caption{The critical row and the critical diagonal  \lremind{boston}}
\label{boston}
\end{figure}
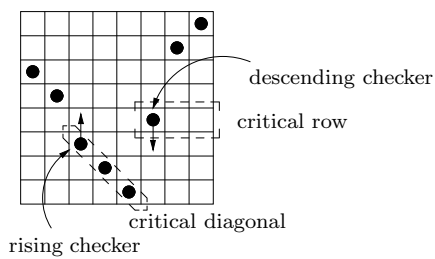

{\em Phase 1}\notation{Phase 1} depends on the answers to the two questions:
 Where (if anywhere) is the white checker in the critical row?
 Where (if anywhere) is the highest white checker in the critical
diagonal?  (There may be several white checkers in the critical
diagonal.)

Based on the answer to these questions, the pair of white checkers
either swap\notation{swap} rows (i.e. if the white checkers start at $(r_1,c_1)$ and
$(r_2, c_2)$, then they will end at $(r_2, c_1)$ and $(r_1, c_2)$), or
they stay where they are, according to Table~\ref{keywest}.  The central entry
of the table
is the only time when there is a possibility for choice: the pair of
white checkers can stay, or {\em if there are no white checkers
  in the rectangle between them} they can swap. 
Call white checkers in this rectangle {\em blockers}.\notation{blocker}
See Figure~\ref{blockerexample}
for an example of a blocker.

\begin{table}
\begin{tabular}{|c|c||c|c|c|}
\hline
\multicolumn{2}{|c||}{} & \multicolumn{3}{c|}{\em White checker in critical row?} \\  \cline{3-5}
\multicolumn{2}{|c||}{} 
 & yes, in  descending & yes, elsewhere & no  \\ 
\multicolumn{2}{|c||}{} 
 & checker's square & & 
\\ \hline \hline
{\em Top white} &  yes, in rising  &  swap & swap & stay${}^{\dagger}$ \\
{\em checker} &  checker's square &  & & \\ \cline{2-5}
{\em in } & yes,  & swap & swap if no blocker  & stay \\ 
{\em critical} &  elsewhere &  &  {\bf or} stay  &  \\ \cline{2-5}
{\em diagonal?} & no & stay & stay & stay \\ \hline
\end{tabular}
\caption{Phase 1 of the white checker moves (see Figure~\ref{durham} for a pictorial
description) \lremind{keywest}} \label{keywest}
\end{table}

\begin{figure}[ht]
  \begin{center}     \setlength{\unitlength}{0.00083333in}
\begingroup\makeatletter\ifx\SetFigFont\undefined%
\gdef\SetFigFont#1#2#3#4#5{%
  \reset@font\fontsize{#1}{#2pt}%
  \fontfamily{#3}\fontseries{#4}\fontshape{#5}%
  \selectfont}%
\fi\endgroup%
{\renewcommand{\dashlinestretch}{30}
\begin{picture}(1362,851)(0,-10)
\put(387,449){\ellipse{74}{74}}
\put(537,599){\ellipse{74}{74}}
\put(87,449){\blacken\ellipse{74}{74}}
\put(87,449){\ellipse{74}{74}}
\put(387,599){\blacken\ellipse{74}{74}}
\put(387,599){\ellipse{74}{74}}
\put(275,261){\ellipse{74}{74}}
\put(200,336){\blacken\ellipse{74}{74}}
\put(200,336){\ellipse{74}{74}}
\put(537,749){\blacken\ellipse{74}{74}}
\put(537,749){\ellipse{74}{74}}
\path(1287,599)(575,599)
\blacken\path(635.000,614.000)(575.000,599.000)(635.000,584.000)(635.000,614.000)
\path(1287,449)(425,449)
\blacken\path(485.000,464.000)(425.000,449.000)(485.000,434.000)(485.000,464.000)
\path(1287,74)(275,74)(275,224)
\blacken\path(290.000,164.000)(275.000,224.000)(260.000,164.000)(290.000,164.000)
\path(12,824)(612,824)(612,224)
	(12,224)(12,824)
\path(162,824)(162,224)
\path(312,824)(312,224)
\path(462,824)(462,224)
\path(12,674)(612,674)
\path(12,524)(612,524)
\path(12,374)(612,374)
\dashline{60.000}(275,261)(537,261)(537,599)
	(275,599)(275,261)
\put(1362,561){\makebox(0,0)[lb]{\smash{{{\SetFigFont{8}{9.6}{\rmdefault}{\mddefault}{\updefault}white checker in critical row}}}}}
\put(1362,411){\makebox(0,0)[lb]{\smash{{{\SetFigFont{8}{9.6}{\rmdefault}{\mddefault}{\updefault}blocker}}}}}
\put(1362,36){\makebox(0,0)[lb]{\smash{{{\SetFigFont{8}{9.6}{\rmdefault}{\mddefault}{\updefault}top white checker in critical diagonal}}}}}
\end{picture}
}   \end{center}
\caption{Example of a blocker \lremind{blockerexample}}
\label{blockerexample}
\end{figure}
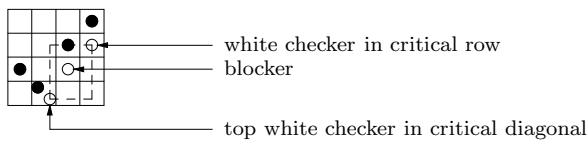

{\em Phase 2}\notation{Phase 2} is a ``clean-up phase'':  if any white checkers
are not {happy}, then move them by sliding them either left
or up so that they become happy.  This is always possible, in a unique way.

% Thanks Allen!

If $\cb$ corresponds to the configuration of
black and white checkers before the move, we denote the one or two
possible configurations after the move by $\cbn$ 
and/or $\cbs$.\notation{$\cn, \cs$}
Examples of the ten cases are given in Figure~\ref{durham}.

\begin{figure}[ht]
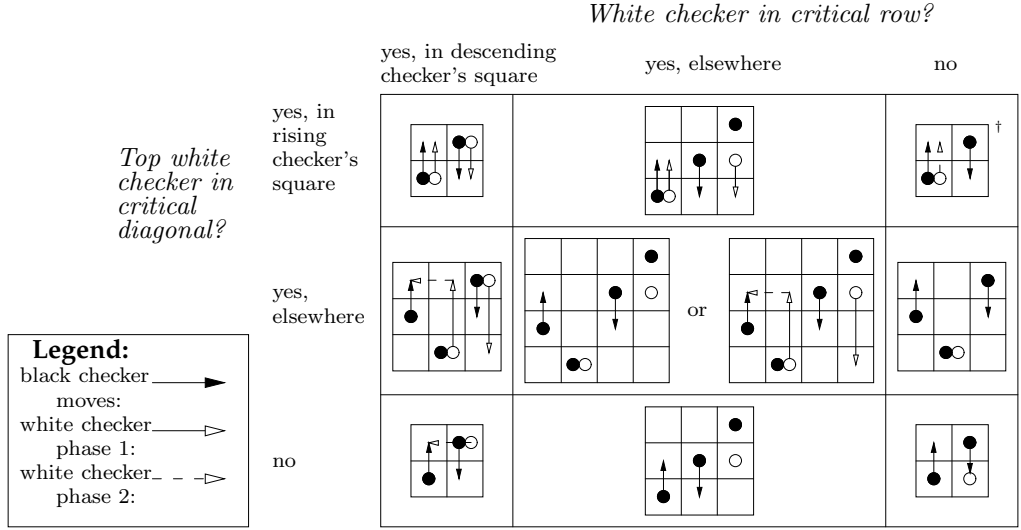

\begin{center}\include{durham}\end{center}
\caption{Examples of the nine cases 
(case $\dagger$ is discussed in Sect.~\ref{bijst})
\lremind{durham}}
\label{durham}
\end{figure}

At the end of the checker game, the black checkers are in position
$\bullet_{\final}$, and as the white checkers are happy, they must lie
on a subset of the black checkers.  The corresponding output is again
a subset $S$ of $\{ 1, \dots, n \}$ of size $k$.    Denote
this final configuration of white checkers by $\circ_S$.\notation{$\circ_S$}

Hence to any two subsets $A, B$ of $\{ 1, \dots, n \}$ of size $k$,
we can associate a formal sum of subsets of size $k$, corresponding
to checker games.
Figure~\ref{g24} illustrates the two checker games starting
with $A=B=\{ 2,4 \}$ for $k=2$, $n=4$.  The output consists of
the two subsets $\{ 1,4 \}$ and $\{ 2,3 \}$.  
Figure~\ref{mult2} illustrates all checker games starting with
$A=B= \{2, 4, 6 \}$ for $k=3$, $n=6$.  The output
is $\{ 2,3,4 \} + 2 \{ 1,3,5 \} + \{ 1, 2, 6 \}$.

\begin{figure}[ht]
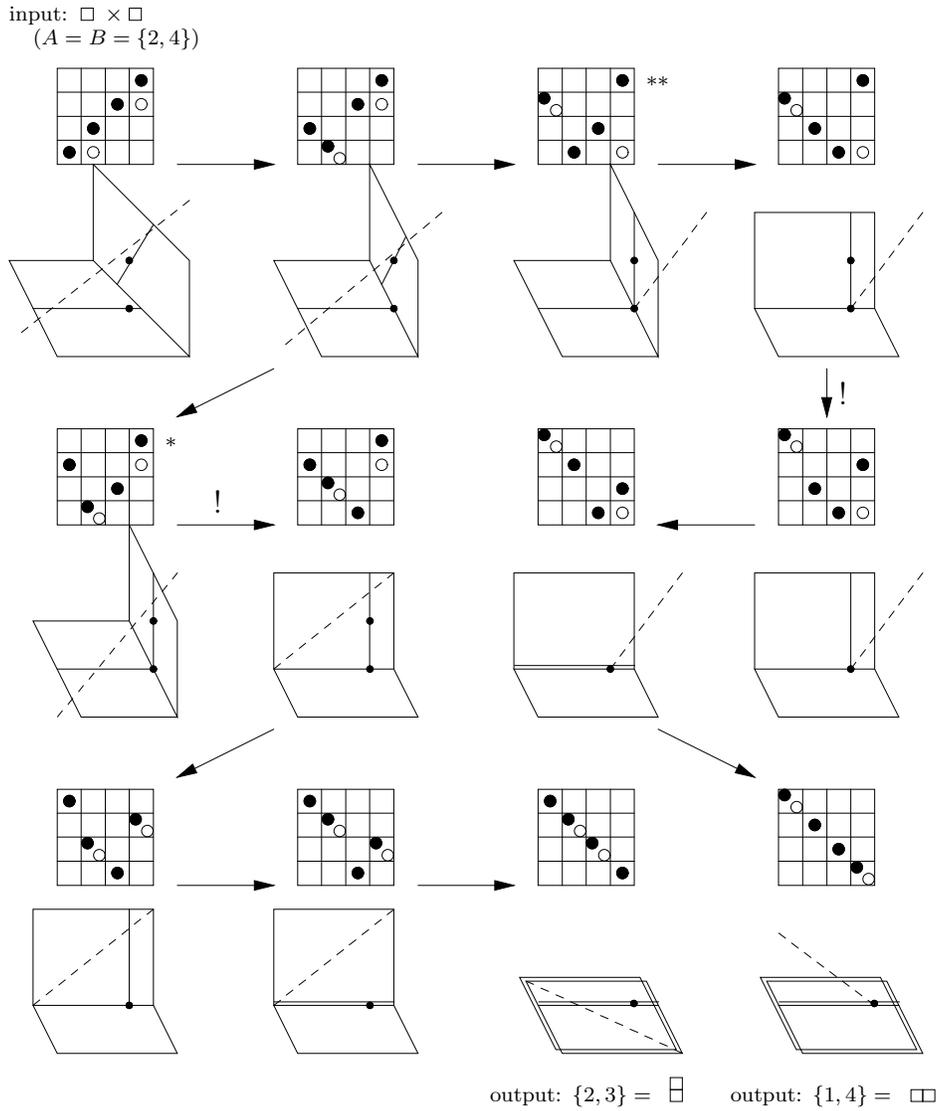

\begin{center} \include{g24} \end{center}
\caption{Two checker games with the same starting position,
and  the geometric
interpretation of Section~\ref{checkergeometry}  (compare to Figure~\ref{butte};
checker configurations * and ** are discussed in Sect.~\ref{cap},
and the moves labeled ``!'' are discussed in Sect.~\ref{bijst})
\lremind{g24}} \label{g24}
\end{figure}

\begin{figure}[ht]
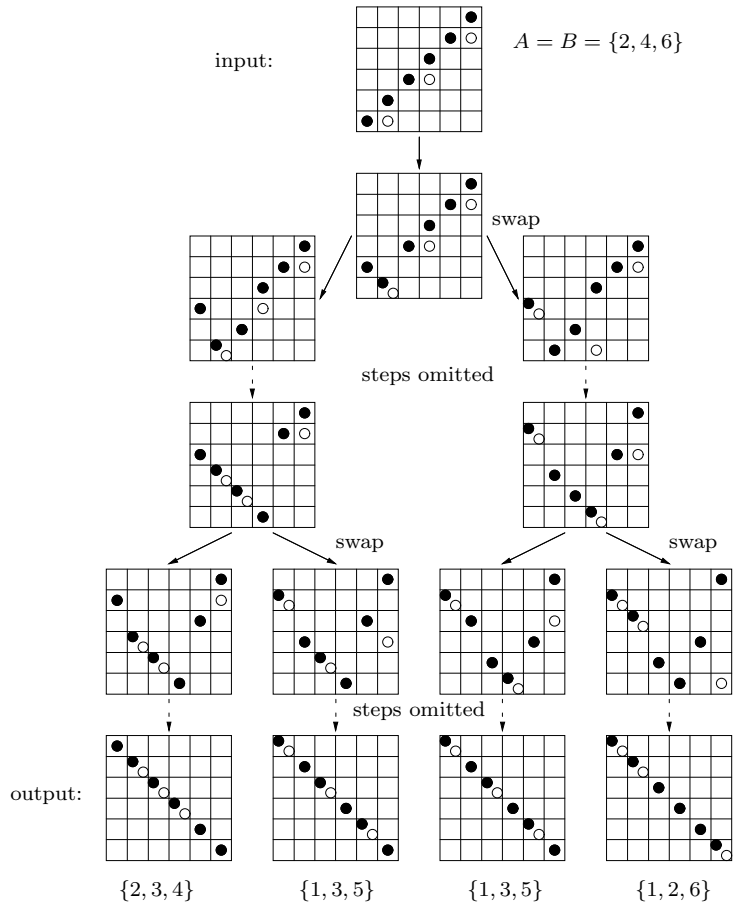

\begin{center} \include{mult2} \end{center}
\caption{Computing $c^{(3,2,1)}_{(2,1),(2,1)} =2$
using $k=3$, $n=6$; some intermediate steps are omitted \lremind{mult2}}
\label{mult2}
\end{figure}

\bpoint{Littlewood-Richardson coefficients in terms of checker games}

Littlewood-Richardson rules give combinatorial descriptions
of Littlewood-Richardson coefficients $c_{\al \be}^{\ga}$, defined
by
$$
s_\al s_\be = \sum_{\ga} c_{\al \be}^{\ga} s_{\ga}
$$
where $\al$ and $\be$ are partitions, and $s_{\bu}$
are Schur functions.
The checker games will compute
compute the sum on the right modulo the ideal 
(in the ring of symmetric functions) defining the cohomology
ring of the Grassmannian $G(k,n)$.

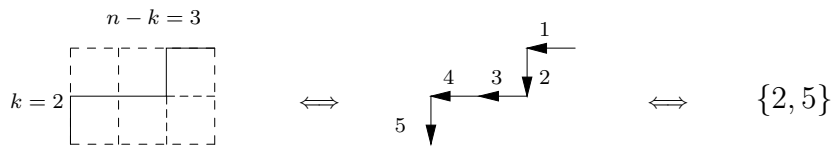
\begin{figure}[ht]
\begin{center} \setlength{\unitlength}{0.00083333in}
\begingroup\makeatletter\ifx\SetFigFont\undefined%
\gdef\SetFigFont#1#2#3#4#5{%
  \reset@font\fontsize{#1}{#2pt}%
  \fontfamily{#3}\fontseries{#4}\fontshape{#5}%
  \selectfont}%
\fi\endgroup%
{\renewcommand{\dashlinestretch}{30}
\begin{picture}(4650,861)(0,-10)
\dashline{60.000}(375,612)(1275,612)(1275,12)
	(375,12)(375,612)
\dashline{60.000}(975,12)(975,312)(1275,312)
\dashline{60.000}(675,612)(675,12)
\path(375,12)(375,312)(975,312)
	(975,612)(1275,612)
\blacken\path(2745.000,342.000)(2625.000,312.000)(2745.000,282.000)(2745.000,342.000)
\path(2625,312)(2925,312)
\blacken\path(3045.000,342.000)(2925.000,312.000)(3045.000,282.000)(3045.000,342.000)
\path(2925,312)(3225,312)
\blacken\path(3195.000,432.000)(3225.000,312.000)(3255.000,432.000)(3195.000,432.000)
\path(3225,312)(3225,612)
\blacken\path(3345.000,642.000)(3225.000,612.000)(3345.000,582.000)(3345.000,642.000)
\path(3225,612)(3525,612)
\blacken\path(2595.000,132.000)(2625.000,12.000)(2655.000,132.000)(2595.000,132.000)
\path(2625,12)(2625,312)
\put(0,237){\makebox(0,0)[lb]{\smash{{{\SetFigFont{8}{9.6}{\rmdefault}{\mddefault}{\updefault}$k=2$}}}}}
\put(600,762){\makebox(0,0)[lb]{\smash{{{\SetFigFont{8}{9.6}{\rmdefault}{\mddefault}{\updefault}$n-k=3$}}}}}
\put(1800,237){\makebox(0,0)[lb]{\smash{{{\SetFigFont{8}{9.6}{\rmdefault}{\mddefault}{\updefault}$\Longleftrightarrow$}}}}}
\put(4650,237){\makebox(0,0)[lb]{\smash{{{\SetFigFont{12}{14.4}{\rmdefault}{\mddefault}{\updefault}$\{ 2,5 \}$}}}}}
\put(3975,237){\makebox(0,0)[lb]{\smash{{{\SetFigFont{8}{9.6}{\rmdefault}{\mddefault}{\updefault}$\Longleftrightarrow$}}}}}
\put(2400,87){\makebox(0,0)[lb]{\smash{{{\SetFigFont{8}{9.6}{\rmdefault}{\mddefault}{\updefault}$5$}}}}}
\put(2700,387){\makebox(0,0)[lb]{\smash{{{\SetFigFont{8}{9.6}{\rmdefault}{\mddefault}{\updefault}$4$}}}}}
\put(3000,387){\makebox(0,0)[lb]{\smash{{{\SetFigFont{8}{9.6}{\rmdefault}{\mddefault}{\updefault}$3$}}}}}
\put(3300,387){\makebox(0,0)[lb]{\smash{{{\SetFigFont{8}{9.6}{\rmdefault}{\mddefault}{\updefault}$2$}}}}}
\put(3300,687){\makebox(0,0)[lb]{\smash{{{\SetFigFont{8}{9.6}{\rmdefault}{\mddefault}{\updefault}$1$}}}}}
\end{picture}
} \end{center}
\caption{The bijection between $Rec_{k,n-k}$ and size $k$
subsets of $\{ 1, \dots, n \}$ \lremind{rec}
} \label{rec}
\end{figure}

More precisely, let $Rec_{k,n-k}$\notation{$Rec_{k,n-k}$} be the set of subdiagrams 
of the rectangle with $k$
rows and $n-k$ columns.  Throughout this article, we identify $Rec_{k,n-k}$ with subsets of
$\{1, \dots, n\}$ of size $k$ using the well-known bijection
(see Figure~\ref{rec}).
Fix $\al$ and $\be$.

\tpoint{Theorem (\LR, first version)} 
\label{clr} \lremind{clr}
{\em
\begin{enumerate}
\item[(a)] 
${\displaystyle \sum_{ \text{games $G$}} s_{\text{output($G$)}}
= \sum_{\ga \in Rec_{k,n-k}} c_{\al \be}^{\ga} s_{\ga}}$,
where the left sum is over all checker games with input $\al$ and $\be$,
and $\text{output($G$)}$ is the output of checker game $G$.
\item[(b)]
Hence if $\ga \in Rec_{k,n-k}$, then
the integer $c_{\al \be}^{\ga}$ is the number
of checker games starting with configuration $\circ_{\al, \be}\bullet_{\init}$
and ending with configuration $\circ_{\ga}\bullet_{\final}$.
\end{enumerate}
}

For example, Figure~\ref{g24} computes $s_{(1)}^2  = s_{(2)} + s_{(1,1)}$.
Figure~\ref{mult2} computes $c^{(3,2,1)}_{(2,1),(2,1)} =2$
using $k=3$, $n=6$.

Theorem~\ref{clr} follows immediately from Theorem~\ref{glr2} (\LR, second version).

\epoint{Bijection to tableaux} 
\label{bijst} \lremind{bijst}
\secretnote{Anders e-mail May 21 '02; placing $\gamma-\alpha$ in
  $\beta$.}  There is a straightforward bijection to tableaux (using
the tableaux description of \cite[Cor.~5.1.2]{fulton}).  Whenever
there is move described by a $\dagger$ in Figure~\ref{durham} (see
also Table~\ref{keywest}), where the ``rising'' white checker is the
$r^{\text{th}}$ white checker (counting by row) and the
$c^{\text{th}}$ (counting by column), place an $r$ in row $c$ of the
tableau.  

For example, in Figure~\ref{g24}, the left-most output corresponds to
the (one-cell) tableau ``2'', and the right-most output corresponds to
the tableau ``1''.  The moves where the tableaux are filled are marked
with ``!''.

\epoint{Remarks}  
{\em (a)} Like Pieri's formula and Monk's formula, this rule most naturally
gives all terms in an intersection all at once (part (a)),
but the individual coefficients can be easily extracted (part (b)).

\noindent
{\em (b)} A derivation of Pieri's formula from the
\LR{}
is left as an exercise to
the reader.  Note that Pieri's original proof 
was also by degeneration methods.

\noindent
{\em (c)} Some properties of Littlewood-Richardson coefficients
clearly follow from the \LR, while others do not.  For example, it is
not clear why $c_{\al \be}^{\ga} = c_{\be \al}^{\ga}$.  However, it
can be combinatorially shown (most easily via the link to puzzles,
Sect.~\ref{checkerpuzzle}) that (i) the rule is independent of the
choice of $n$ and $k$ (i.e. the computation of $c_{\al \be}^{\ga}$ is
independent of any $n$ and $k$ such that $\ga \in Rec_{k,n-k}$), and
(ii) the ``triality'' $c_{\al \be}^{\ga} = c_{\be
  \ga^{\vee}}^{\al^{\vee}}$ for $\al, \be, \ga \in Rec_{k,n-k}$ holds.

\noindent
{\em (d)}
We will need the following  combinatorial observation.  

\tpoint{Lemma} \lremind{midsort} \label{midsort} {\em Suppose at some
  point in the algorithm, the descending  checker  is in column $c$.
 Suppose the white checkers are at $(r_1, c_1)$, \dots, $(r_k, c_k)$
 with $c_1 < \dots < c_k$.  Then $(r_1, \dots, r_k)$ is
increasing for $c_i \leq c$ and decreasing for $c_i \geq c$.}

This follows from a
straightforward induction showing that this property is preserved 
by each move.  
We say that configuration $\cb$ with this property is {\em
  mid-sort}.\notation{mid-sort}  For example, the white checkers of
Figures~\ref{mult2} and~\ref{dc} are mid-sort.

\section{Describing the geometry of the Grassmannian and flag variety
with checkers} \label{checkergeometry}
\lremind{checkergeometry}
In this section, we interpret 
checker configurations geometrically, and state the main technical
result of the paper, Theorem~\ref{glr3} (the \LR, final version).
{Black checkers} will correspond to the relative position of
the two flags $\fm_{\bu}$ and $\ff_{\bu}$, and {white checkers} will
correspond to the  position of the $k$-plane relative to
the flags.

\bpoint{The relative position of two flags, given by black checkers;
the variety $X_{\bullet}$} \label{domdefhere}\lremind{domdefhere}
Given two flags $\fm_{\bu}$ and $\ff_{\bu}$,\notation{$\fm_{\bu}, \ff_{\bu}$} construct an $n \times n$ ``rank'' table of
the numbers $\dim \fm_i \cap \ff_j$.  Up to the action of $GL(n)$, the
two flags are specified by this table.

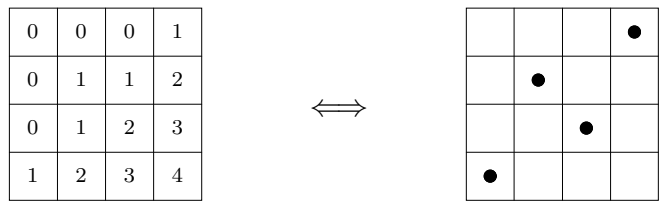
\begin{figure}[ht]
\begin{center} \setlength{\unitlength}{0.00083333in}
\begingroup\makeatletter\ifx\SetFigFont\undefined%
\gdef\SetFigFont#1#2#3#4#5{%
  \reset@font\fontsize{#1}{#2pt}%
  \fontfamily{#3}\fontseries{#4}\fontshape{#5}%
  \selectfont}%
\fi\endgroup%
{\renewcommand{\dashlinestretch}{30}
\begin{picture}(4074,1239)(0,-10)
\path(12,1212)(1212,1212)(1212,12)
	(12,12)(12,1212)
\path(312,1212)(312,12)
\path(612,1212)(612,12)
\path(912,1212)(912,12)
\path(12,912)(1212,912)
\path(12,612)(1212,612)
\path(12,312)(1212,312)
\path(2862,1212)(4062,1212)(4062,12)
	(2862,12)(2862,1212)
\path(3162,1212)(3162,12)
\path(3462,1212)(3462,12)
\path(3762,1212)(3762,12)
\path(2862,912)(4062,912)
\path(2862,612)(4062,612)
\path(2862,312)(4062,312)
\put(3912,1062){\blacken\ellipse{76}{76}}
\put(3912,1062){\ellipse{76}{76}}
\put(3012,162){\blacken\ellipse{76}{76}}
\put(3012,162){\ellipse{76}{76}}
\put(3612,462){\blacken\ellipse{76}{76}}
\put(3612,462){\ellipse{76}{76}}
\put(3312,762){\blacken\ellipse{76}{76}}
\put(3312,762){\ellipse{76}{76}}
\put(1024,424){\makebox(0,0)[lb]{\smash{{{\SetFigFont{8}{9.6}{\rmdefault}{\mddefault}{\updefault}$3$}}}}}
\put(1024,124){\makebox(0,0)[lb]{\smash{{{\SetFigFont{8}{9.6}{\rmdefault}{\mddefault}{\updefault}$4$}}}}}
\put(424,124){\makebox(0,0)[lb]{\smash{{{\SetFigFont{8}{9.6}{\rmdefault}{\mddefault}{\updefault}$2$}}}}}
\put(1024,724){\makebox(0,0)[lb]{\smash{{{\SetFigFont{8}{9.6}{\rmdefault}{\mddefault}{\updefault}$2$}}}}}
\put(1024,1024){\makebox(0,0)[lb]{\smash{{{\SetFigFont{8}{9.6}{\rmdefault}{\mddefault}{\updefault}$1$}}}}}
\put(724,724){\makebox(0,0)[lb]{\smash{{{\SetFigFont{8}{9.6}{\rmdefault}{\mddefault}{\updefault}$1$}}}}}
\put(424,724){\makebox(0,0)[lb]{\smash{{{\SetFigFont{8}{9.6}{\rmdefault}{\mddefault}{\updefault}$1$}}}}}
\put(1887,537){\makebox(0,0)[lb]{\smash{{{\SetFigFont{12}{14.4}{\rmdefault}{\mddefault}{\updefault}$\Longleftrightarrow$}}}}}
\put(124,124){\makebox(0,0)[lb]{\smash{{{\SetFigFont{8}{9.6}{\rmdefault}{\mddefault}{\updefault}$1$}}}}}
\put(724,124){\makebox(0,0)[lb]{\smash{{{\SetFigFont{8}{9.6}{\rmdefault}{\mddefault}{\updefault}$3$}}}}}
\put(124,424){\makebox(0,0)[lb]{\smash{{{\SetFigFont{8}{9.6}{\rmdefault}{\mddefault}{\updefault}$0$}}}}}
\put(724,424){\makebox(0,0)[lb]{\smash{{{\SetFigFont{8}{9.6}{\rmdefault}{\mddefault}{\updefault}$2$}}}}}
\put(424,424){\makebox(0,0)[lb]{\smash{{{\SetFigFont{8}{9.6}{\rmdefault}{\mddefault}{\updefault}$1$}}}}}
\put(124,724){\makebox(0,0)[lb]{\smash{{{\SetFigFont{8}{9.6}{\rmdefault}{\mddefault}{\updefault}$0$}}}}}
\put(724,1024){\makebox(0,0)[lb]{\smash{{{\SetFigFont{8}{9.6}{\rmdefault}{\mddefault}{\updefault}$0$}}}}}
\put(424,1024){\makebox(0,0)[lb]{\smash{{{\SetFigFont{8}{9.6}{\rmdefault}{\mddefault}{\updefault}$0$}}}}}
\put(124,1024){\makebox(0,0)[lb]{\smash{{{\SetFigFont{8}{9.6}{\rmdefault}{\mddefault}{\updefault}$0$}}}}}
\end{picture}
} \end{center}
%\begin{tabular}{|c|c|c|c|} \hline
%0 & 0 & 0 & 1 \\ \hline
%0 & 1 & 1 & 2 \\ \hline
%0 & 1 & 2 & 3 \\ \hline
%1 & 2 & 3 & 4 \\ \hline
%\end{tabular}
%$\quad \quad \Longleftrightarrow \quad \quad$
%\begin{tabular}{|c|c|c|c|} \hline
% &  &  & $\bullet$ \\ \hline
% & $\bullet$ &  &  \\ \hline
% &  & $\bullet$ &  \\ \hline
%$\bullet$ &  &  &  \\ \hline
%\end{tabular}
\caption{The relative position of two flags, given by numbers
  on an $n \times n$ board, and by a configuration of black checkers
  \lremind{idaho}}
\label{idaho}
\end{figure}

This data is equivalent to the data of $n$ black checkers on the
checkerboard such that no two are in the same row or column.  The
bijection is given as follows.  
We say a square $(i_1,j_1)$ {\em dominates}\notation{dominate} another square
$(i_2, j_2)$ if $i_1 \geq i_2$ and $j_1 \geq j_2$.
Given the black checkers, $\dim \fm_i
\cap \ff_j$ is given by the number of black checkers dominated by
square $(i,j)$.  An example of the bijection is given in Figure~\ref{idaho}.

Note that each square in the table corresponds to a vector
space, whose dimension is the number of black checkers dominated by
that square.  The vector space is the span of the vector
spaces corresponding to the black checkers it dominates.

Let $X_{\bullet}$ be the (locally closed) subvariety of $Fl(n) \times
Fl(n)$ corresponding to flags $(\fm_{\bu}, \ff_{\bu})$ in relative
position given by black checker configuration
$\bullet$.\notation{$X_{\bullet}$}  The variety
$X_{\bullet}$ is smooth, and its codimension in $Fl(n)
\times Fl(n)$ is the number of pairs of distinct black checkers $a$ and
$b$ such that $a$ dominates $b$.  (This is a straightforward exercise;
it also follows quickly from Sect.~\ref{bsv}.)
This sort of construction is common in the literature.
%\remind{Later look at Fulton's paper in Duke (his first?) for the
%same idea.} 

\epoint{Side remark: Schubert varieties of the flag variety} For a
fixed $\ff_{\bu}$ and fixed configuration $\bullet$ of black checkers,
the set of $\fm_{\bu}$ in $Fl(n)$ such that the relative position of
$\fm_{\bu}$ and $\ff_{\bu}$ is given by $\bullet$ is an open Schubert
cell $\Sv^{\ci}_{\bullet}$ (i.e. the second projection $pr_2:
X_{\bullet} \rightarrow Fl(n)$ is a fibration by Schubert cells).
This set is a $GL_n$-orbit (what some authors call a ``double Schubert
cell'').  Schubert cells are usually indexed by permutations; the
bijection between checker configurations and permutations was given in
Section~\ref{seattle}.  For example, the permutation corresponding to
Figure~\ref{idaho} is $1324$.

\epoint{The specialization order} 
\label{madison} \lremind{madison}
Given a point $p$ of $Fl(n)$ (parametrizing $\fm_{\bu}$) in the dense
open Schubert cell (with respect to a fixed reference flag
$\ff_{\bu}$), the {specialization order} (Sect.~\ref{seattle}) can
be interpreted as a sequence in $Fl(n)$, consisting of a chain of
$\binom n 2$ $\proj^1$'s, starting at $p$ and ending with the ``most
degenerate'' point of $Fl(n)$ (corresponding to the reference flag
$\ff_{\bu}$).

We first describe the chain informally.  Each $\proj^1$ corresponds to
a move of black checkers.  All but one point of the $\proj^1$ lies in
one open stratum $X_{\bullet}$; the remaining point (where the
$\proj^1$ meets the next component of the one-parameter degeneration)
lies on an open stratum $X_{\bn}$ of dimension one lower.  If the move
corresponds to the descending checker in row $r$ dropping one row (and
another checker to the left rising one row), then all components of
the flags $\ff_{\bu}$ and $\fm_{\bu}$ except $\fm_r$ are held fixed.
For example, the geometry corresponding to the specialization order
for $n=4$ is shown in Figure~\ref{butte}.

%\begin{figure}
%remind{Show the six moves of the black checkers for $n=4$,
%along with the corresponding permutations.}
%\caption{Geometric interpretation of the {\em specialization order} for $n=4$
%(see also Figure~\ref{butte})
%\lremind{butte2} remind{clean this up; figure sketched}}
%\label{butte2}
%\end{figure} 

More precisely, there is a $\proj^1$-fibration $X_{\bullet} \cup
X_{\bn} \rightarrow X_{\bn}$; the $\proj^1$ corresponding to
a point of $X_{\bullet}$ described in the previous paragraph is a
fiber.  A useful alternate description of $X_{\bullet} \cup X_{\bn}$ is
given in Section~\ref{portland}.

\bpoint{The position of a $k$-plane relative to two given flags,
in terms of white checkers; two-flag Schubert varieties $Y_{\cb}$
and $X_{\cb}$} \label{eugene} \lremind{eugene}
Suppose two flags $\fm_{\bu}$ and $\ff_{\bu}$ are in relative
position given by black checker configuration $\bullet$.  The position
of a $k$-plane $V$ relative to the two flags, up to the action of the
subgroup of $GL(n)$ fixing the two flags, is determined by the table
of numbers $\dim V \cap \fm_i \cap \ff_j$.

This data is equivalent to the data of $k$ white checkers on the
checkerboard that are {happy} (see Sect.~\ref{rulewhite}), with no two 
in the same row or column.
The bijection is given as follows: 
$\dim V \cap
\fm_i \cap \ff_j$ is 
the number of white checkers in
squares dominated by $(i,j)$.  See Figure~\ref{g24} for example; the
line $\proj V$ is depicted as a dashed line.

If $\fm_{\bu}$ and $\ff_{\bu}$ are two flags whose
relative position is given by $\bullet$, 
let the {\em open two-flag Schubert variety} $Y_{\cb} \subset G(k,n)$ be 
the set of $k$-planes
whose position relative to the flags is given
by $\cb$;  
define the {\em closed two-flag Schubert variety} $\overline{Y}_{\cb}$ to be $\Cl_{G(k,n)} Y_{\cb}$.
Let $X_{\cb}$ (resp. $\overline{X}_{\cb}$) in $G(k,n) \times X_{\bullet}$
be the {\em universal open (resp. closed) two-flag Schubert 
variety}.\notation{open two-flag
Schubert variety $Y_{\cb}$,
closed two-flag
Schubert variety $\overline{Y}_{\cb}$,
universal two-flag Schubert variety,
$X_{\cb}$, $\overline{X}_{\cb}$}
%\remind{Ask Sara Billey about this.} 
(Similar constructions are common in the literature.)

Note that
(i) $X_{\cb} \rightarrow X_{\bullet}$ is
a $Y_{\cb}$-fibration; (ii) $\overline{X}_{\cb}
\rightarrow X_{\bullet}$ is an $\overline{Y}_{\cb}$-fibration, and proper; 
(iii)
$G(k,n)$ is the disjoint union of the $Y_{\cb}$ (for fixed
$\bullet$ and $\ff_{\bu}$, $\fm_{\bu}$); (iv) $G(k,n) \times Fl(n) \times Fl(n)$ is the
disjoint union of the $X_{\cb}$.  {\em Caution:}
the disjoint unions of (iii) and (iv) are not in general stratifications;
Section~\ref{cap} (a) provides a counterexample.

\tpoint{Lemma} {\em The variety $Y_{\cb}$ is irreducible and smooth; its dimension is the sum over all white checkers
$w$ of the number of black checkers $w$ dominates minus
the number of white checkers $w$ dominates (including itself).}
\label{dimY} \lremind{dimY}

The proof is straightforward and hence omitted.

\tpoint{Proposition (Initial position of white checkers)} 
\label{ip} \lremind{ip}
{\em Suppose $A=\{ a_1, \dots, a_k \}$ and $B = \{ b_1, \dots, b_k \}$ are two
subsets of $\{ 1, \dots , n \}$ of size $k$, and $\fm_{\bu}$ and $\ff_{\bu}$ are two
transverse (opposite) flags (i.e. with relative position given by
$\bullet_{\init}$).   Then $\Sv_A(\fm_{\bu}) \cap
\Sv_B(\ff_{\bu})$ is the closed two-flag Schubert variety $\overline{Y}_{\circ_{A,B}\bullet_{\init}}$.}

In the literature, these intersections are known as
{\em Richardson varieties}~\cite{richardson}; see
\cite{lakshmibai} for more discussion and references.
They were also called {\em skew Schubert varieties}
by Stanley \cite{stanleyskew}.

{\em Proof.}
We deal first with the case of characteristic 0.
By the Kleiman-Bertini theorem \cite{kleimanbertini}, 
$\Sv_A(\fm_{\bu}) \cap
\Sv_B(\ff_{\bu})$ is reduced of the expected dimension.  
The generic point of any of its components lies in
$Y_{\circ_1 \bullet_{\init} }$ for some configuration $\circ_1$ 
of white checkers, where the first coordinates
of the white checkers of $\circ_1$ are given by the set $A$
and the second coordinates are given by the set $B$.
A short calculation using Lemma~\ref{dimY} yields
$\dim Y_{\circ_1 \bullet_{\init} } \leq
\dim Y_{ \circ_{A,B} \bullet_{\init}}$, 
with equality holding if and only
if $\circ_1 = \circ_{A,B}$.  (Reason:
the  sum over all white checkers
$w$ in $\circ_1$ of the number of black checkers $w$ dominates
is $\sum_{a \in A} a + \sum_{b \in B} b - kn$, which is independent of $\circ_1$, so 
$\dim Y_{ \circ_1 \bullet_{\init}}$ is maximized when
there is no white checker dominating another, which is
the definition of $\circ_{A,B}$.)  
Then it can be checked directly that 
$\dim Y_{ \circ_{A,B}\bullet_{\init}} = \dim \Sv_A \cap \Sv_B$.
As $Y_{ \circ_{A,B}\bullet_{\init}}$
is irreducible, the result in characteristic 0 follows.  

In positive characteristic, the same argument shows that the cycle
$\Sv_A(\fm_{\bu}) \cap \Sv_B(\ff_{\bu})$ is some positive multiple of
the the cycle $\overline{Y}_{\circ_{A,B}\bullet_{\init}}$.  
It is an easy exercise to show that the intersection is transverse,
i.e. that this multiple is 1.  It will be easier still to conclude
the proof combinatorially; we will do this in Section~\ref{wimbledon}.

\bpoint{The \LR: deforming cycles in
  the Grassmannian}  \label{baltimore} \lremind{baltimore}
We can now give a geometric interpretation of
Theorem~\ref{clr} (which will be the first version of the \LR).
We wish to compute the class (in $H_*(G(k,n))$) of the intersection of
two Schubert cycles.  By the Kleiman-Bertini theorem
\cite{kleimanbertini}
(or\secretnote{There is a macro \sigkb which links to the 
generic smoothness theorem in \cite{si}.} 
our Grassmannian Kleiman-Bertini 
theorem \cite[\sigkb]{si} in positive characteristic), 
this is the class of the intersection
of two Schubert varieties with respect to two general (transverse) flags, which by
Proposition~\ref{ip} is  $[\overline{Y}_{\circ_{A,B} \bullet_{\init} }]$.
By a sequence of codimension 1 degenerations 
(corresponding to the {specialization order}), we degenerate
the two flags until they are equal.  The base of each degeneration is
 a $\proj^1$ in $Fl(n)$, where $\A^1 \subset X_{\bullet}$
and $\{ \infty \} \in X_{\bn}$.  
We have a $\overline{Y}_{\cb}$-fibration above $\A^1$.
The fiber
above $\infty$ of the 
closure of this fibration (in $G(k,n) \times \proj^1$) turns
out to be one of $\overline{Y}_{\cbn }$, $\overline{Y}_{\cbs }$, 
or $\overline{Y}_{\cbn } \cup
\overline{Y}_{\cbs}$ (Theorem~\ref{glr2}; $\cn$ and $\cs$ were defined 
in Sect.~\ref{rulewhite}); the components appear with multiplicity 1.
In other words, the two-flag Schubert variety degenerates
to another two-flag Schubert variety, or to the union of two.
Figure~\ref{g24} illustrates the geometry
behind the computation of $s_{(1)}^2 = s_{(1,1)} + s_{(2)}$ in $G(2,4)$.

\point % point needed because I refer back to the previous section
\label{hamilton} \lremind{hamilton}
More precisely, assume the black checker configuration $\bullet$ is in
the specialization order, and the white checkers are {\em mid-sort}
(in the sense of Lemma~\ref{midsort}).
Consider the diagram:\lremind{bloop}
\begin{equation}
\label{bloop}
\begin{array}{ccccc}
\overline{X}_{\cb} := \Cl_{G(k,n) \times X_{\bullet}} X_{\cb}  & \hookrightarrow & 
\Cl_{G(k,n) \times (X_{\bullet} \cup X_{\bn})} X_{\cb}
& \hookleftarrow & D_X \\
\downarrow  & & \downarrow & & \downarrow \\
X_{\bullet} & \hookrightarrow & X_{\bullet} \cup X_{\bn} & \hookleftarrow
& X_{\bn}.
\end{array}
\end{equation}
The Cartier divisor $D_X$ is defined by fibered product.\notation{$D_X$}
Note that the vertical morphisms are proper, 
the vertical morphism on the left is a $\overline{Y}_{\cb}$-fibration, 
the horizontal inclusions
on the left are open immersions, and the horizontal inclusions on the right
are closed immersions. 

The geometric constructions described in Section~\ref{baltimore} are
obtained from \eqref{bloop} by base changing via $\proj^1 \rightarrow
X_{\bullet} \cup X_{\bn}$ (described in Sect.~\ref{madison}), to obtain:
\lremind{bloop2}
\begin{equation}
\label{bloop2}
\begin{array}{ccccc}
\Y_{\cb} & \hookrightarrow & 
\Cl_{G(k,n) \times \proj^1} \Y_{\cb}
& \hookleftarrow & D_{\Y} \\
 \downarrow & & \downarrow & & \downarrow \\
\A^1 & \hookrightarrow & \proj^1 & \hookleftarrow
& \{ \infty \}.
\end{array}
\end{equation}
Here $\Y_{\cb}$,  $\Cl_{G(k,n) \times \proj^1} \Y_{\cb}$, and $D_{\Y}$ are defined
by fibered product (or restriction) from \eqref{bloop}.
Again, the vertical morphisms are proper, 
the vertical morphism on the left is a $\overline{Y}_{\cb}$-fibration, 
the horizontal inclusions
on the left are open immersions, and the horizontal inclusions on the right
are closed immersions.  

The informal statement of Section~\ref{baltimore} can now be
made precise:

\tpoint{Theorem (\LR, second version)}\label{glr2} \lremind{glr2} 
{\em $D_{\Y} = \overline{Y}_{
    \cbn}$, $\overline{Y}_{\cbs}$, or
 $\overline{Y}_{\cbn } \cup \overline{Y}_{\cbs }$.}

By base change from \eqref{bloop} to \eqref{bloop2}, Theorem~\ref{glr2} is a consequence of the following,
which 
is proved in Section~\ref{thmpf}.
(The notation $\Y_{\cb}$ and $D_{\Y}$ will not
be used hereafter.)

\tpoint{Theorem (\LR, final version)}\label{glr3} \lremind{glr3}
 {\em $D_X = \overline{X}_{\cbn}$,
$\overline{X}_{\cbs}$,
  or $\overline{X}_{ \cbn} \cup \overline{X}_{ \cbs}$.}

In other words, a particular divisor $D_X$ (corresponding to 
$X_{\bn} \subset \Xbb$) on a universal two-flag Schubert variety
is another such variety, or the union of two such varieties.

\bpoint{Enumerative problems and checker tournaments}
\lremind{wimbledon}
\label{wimbledon} Suppose $\Sv_{\al_1}$, \dots, $\Sv_{\al_{\ell}}$ are
Schubert classes on $G(k,n)$ of total codimension $\dim G(k,n)$.  Then
(the degree of) their intersection --- the solution to an enumerative
problem by the 
Kleiman-Bertini theorem \cite{kleimanbertini} (or 
our Grassmannian Kleiman-Bertini theorem \cite[\sigkb]{si} 
in positive characteristic)
--- can
clearly be inductively computed using the \LR.  Hence Schubert
problems can be solved by counting {\em checker tournaments} of
$\ell-1$ games, where the input to the first game is $\al_1$ and
$\al_2$, and for $i>1$ the input to the $i^{\text{th}}$ game is
$\al_{i+1}$ and the output of the previous game.\notation{checker
  tournaments} (The outcome of each checker tournament will always be
the same --- the class of a point.)

\noindent 
{\em Conclusion of proof of Proposition~\ref{ip} in positive characteristic.}
We will show that the multiplicity with which 
$\overline{Y}_{\circ_{A,B}\bullet_{\init}}$ appears
in 
$\Sv_A(\fm_{\bu}) \cap \Sv_B(\ff_{\bu})$ is 1.
(We will not use the Grassmannian Kleiman-Bertini Theorem
\cite[\sigkb]{si} as its proof relies on Prop.~\ref{ip}.)

Choose $C = \{c_1, \dots, c_k \}$ such that $\dim \Sv_A \cup \Sv_B
\cup \Sv_C = 0$ (where $\cup$ is the cup product in cohomology) and
$\deg \Sv_A \cup \Sv_B \cup \Sv_C > 0$.  In characteristic 0, the
above discussion shows that $\deg \Sv_A \cup \Sv_B \cup \Sv_C$ is the
number of checker tournaments with inputs $A$, $B$, $C$.  In positive
characteristic, the above discussion shows that if the multiplicity is
greater than one, then $\deg \Sv_A \cup \Sv_B \cup \Sv_C$ is {\em
  strictly less} than the same number of checker games.  But $\deg
\Sv_A \cup \Sv_B \cup \Sv_C$ is independent of characteristic,
yielding a contradiction. \epf

\bpoint{Cautions} \label{cap} \lremind{cap} (a) The degenerations used
in the \LR{} follow the specialization order.  An arbitrary path through
the Bruhat order will not work in general.  For example, if $\cb$ is
as shown on the left of Figure~\ref{capfig}, then $X_{\cb}$
corresponds to points $p_1$ and $p_2$ in $\proj^3$, lines $\ell_1$ and
$\ell_2$ through $p_1$ such that $\ell_1$, $\ell_2$, and $p_2$ span
$\proj^3$, and a point $q \in \ell_1 - p_1$.  Then for example $\fm_3
= \spam( \ell_2,p_2)$ and the line corresponding to the point of
$G(2,4)$ is $\spam( q, p_2)$.  The degeneration shown in
Figure~\ref{capfig} (to $\bullet'$, say) corresponds to letting $p_2$
tend to $p_1$, and remembering the line $\ell_3$ of approach.  Then
 the divisor on $\Cl_{G(k,n) \times \left( \Xbb \right)} X_{\cb}$ 
corresponding to $X_{\bullet'}$ 
parametrizes lines through $p_1$ contained in
$\spam (\ell_1,\ell_3)$; this is not of the form $X_{\circ' \bullet'}$
for any $\circ'$.

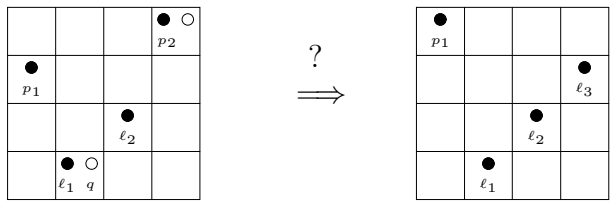
\begin{figure}[ht]
\begin{center} \setlength{\unitlength}{0.00083333in}
\begingroup\makeatletter\ifx\SetFigFont\undefined%
\gdef\SetFigFont#1#2#3#4#5{%
  \reset@font\fontsize{#1}{#2pt}%
  \fontfamily{#3}\fontseries{#4}\fontshape{#5}%
  \selectfont}%
\fi\endgroup%
{\renewcommand{\dashlinestretch}{30}
\begin{picture}(3774,1239)(0,-10)
\put(987,1137){\blacken\ellipse{74}{74}}
\put(987,1137){\ellipse{74}{74}}
\put(1137,1137){\ellipse{74}{74}}
\put(387,237){\blacken\ellipse{74}{74}}
\put(387,237){\ellipse{74}{74}}
\put(537,237){\ellipse{74}{74}}
\put(3312,537){\blacken\ellipse{74}{74}}
\put(3312,537){\ellipse{74}{74}}
\put(3612,837){\blacken\ellipse{74}{74}}
\put(3612,837){\ellipse{74}{74}}
\put(2712,1137){\blacken\ellipse{74}{74}}
\put(2712,1137){\ellipse{74}{74}}
\put(3012,237){\blacken\ellipse{74}{74}}
\put(3012,237){\ellipse{74}{74}}
\put(162,837){\blacken\ellipse{74}{74}}
\put(162,837){\ellipse{74}{74}}
\put(762,537){\blacken\ellipse{74}{74}}
\put(762,537){\ellipse{74}{74}}
\path(612,1212)(612,12)
\path(912,1212)(912,12)
\path(12,912)(1212,912)
\path(3162,1212)(3162,12)
\path(2562,312)(3762,312)
\path(2562,612)(3762,612)
\path(2562,912)(3762,912)
\path(2862,1212)(2862,12)
\path(3462,1212)(3462,12)
\path(2562,1212)(3762,1212)(3762,12)
	(2562,12)(2562,1212)
\path(12,612)(1212,612)
\path(312,1212)(312,12)
\path(12,312)(1212,312)
\path(12,1212)(1212,1212)(1212,12)
	(12,12)(12,1212)
\put(947,987){\makebox(0,0)[lb]{\smash{{{\SetFigFont{5}{6.0}{\rmdefault}{\mddefault}{\updefault}$p_2$}}}}}
\put(102,687){\makebox(0,0)[lb]{\smash{{{\SetFigFont{5}{6.0}{\rmdefault}{\mddefault}{\updefault}$p_1$}}}}}
\put(2652,987){\makebox(0,0)[lb]{\smash{{{\SetFigFont{5}{6.0}{\rmdefault}{\mddefault}{\updefault}$p_1$}}}}}
\put(502,87){\makebox(0,0)[lb]{\smash{{{\SetFigFont{5}{6.0}{\rmdefault}{\mddefault}{\updefault}$q$}}}}}
\put(707,387){\makebox(0,0)[lb]{\smash{{{\SetFigFont{5}{6.0}{\rmdefault}{\mddefault}{\updefault}$\ell_2$}}}}}
\put(3557,687){\makebox(0,0)[lb]{\smash{{{\SetFigFont{5}{6.0}{\rmdefault}{\mddefault}{\updefault}$\ell_3$}}}}}
\put(3257,387){\makebox(0,0)[lb]{\smash{{{\SetFigFont{5}{6.0}{\rmdefault}{\mddefault}{\updefault}$\ell_2$}}}}}
\put(2957,87){\makebox(0,0)[lb]{\smash{{{\SetFigFont{5}{6.0}{\rmdefault}{\mddefault}{\updefault}$\ell_1$}}}}}
\put(327,87){\makebox(0,0)[lb]{\smash{{{\SetFigFont{5}{6.0}{\rmdefault}{\mddefault}{\updefault}$\ell_1$}}}}}
\put(1812,612){\makebox(0,0)[lb]{\smash{{{\SetFigFont{12}{14.4}{\rmdefault}{\mddefault}{\updefault}$\Longrightarrow$}}}}}
\put(1887,837){\makebox(0,0)[lb]{\smash{{{\SetFigFont{12}{14.4}{\rmdefault}{\mddefault}{\updefault}?}}}}}
\end{picture}
} \end{center}
\caption{The dangers of straying from the specialization order
  \lremind{capfig}}
\label{capfig}
\end{figure}

(b)  Unlike the variety $\overline{X}_{\bullet} = \Cl_{Fl(n)\times Fl(n)} X_{\bullet}$,
the variety $\overline{X}_{\cb}$ cannot be  defined numerically, i.e. in general
$\overline{X}_{\cb}$ will be only one irreducible component of
$$X'_{\circ \bullet} := \{ (V, \ff_{\bu}, \fm_{\bu}) \in G(k,n) \times X_{\bullet} \subset G(k,n) \times Fl(n) \times Fl(n) : \dim V \cap \ff_i \cap \fm_j \geq \gamma^{i,j}_{\circ} \}$$
where $\ga^{i,j}_{\circ}$ is the number of white checkers dominated by $(i,j)$.  For example, in Figure~\ref{g24}, if $\cb$ is the configuration marked ``*'' and $\circ' \bullet$ is the configuration marked ``**'', then $X'_{\cb} =
\overline{X}_{\cb} \cup
\overline{X}_{\circ' \bullet}$.

\section{Application:  Littlewood-Richardson rules}
\label{app2}
\lremind{app2} In this section, we discuss the bijection between
checkers, the classical Littlewood-Richardson rule involving tableaux,
and puzzles.  We extend the checker and puzzle rules to
$K$-theory, proving a conjecture of Buch.  (We rely on
Buch's extension of
the tableau rule \cite{buch}.)  We then describe
progress of extending this method to the flag manifold (the open question
of a Littlewood-Richardson rule for Schubert polynomials).  Finally,
we conclude with open questions.

\bpoint{Checkers, puzzles, tableaux}
\label{cpt} \lremind{cpt}
A bijection between checker games and puzzles is given in
Section~\ref{checkerpuzzle}.   Combining this with
Tao's ``proof-without-words'' of a bijection
between puzzles and tableaux (given in Figure~\ref{punt})
yields the bijection between checker-games and tableaux:

\tpoint{Theorem (bijection from checker games to tableaux)} {\em The
  construction of Section~\ref{bijst} gives a bijection from checker
  games to tableaux.}

There is undoubtedly a simpler direct proof (given the elegance of this map,
and the inelegance of the bijection from checkers to puzzles).

\begin{figure}[ht]
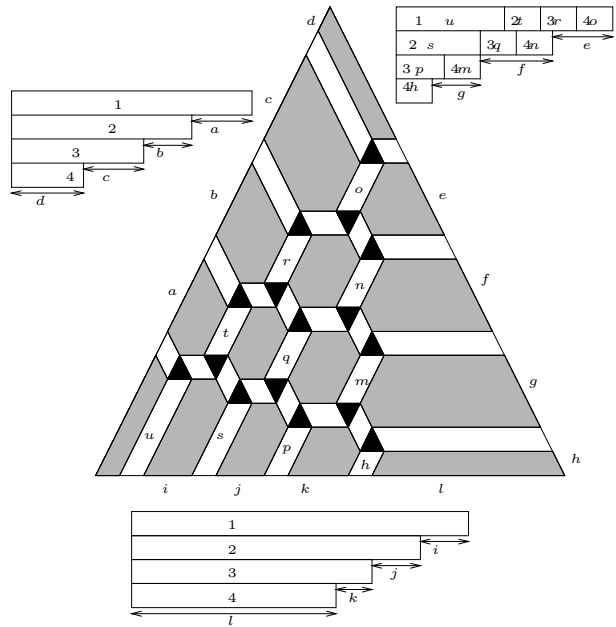

\begin{center}
\include{terry4}
\end{center}
\caption{Tao's ``proof without words'' of the bijection between puzzles and tableaux ($1$-triangles are depicted as black, regions of
$0$-rectangles are grey, and regions of rhombi are white)\lremind{punt, fig is terry4}}
\label{punt}
\end{figure}

Hence checker games give the first geometric interpretation of tableaux and
puzzles; indeed there is a bijection between tableaux/puzzles
and solutions of the corresponding triple-intersection Schubert problem,
once branch paths are chosen \cite[Sect.~4.3]{si}, 
\cite{sv}.\secretnote{Numerical 
Solutions, likely numsol}  

Note that to each puzzle, there are three possible checker-games,
depending on the orientation of the puzzle.  These correspond to three
degenerations of three general flags.  It would be interesting to
relate these three degenerations.

\bpoint{$K$-theory: checkers, puzzles, tableaux} 
\label{kcpt} \lremind{kcpt}
Buch \cite{bpc}
has conjectured that checker-game analysis can be extended to
$K$-theory or the Grothendieck ring (see \cite{buch} for background on
the $K$-theory of the Grassmannian).  Precisely, the rules for checker
moves are identical, except there is a new term in the middle square
of Table~\ref{keywest} (the case where there is a choice of moves), of
one lower dimension, with a minus sign.  If the two white checkers in
question are at $(r_1,c_1)$ and $(r_2,c_2)$, with $r_1>r_2$ and
$c_1<c_2$, then they move to $(r_2,c_1)$ and $(r_1-1,c_2)$ (see
Figure~\ref{sacramento}).  Call this a {\em sub-swap}; denote the
resulting configuration $\css \bn$.  Note that by Lemma~\ref{dimY},
$\dim \overline{Y}_{\css \bn} = \dim \overline{Y}_{\cb}-1$.

\begin{figure}[ht]
\begin{center}\setlength{\unitlength}{0.00083333in}
\begingroup\makeatletter\ifx\SetFigFont\undefined%
\gdef\SetFigFont#1#2#3#4#5{%
  \reset@font\fontsize{#1}{#2pt}%
  \fontfamily{#3}\fontseries{#4}\fontshape{#5}%
  \selectfont}%
\fi\endgroup%
{\renewcommand{\dashlinestretch}{30}
\begin{picture}(1362,939)(0,-10)
\put(1238,574){\ellipse{76}{76}}
\put(825,124){\ellipse{76}{76}}
\put(750,124){\blacken\ellipse{76}{76}}
\put(750,124){\ellipse{76}{76}}
\put(563,349){\blacken\ellipse{76}{76}}
\put(563,349){\ellipse{76}{76}}
\put(1013,574){\blacken\ellipse{76}{76}}
\put(1013,574){\ellipse{76}{76}}
\put(1238,799){\blacken\ellipse{76}{76}}
\put(1238,799){\ellipse{76}{76}}
\path(563,387)(563,574)
\blacken\path(578.000,514.000)(563.000,574.000)(548.000,514.000)(578.000,514.000)
\path(1013,537)(1013,349)
\blacken\path(998.000,409.000)(1013.000,349.000)(1028.000,409.000)(998.000,409.000)
\path(825,162)(825,574)
\whiten\path(840.000,514.000)(825.000,574.000)(810.000,514.000)(840.000,514.000)
\dashline{60.000}(825,574)(563,574)
\whiten\path(623.000,589.000)(563.000,574.000)(623.000,559.000)(623.000,589.000)
\path(450,912)(1350,912)(1350,12)
	(450,12)(450,912)
\path(900,912)(900,12)
\path(675,912)(675,12)
\path(1125,912)(1125,12)
\path(450,687)(1350,687)
\path(450,462)(1350,462)
\path(450,237)(1350,237)
\path(1238,537)(1238,331)
\whiten\path(1223.000,391.000)(1238.000,331.000)(1253.000,391.000)(1223.000,391.000)
\put(0,462){\makebox(0,0)[lb]{\smash{{{\SetFigFont{5}{6.0}{\rmdefault}{\mddefault}{\updefault}$(-1) \times$}}}}}
\end{picture}
}\end{center}
\caption{Buch's ``sub-swap'' case for the
$K$-theory geometric Littlewood-Richardson rule (cf. Figure~\ref{durham})
\lremind{sacramento}}
\label{sacramento}
\end{figure}
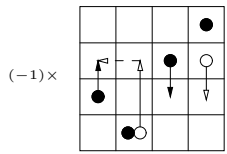

\tpoint{Theorem ($K$-theory Geometric Littlewood-Richardson rule)}
{\em Buch's rule describes multiplication in the Grothendieck
ring of $G(k,n)$.} \label{kglr} \lremind{kglr}

\bpf We give a bijection from $K$-theory checker games to Buch's
``set-valued tableaux'' (certain tableaux whose entries are sets of
consecutive integers, \cite{buch}), generalizing the bijection of
Section~\ref{bijst}.  Each white checker now has a memory of certain
earlier rows.  Each time there is a sub-swap, where a checker rises
from being the $r^{\text{th}}$ white checker to being the
$(r-1)^{\text{st}}$ (counting by row), that checker adds to its memory
that it had once been the $r^{\text{th}}$ checker (by row).  Whenever
there is move described by a $\dagger$ in Figure~\ref{durham}, where
the white checker is the $r^{\text{th}}$ by row and the
$c^{\text{th}}$ by column, in row $c$ of the
tableau place the set consisting of $r$ and all remembered earlier rows.
Then erase the memory of that white checker.
(The reader may verify that in Figure~\ref{g24} ($G(2,4)$), 
the result is an additional set-valued tableau, with a single cell
containing the set $\{ 1, 2 \}$.)

The proof
that this is a bijection is straightforward and left to the
reader.  \epf

This result suggests that Buch's rule reflects a geometrically
stronger fact, 
extending 
the final version of the \LR{} \ref{glr3}.

\tpoint{Conjecture ($K$-theory \LR, geometric form, with A. Buch)}
\label{bvconjecture} \lremind{bvconjecture}
{\em \begin{enumerate}
\item[(a)]  In the Grothendieck ring, $[X_{\cb}] = [\overline{X}_{\cbn}]$,
$[\overline{X}_{\cbs}]$,
  or $[\overline{X}_{ \cbn}]+ [\overline{X}_{ \cbs} ]
                      - [\overline{X}_{ \css \bn}  ]$.
\item[(b)] Scheme-theoretically, 
  $D_X = \overline{X}_{\cbn}$,
$\overline{X}_{\cbs}$,
  or $\overline{X}_{ \cbn} \cup \overline{X}_{ \cbs}$.
In the latter case, the scheme-theoretic intersection
$\overline{X}_{ \cbn} \cap \overline{X}_{ \cbs}$
is a translate of $\overline{X}_{ \css \bn}$. 
\end{enumerate}
}

Part (a) clearly follows from part (b).

Knutson has speculated that the total space of the degeneration
is Cohen-Macaulay; this would imply the conjecture.

The $K$-theory \LR~\ref{kglr} can be extended to puzzles. 

\tpoint{Theorem ($K$-theory Puzzle Littlewood-Richardson rule)} 
\lremind{kglrp} \label{kglrp}
{\em The
  $K$-theory Littlewood-Richardson coefficient corresponding to subsets
  $\al$, $\be$, $\ga$ is the number of puzzles with sides given by
  $\al$, $\be$, $\ga$ completed with the pieces shown in
  Figure~\ref{sanluisobispo}.  There is a factor of $-1$ for each
$K$-theory piece in the puzzle.}

\begin{figure}[ht]
\begin{center}\setlength{\unitlength}{0.00083333in}
\begingroup\makeatletter\ifx\SetFigFont\undefined%
\gdef\SetFigFont#1#2#3#4#5{%
  \reset@font\fontsize{#1}{#2pt}%
  \fontfamily{#3}\fontseries{#4}\fontshape{#5}%
  \selectfont}%
\fi\endgroup%
{\renewcommand{\dashlinestretch}{30}
\begin{picture}(4824,1110)(0,-10)
\path(12,387)(612,387)(312,912)(12,387)
\path(1212,387)(1812,387)(1512,912)(1212,387)
\path(3612,1062)(4812,1062)(4212,12)(3612,1062)
\path(2412,537)(2712,1062)(3012,537)
	(2712,12)(2412,537)
\put(237,312){\makebox(0,0)[lb]{\smash{{{\SetFigFont{10}{12.0}{\rmdefault}{\mddefault}{\updefault}$0$}}}}}
\put(87,612){\makebox(0,0)[lb]{\smash{{{\SetFigFont{10}{12.0}{\rmdefault}{\mddefault}{\updefault}$0$}}}}}
\put(462,612){\makebox(0,0)[lb]{\smash{{{\SetFigFont{10}{12.0}{\rmdefault}{\mddefault}{\updefault}$0$}}}}}
\put(1287,612){\makebox(0,0)[lb]{\smash{{{\SetFigFont{10}{12.0}{\rmdefault}{\mddefault}{\updefault}$1$}}}}}
\put(1662,612){\makebox(0,0)[lb]{\smash{{{\SetFigFont{10}{12.0}{\rmdefault}{\mddefault}{\updefault}$1$}}}}}
\put(1437,312){\makebox(0,0)[lb]{\smash{{{\SetFigFont{10}{12.0}{\rmdefault}{\mddefault}{\updefault}$1$}}}}}
\put(2487,762){\makebox(0,0)[lb]{\smash{{{\SetFigFont{10}{12.0}{\rmdefault}{\mddefault}{\updefault}$1$}}}}}
\put(2862,762){\makebox(0,0)[lb]{\smash{{{\SetFigFont{10}{12.0}{\rmdefault}{\mddefault}{\updefault}$0$}}}}}
\put(2862,237){\makebox(0,0)[lb]{\smash{{{\SetFigFont{10}{12.0}{\rmdefault}{\mddefault}{\updefault}$1$}}}}}
\put(2487,237){\makebox(0,0)[lb]{\smash{{{\SetFigFont{10}{12.0}{\rmdefault}{\mddefault}{\updefault}$0$}}}}}
\put(3687,762){\makebox(0,0)[lb]{\smash{{{\SetFigFont{10}{12.0}{\rmdefault}{\mddefault}{\updefault}$1$}}}}}
\put(3987,237){\makebox(0,0)[lb]{\smash{{{\SetFigFont{10}{12.0}{\rmdefault}{\mddefault}{\updefault}$0$}}}}}
\put(4362,237){\makebox(0,0)[lb]{\smash{{{\SetFigFont{10}{12.0}{\rmdefault}{\mddefault}{\updefault}$1$}}}}}
\put(4662,762){\makebox(0,0)[lb]{\smash{{{\SetFigFont{10}{12.0}{\rmdefault}{\mddefault}{\updefault}$0$}}}}}
\put(3912,987){\makebox(0,0)[lb]{\smash{{{\SetFigFont{10}{12.0}{\rmdefault}{\mddefault}{\updefault}$0$}}}}}
\put(4362,987){\makebox(0,0)[lb]{\smash{{{\SetFigFont{10}{12.0}{\rmdefault}{\mddefault}{\updefault}$1$}}}}}
\end{picture}
}\end{center}
\caption{The $K$-theory puzzle pieces
\lremind{sanluisobispo}}
\label{sanluisobispo}
\end{figure}
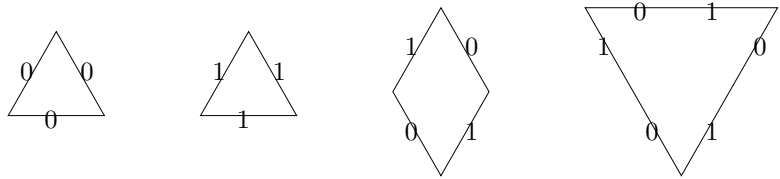

The first three pieces of Figure~\ref{sanluisobispo} are the usual
puzzle pieces of \cite{ktw,kt}; they may be rotated.  The fourth piece
is new; it may not be rotated.  Tao had earlier,
independently, discovered this piece \cite{tpc}.

Theorem~\ref{kglrp} may be proved via the $K$-theory \LR~\ref{kglr}, or
by generalizing Tao's proof of Figure~\ref{punt}.  Both proofs
are omitted.

As a consequence, we immediately have:

\tpoint{Corollary (triality of $K$-theory Littlewood-Richardson
  coefficients)}  {\em If 
  $K$-theory Littlewood-Richardson coefficients are denoted $C^{\cdot}_{\cdot \cdot}$, 
$C_{\al \be}^{\ga^{\vee}} = C_{\be
    \ga}^{\al^{\vee}} = C_{ \ga \al}^{\be^{\vee}}$.}

This is immediate in cohomology, but not obvious in the Grothendieck
ring.  The following direct proof is due to Buch (cf. \cite[p.~30]{buch}).
%\remind{Ask Anders if he minds, and if I've transcribed
%this correctly.  What's the symbol for $c$ in
%$K$-theory? Answer:  usually $c$, but I don't want to do this.}

\bpf
Let $\rho : G(d,n) \rightarrow pt$ be the map to a point.  Define a
pairing on $K_0(X)$ by $(a,b) := \rho_*(a \cdot b)$.  This pairing is
perfect, but (unlike for cohomology) the Schubert structure sheaf basis
is not dual to itself.  However, if $t$ denotes the top exterior power
of the tautological subbundle on $G(k,n)$, then the dual basis to the
structure sheaf basis is
$$    \{ t \oh_Y : Y \text{ is a Schubert variety in } G(k,n) \}.$$
More precisely, the structure sheaf for a 
partition $\lambda = (\lambda_1, \dots, \lambda_k)$ is
dual to $t$ times the structure 
sheaf for $\lambda^\vee = (n-k- \lambda_k, ..., n-k- \lambda_1)$.
(For more details, see \cite[Sect.~8]{buch}; 
this property is special for Grassmannians.)
Hence  $ \rho_*(t \oh_{\lambda} \oh_\mu \oh_nu)             
= C^{\nu^\vee}_{\lambda \mu} = C^{\lambda^\vee}_{\mu \nu} =
C^{\mu^\vee}_{\nu \lambda}$.
\epf

\bpoint{Toward a \LR{} for flag manifolds} \label{flint}
\lremind{flint} The same methods can be applied to flag manifolds, in
the hope of addressing the important open problem of finding a
Littlewood-Richardson rule in this context (i.e. structure
coefficients for the multiplication of Schubert polynomials), see
\cite{stanley}, \cite[p.~180]{fulton},
\cite[Sect.~9.10]{pr}.\secretnote{See Jan. '03 e-mail from Sara} This
problem has already motivated a great deal of remarkable work.

 We
describe some partial results here, informally.
Define two-flag Schubert varieties $X_{\de}$ of flag manifolds analogously to
$X_{\circ \bullet}$, as the locus in $Fl(n) \times
Fl(n) \times Fl(n)$ parametrizing (i) two flags $\ff_{\bu}$ and
$\fm_{\bu}$ in relative position given by $\bullet$, and (ii) a third
flag $F_{\bu}$ such that $\dim \ff_i \cap \fm_j \cap F_k = \de_{ijk}$.
The configuration $\bullet$ is obtained from the rank table $(\de_{ijn}){i,j}$
by the bijection of Section~\ref{domdefhere}.
The choice of the data $\de_{ijk}$ for which this data is non-empty can
be usefully summarized in two (equivalent) ways.   It is not yet clear
which is the more convenient notation.

First, these $\de$ correspond to configurations of checkers on
an $n \times n$ board, where there are $i$ checkers labeled $i$ for $1
\leq i \leq n$.    No two checkers with the same label are in the same
column or the same row.  An $i$-checker ($i<n$) is {\em happy}  if
there is an $(i+1)$-checker to its left (or in the same square) and an $(i+1)$-checker above
it (or in the same square).
(See for example Figure~\ref{sanmateo}.)  The bijection to $\de$ is given as follows:
$\de_{ijk}$ is the number of $k$-checkers dominated by $(i,j)$.
The morphism to $X_{\bullet}$ corresponds to forgetting all but the
$n$-checkers.

Second, these $\de$ correspond to ``wiring diagrams'', with
$n$ wires entering the board from below and leaving the board on the
right.  The bijection to $\de$ is given as follows:  $\de_{ijk}$ is
the number of wires numbered at least $k$ dominated by $(i,j)$.
See Figure~\ref{sanmateo2} for an example.

\begin{figure}[ht]
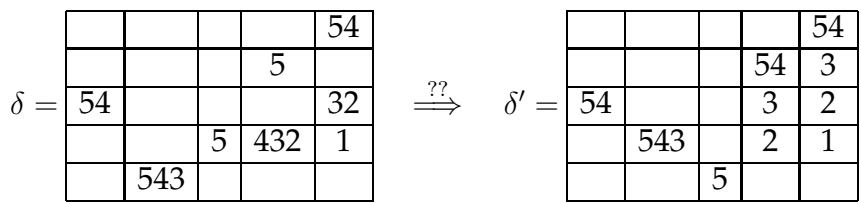

\begin{center}
$\de = $ 
\begin{tabular}{|c|c|c|c|c|} \hline
& & & & 54 \\ \hline
& & & 5 & \\ \hline
54 & & & & 32 \\ \hline
& & 5 & 432 & 1 \\ \hline
& 543 & & & \\ \hline
\end{tabular} 
$\quad \stackrel {??} {\Longrightarrow} \quad \de' =$
\begin{tabular}{|c|c|c|c|c|} \hline
& & & & 54 \\ \hline
& & & 54 & 3 \\ \hline
54 & & & 3 & 2 \\ \hline
&  543 &  & 2 & 1 \\ \hline
& & 5 & &  \\ \hline
\end{tabular} 
\end{center}
\caption{Is the two-flag Schubert variety on the right in the
closure of the one on the left?\lremind{sanmateo}}
\label{sanmateo}
\end{figure}

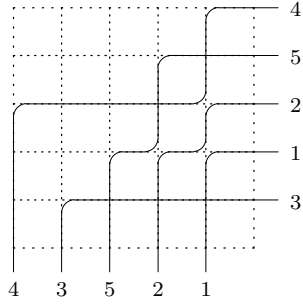
\begin{figure}[ht]
\begin{center}
\setlength{\unitlength}{0.00083333in}
\begingroup\makeatletter\ifx\SetFigFont\undefined%
\gdef\SetFigFont#1#2#3#4#5{%
  \reset@font\fontsize{#1}{#2pt}%
  \fontfamily{#3}\fontseries{#4}\fontshape{#5}%
  \selectfont}%
\fi\endgroup%
{\renewcommand{\dashlinestretch}{30}
\begin{picture}(1760,1875)(0,-10)
\put(1153.000,1313.000){\arc{164.596}{0.0852}{1.4856}}
\put(1317.833,1148.167){\arc{166.405}{3.2359}{4.6181}}
\put(1317.833,848.167){\arc{166.405}{3.2359}{4.6181}}
\put(1153.000,1013.000){\arc{164.596}{0.0852}{1.4856}}
\put(853.000,1013.000){\arc{164.596}{0.0852}{1.4856}}
\put(1017.833,848.167){\arc{166.405}{3.2359}{4.6181}}
\put(1017.833,1448.167){\arc{166.405}{3.2359}{4.6181}}
\put(1317.833,1748.167){\arc{166.405}{3.2359}{4.6181}}
\put(117.833,1148.167){\arc{166.405}{3.2359}{4.6181}}
\put(717.833,848.167){\arc{166.405}{3.2359}{4.6181}}
\put(417.833,548.167){\arc{166.405}{3.2359}{4.6181}}
\dottedline{60}(35,1831)(1535,1831)(1535,331)
	(35,331)(35,1831)
\dottedline{45}(35,1531)(1535,1531)
\dottedline{45}(335,1831)(335,331)
\dottedline{45}(35,1231)(1535,1231)
\dottedline{45}(35,931)(1535,931)
\dottedline{45}(35,631)(1535,631)
\dottedline{45}(635,1831)(635,331)
\dottedline{45}(935,1831)(935,331)
\dottedline{45}(1235,1831)(1235,331)
\put(1200,31){\makebox(0,0)[lb]{\smash{{{\SetFigFont{8}{9.6}{\rmdefault}{\mddefault}{\updefault}$1$}}}}}
\put(900,31){\makebox(0,0)[lb]{\smash{{{\SetFigFont{8}{9.6}{\rmdefault}{\mddefault}{\updefault}$2$}}}}}
\put(600,31){\makebox(0,0)[lb]{\smash{{{\SetFigFont{8}{9.6}{\rmdefault}{\mddefault}{\updefault}$5$}}}}}
\put(300,31){\makebox(0,0)[lb]{\smash{{{\SetFigFont{8}{9.6}{\rmdefault}{\mddefault}{\updefault}$3$}}}}}
\put(0,31){\makebox(0,0)[lb]{\smash{{{\SetFigFont{8}{9.6}{\rmdefault}{\mddefault}{\updefault}$4$}}}}}
\put(1760,1776){\makebox(0,0)[lb]{\smash{{{\SetFigFont{8}{9.6}{\rmdefault}{\mddefault}{\updefault}$4$}}}}}
\put(1760,1476){\makebox(0,0)[lb]{\smash{{{\SetFigFont{8}{9.6}{\rmdefault}{\mddefault}{\updefault}$5$}}}}}
\put(1760,1176){\makebox(0,0)[lb]{\smash{{{\SetFigFont{8}{9.6}{\rmdefault}{\mddefault}{\updefault}$2$}}}}}
\put(1760,876){\makebox(0,0)[lb]{\smash{{{\SetFigFont{8}{9.6}{\rmdefault}{\mddefault}{\updefault}$1$}}}}}
\put(1760,576){\makebox(0,0)[lb]{\smash{{{\SetFigFont{8}{9.6}{\rmdefault}{\mddefault}{\updefault}$3$}}}}}
\path(1235,181)(1235,856)
\path(1310,931)(1685,931)
\path(935,181)(935,856)
\path(1010,931)(1160,931)
\path(1235,1006)(1235,1156)
\path(1310,1231)(1685,1231)
\path(1010,631)(1160,631)
\path(1310,631)(1685,631)
\path(635,181)(635,556)
\path(935,1006)(935,1156)
\path(935,1306)(935,1456)
\path(1010,1531)(1160,1531)
\path(1310,1531)(1685,1531)
\path(1235,1306)(1235,1756)
\path(1310,1831)(1685,1831)
\path(110,1231)(1160,1231)
\path(35,1156)(35,181)
\path(410,631)(860,631)
\path(710,931)(860,931)
\path(335,556)(335,181)
\path(635,856)(635,706)
\path(635,706)(635,556)
\path(860,631)(1010,631)
\path(1160,631)(1310,631)
\path(935,1306)(935,1156)
\path(1160,1531)(1310,1531)
\end{picture}
}
\end{center}
\caption{The ``wiring diagram'' corresponding to the first checker configuration of Figure~\ref{sanmateo}\lremind{sanmateo2}}
\label{sanmateo2}
\end{figure}

There is an analogue of Proposition~\ref{ip}, describing the
intersection of two Schubert varieties (with respect to transverse flags
$\ff_{\bu}$ and $\fm_{\bu}$) as $\overline{X}_{\de}$ for an
explicitly given $\de$; call the set of $\de$ obtained in such a way
the {\em initial positions}.

\tpoint{Conjecture (Existence of a \LR{} for flag manifolds)} 
\label{conjecture} \lremind{conjecture}
{\em There exists a subset $M$ of the $\{ (\de_{ijk})_{ijk} \}$, where
  the configuration of $n$-checkers  is in the specialization order (analogous to ``mid-sort''),
 containing the set of initial positions,
  such that the divisor $D_X$ on $\overline{X}_{\de}$ corresponding to the
next element of the specialization order is the union of varieties
  of the form $\overline{X}_{\de'}$, each
  appearing with multiplicity 1, where $(\de'_{ijk})_{ijk} \in M$.}

This conjecture (which parallels Theorem~\ref{glr3}) looks quite weak,
as it specifies neither (i) the ``mid-sort'' property $M$ nor (ii) how
to determine $\de'$.  However, it suffices to obtain almost all of the
applications described in \cite{si} (e.g. reality, number of solutions
in positive characteristic, numerical calculation of solutions to all
Schubert problems, and more).\secretnote{Refer to entire article.}
Furthermore, it implies the existence of a Littlewood-Richardson rule,
and an answer to (ii) would give an explicit rule.

\tpoint{Proposition} {\em Conjecture~\ref{conjecture} is true for $n \leq 5$.}
\label{galveston} \lremind{galveston}

This is a considerable amount of evidence, involving $\sum_{n=1}^5
\binom n 2 \sum_{\al, \be, \ga} c^{\ga}_{\al \be}$ degenerations.
However, there are two serious reasons to remain suspicious: (a) 
Knutson's puzzle variant computes Littlewood-Richardson coefficients
for $n<5$ but fails for $n=5$ (and puzzles are related to checkers via
Sect.~\ref{checkerpuzzle}), and (b) all Littlewood-Richardson coefficients
are $0$ or $1$ for $n<6$.

{\noindent {\em Sketch of proof of Proposition~\ref{galveston} for $n \leq 4$.  }} We build $M$
inductively, starting with the set of initial positions.  Successively,
for each element $\de$ of $M$, we find a list of possible divisors
$\overline{X}_{\de'}$ on $\overline{X}_{\de}$ (above $X_{\bn}$) as follows.  If the
$n$-checkers of $\de$ are in position $\bullet$, the $n$-checkers of
$\de'$ are in position $\bn$.  Then given the positions of the
$k$-checkers of $\de$ and the $(k+1)$-checkers of $\de'$ (for $1 \leq
k<n$), we list the (finite number of) possibilities of the choice
of positions of the choices of $k$-checkers of $\de'$ subject to two conditions:  the $k$-checkers  are
happy, and $\de'_{ijk} \geq \de_{ijk}$.    We then discard all $\de'$
such that $\dim X_{\de'} < \dim X_{\de}-1$.  Then each component of $D_X$
lies in $D_{\de'}$ for some remaining $\de'$.  If any $X_{\de'}$ has
dimension at least $\dim X_{\de}$, we stop, and the process fails.
Otherwise, $D_X$ must consist of a union of these $X_{\de'}$,
possibility with multiplicity.  We then add these $\de'$ to $M$.
(This is a greedy algorithm, which sometimes includes $\de'$ which
do not correspond to components of $D_X$.)
Thus the
conjecture is true for this $\de$ (with this choice of $M$) with the
possible exception of the multiplicity 1 claim.  We now repeat the
process with this enlarged $M$.

For $n \leq 4$, one checks (most easily by computer) that this process never fails.  Moreover,
as all Littlewood-Richardson coefficients for $n \leq 4$ are 1, the
multiplicity must be 1 at each stage. \epf

The author has a full description of the degenerations which actually
appear for $n \leq 4$, available upon request.  However, it is not
clear how to generalize this to a conjectural Littlewood-Richardson
rule.

For $n=5$, this process fails at six cases (where $\dim X_{\de'} =
\dim X_{\de}$).  However, in three of the six cases, it can be shown
that $\overline{X}_{\de}$ does not meet $X_{\de'}$ (one is shown in
Figure~\ref{sanmateo}; the other two are identical except for the
position of the 1-checker), and in the other three cases, ${\de}$ itself
can be removed from $M$ (i.e. $\de$ was falsely included in $M$ by the
greediness of the algorithm).

For larger $n$, this greedy algorithm will certainly produce even worse
pathologies (i.e. arbitrarily many cases where $\dim X_{\de'} - \dim X_{\de}$ is
arbitrarily large).

\bpoint{Questions}
\label{open} \lremind{open}  One motivation for the \LR{} is that it should
generalize well to other important geometric situations (as it has in
$K$-theory and at least partially for the flag manifold).
We now briefly describe some potential applications; some 
are work in progress.

\noindent
{\em (a)} These methods may apply to other groups where
Littlewood-Richardson rules are not known.  For example, for the
symplectic Grassmannian, there are only rules in the Lagrangian and
Pieri cases.  L. Mihalcea has made progress in finding a Geometric
Littlewood-Richardson rule in the Lagrangian case, and has suggested
that a similar algorithm should exist in general.

\noindent  %kpc
{\em (b)} An important intermediate stage between the Grassmannian and
the full flag manifold is the two-step partial flag manifold.
This case has useful applications, for example, to
Grassmannians of other groups.  The preprint \cite{bkt} gives a
connection to Gromov-Witten invariants (and does much more).  The
Littlewood-Richardson behavior of the two-step partial flag manifold
is much better than for full flag manifolds.  Buch, Kresch, and
Tamvakis have suggested that Knutson's proposed partial flag rule
(which fails for flags in general) holds for two-step flags, and have
verified this up to $n=16$ \cite[p.~6]{bkt}.

Is there a good (and
straightforward) checker-rule for such homogeneous 
spaces?\secretnote{vetted with anders}

\noindent
{\em (c)} Can equivariant Littlewood-Richardson coefficients be
understood geometrically in this way?  For example, can 
equivariant puzzles \cite{kt}
be translated to checkers, and can
partially-completed equivariant puzzles thus be given a geometric 
interpretation?  Can
this be combined with Theorem~\ref{kglr} to yield a
Littlewood-Richardson rule in equivariant $K$-theory?  

\noindent
{\em (d)} The quantum cohomology of the Grassmannian can be translated
into classical questions in the enumerative geometry of surfaces.  One
may hope that degeneration methods introduced here and in
\cite{crelle} will apply.  This perspective is being pursued (with
different motivation) independently by I. Coskun (for rational
scrolls), D. Avritzer,  and M. Honsen (on Veronese surfaces).

\noindent
{\em (e)} Is there any relation between the wiring diagrams of
Figure~\ref{sanmateo2} and rc-graphs?  (We suspect not.)

\noindent {\em (f)}
D. Eisenbud and J. Harris \cite{eh} have a
particular (irreducible, one-parameter) path in the flag variety,
whose general point is in the large open Schubert cell, and whose
special point is the smallest cell: consider the osculating flag
$\fm_{\bu}$ to a point $p$ on a rational normal curve, as $p$ tends to
a reference point $q$ with osculating flag $\ff_{\bu}$.  Eisenbud has
asked if the specialization order is some sort of limit (a
``polygonalization'') of such paths. This would provide a single path
that breaks intersections of Schubert cells into their components.
(Of course, the limit cycles could not have multiplicity one in
general.)  Eisenbud and Harris' proof of the Pieri formula is evidence
that this could be true.

Sottile has a precise conjecture generalizing Eisenbud and Harris'
approach to all flag manifolds \cite[Sect.~5]{sottileMMJ}.  He has
generalized this further: one replaces the rational normal curve by
the curve $e^{t \eta} X_u(F_{\cdot})$, where $\eta$ is a principal
nilpotent in the Lie algebra of the respective algebraic group, and
the limit is then $\lim_{t \rightarrow 0} e^{t \eta} X_u(F_{\cdot})
\cap X_w$, where $X_w$ is given by the flag fixed by $\lim_{t
\rightarrow 0} e^{t \eta}$, \cite{spc}.
Eisenbud's question in this context then involves polygonalizing or
degenerating this
path.

\section{Bott-Samelson \Qflag Varieties}
\label{bsv} \lremind{bsv}
\point \label{bsvs} \lremind{bsvs}
We will associate a variety to the following data. 
\begin{itemize}
\item $n$ is a positive integer.
\item $\poset$ is a finite subset of the plane (visualized so that downwards
corresponds to increasing the first coordinate and rightwards
corresponds to increasing the second coordinate, in keeping
with the labeling convention for tables),
with the partial order $\prec$ given by {domination} 
(defined in Sect.~\ref{domdefhere}).
\item $\dim:  \poset \rightarrow \{ 0, 1, 2, \dots, n \}$ 
is a morphism of posets (i.e. weakly order-preserving map), denoted {\em dimension}.\notation{dimension}
\item If $[\mathbf{a}, \mathbf{b}]$ is a covering relation (i.e. minimal
  interval) in $\poset$ (i.e. $\mathbf{a}, \mathbf{b} \in \poset$,
  $\mathbf{a} \prec \mathbf{b}$, and there is no $\mathbf{c} \in
  \poset$ such that $\mathbf{a} \prec \mathbf{c} \prec \mathbf{b}$),
  then $\dim \mathbf{a} = \dim \mathbf{b} - 1$.  (This condition is
  likely unnecessary.)

%\remind{Anders points out that this is unnecessary.}

\item $\poset$ has a maximum element and a minimum element.
\end{itemize}
We call this data 
a {\em planar poset}, and denote it by $\poset$; the remaining data $( \prec ,  \dim, n)$ will be implicit.\notation{planar poset, $\poset$}
%\remind{Anders points out that this can be defined in terms of posets,
%and then they can be reducible; his example had two high and two low elements
%with each high dominating each low.}

It will be convenient to represent this data as a planar
graph, whose vertices are elements of $\poset$, and whose edges
correspond to covering relations in  $\poset$. 
The interior of the graph is a union of
{quadrilaterals}; call these the {\em quadrilaterals 
of $\poset$}.
An element of $\poset$ at $(i,j)$ is said to be on the {\em southwest border}
(resp. {\em northeast border}) if there are no elements of $\poset$
$(i',j')$ such that $i' \geq i$ and $j' \leq j$ (resp. $i' \leq i$ and
$j' \geq j$); see Figure~\ref{podunk}.\notation{quadrilateral; southwest and northeast border}

Define the {\em Bott-Samelson \Qflag variety} $\PF(\poset)$
to be the variety parametrizing
a $(\dim \mathbf{s})$-plane $V_{\mathbf{s}}$ in $K^n$ for each $\mathbf{s} \in \poset$, with $V_{\mathbf{s}} \subset V_{\mathbf{t}}$
for $\mathbf{s} \prec \mathbf{t}$.  (It is a closed subvariety of 
$\prod_{\mathbf{s} \in \poset} G(\dim \mathbf{s}, n)$.)\notation{Bott-Samelson \Qflag variety, $\PF(\poset)$}
Elements $\mathbf{s}$ of $\poset$ will be written in bold-faced font;
corresponding vector spaces will be denoted $V_{\mathbf{s}}$.

Forgetting all but the vertices on the southwest border yields a
morphism to the flag manifold, and the usual Bott-Samelson variety is a
fiber of this morphism.

\tpoint{Lemma}  \lremind{bssmooth} {\em The Bott-Samelson \Qflag variety $\PF(\poset)$ is smooth.}
\label{bssmooth}

\bpf Consider the planar graph representation of $\poset$ described
above.  The variety parametrizing the subspaces corresponding to the
southwest border of the graph is a partial flag
variety (and hence smooth).  The Bott-Samelson \Qflag variety $\PF(\poset)$
can be expressed as a tower of $\proj^1$-bundles over the partial flag
variety by inductively adding the data of elements of $\mathbf{s} \in S$
corresponding to ``new'' vertices of quadrilaterals (where the other
three vertices $\mathbf{s_1} \prec \mathbf{s_2} \prec \mathbf{s_3}$ are already parametrized, and
$\mathbf{s_1} \prec \mathbf{s} \prec \mathbf{s_3}$).  \epf

For example, Figure~\ref{podunk} illustrates that one particular
Bott-Samelson \Qflag variety is a tower of five $\proj^1$-bundles over
$Fl(4)$; the correspondence of the $\proj^1$-bundles with
quadrilaterals is illustrated by the arrows.

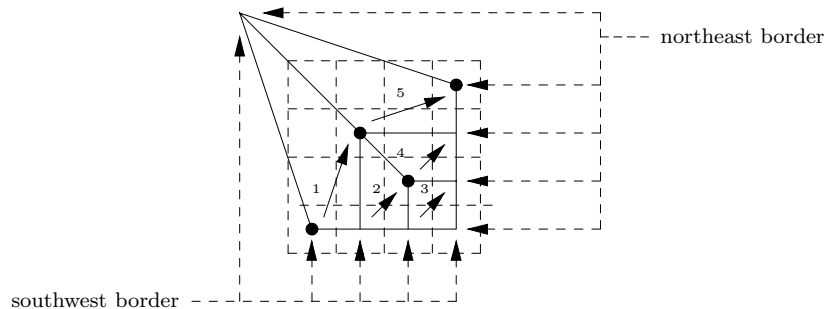
\begin{figure}[ht]
\begin{center} \setlength{\unitlength}{0.00083333in}
\begingroup\makeatletter\ifx\SetFigFont\undefined%
\gdef\SetFigFont#1#2#3#4#5{%
  \reset@font\fontsize{#1}{#2pt}%
  \fontfamily{#3}\fontseries{#4}\fontshape{#5}%
  \selectfont}%
\fi\endgroup%
{\renewcommand{\dashlinestretch}{30}
\begin{picture}(4050,1894)(0,-10)
\put(2775,1387){\blacken\ellipse{74}{74}}
\put(2775,1387){\ellipse{74}{74}}
\put(2475,787){\blacken\ellipse{74}{74}}
\put(2475,787){\ellipse{74}{74}}
\put(2175,1087){\blacken\ellipse{74}{74}}
\put(2175,1087){\ellipse{74}{74}}
\put(1875,487){\blacken\ellipse{74}{74}}
\put(1875,487){\ellipse{74}{74}}
\dashline{60.000}(1725,1537)(2925,1537)(2925,337)
	(1725,337)(1725,1537)
\dashline{60.000}(2025,1537)(2025,337)
\dashline{60.000}(1725,1237)(2925,1237)
\dashline{60.000}(2325,1537)(2325,337)
\dashline{60.000}(2625,1537)(2625,337)
\dashline{60.000}(1800,637)(3000,637)
\dashline{60.000}(1725,937)(2925,937)
\path(1425,1837)(2775,1387)(2775,487)
	(1875,487)(1425,1837)(2175,1087)
	(2475,787)(2775,787)
\path(2175,487)(2175,1087)(2775,1087)
\path(2475,787)(2475,487)
\path(1950,562)(2100,1012)
\blacken\path(2090.513,888.671)(2100.000,1012.000)(2033.592,907.645)(2090.513,888.671)
\path(2250,562)(2400,712)
\blacken\path(2336.360,605.934)(2400.000,712.000)(2293.934,648.360)(2336.360,605.934)
\path(2550,562)(2700,712)
\blacken\path(2636.360,605.934)(2700.000,712.000)(2593.934,648.360)(2636.360,605.934)
\path(2550,862)(2700,1012)
\blacken\path(2636.360,905.934)(2700.000,1012.000)(2593.934,948.360)(2636.360,905.934)
\path(2250,1162)(2700,1312)
\blacken\path(2595.645,1245.592)(2700.000,1312.000)(2576.671,1302.513)(2595.645,1245.592)
\dashline{60.000}(1425,37)(1425,1687)
\blacken\path(1455.000,1567.000)(1425.000,1687.000)(1395.000,1567.000)(1455.000,1567.000)
\dashline{60.000}(1875,37)(1875,412)
\blacken\path(1905.000,292.000)(1875.000,412.000)(1845.000,292.000)(1905.000,292.000)
\dashline{60.000}(2175,37)(2175,412)
\blacken\path(2205.000,292.000)(2175.000,412.000)(2145.000,292.000)(2205.000,292.000)
\dashline{60.000}(2475,37)(2475,412)
\blacken\path(2505.000,292.000)(2475.000,412.000)(2445.000,292.000)(2505.000,292.000)
\dashline{60.000}(2775,37)(2775,412)
\blacken\path(2805.000,292.000)(2775.000,412.000)(2745.000,292.000)(2805.000,292.000)
\dashline{60.000}(3675,1837)(1575,1837)
\blacken\path(1695.000,1867.000)(1575.000,1837.000)(1695.000,1807.000)(1695.000,1867.000)
\dashline{60.000}(3675,1087)(2850,1087)
\blacken\path(2970.000,1117.000)(2850.000,1087.000)(2970.000,1057.000)(2970.000,1117.000)
\dashline{60.000}(3675,787)(2850,787)
\blacken\path(2970.000,817.000)(2850.000,787.000)(2970.000,757.000)(2970.000,817.000)
\dashline{60.000}(3675,487)(2850,487)
\blacken\path(2970.000,517.000)(2850.000,487.000)(2970.000,457.000)(2970.000,517.000)
\dashline{60.000}(3675,1837)(3675,487)
\dashline{60.000}(3675,1387)(2850,1387)
\blacken\path(2970.000,1417.000)(2850.000,1387.000)(2970.000,1357.000)(2970.000,1417.000)
\dashline{60.000}(3675,1687)(3975,1687)
\dashline{60.000}(2775,37)(1125,37)
\put(2550,712){\makebox(0,0)[lb]{\smash{{{\SetFigFont{5}{6.0}{\rmdefault}{\mddefault}{\updefault}$3$}}}}}
\put(2400,937){\makebox(0,0)[lb]{\smash{{{\SetFigFont{5}{6.0}{\rmdefault}{\mddefault}{\updefault}$4$}}}}}
\put(1875,712){\makebox(0,0)[lb]{\smash{{{\SetFigFont{5}{6.0}{\rmdefault}{\mddefault}{\updefault}$1$}}}}}
\put(2250,712){\makebox(0,0)[lb]{\smash{{{\SetFigFont{5}{6.0}{\rmdefault}{\mddefault}{\updefault}$2$}}}}}
\put(2400,1312){\makebox(0,0)[lb]{\smash{{{\SetFigFont{5}{6.0}{\rmdefault}{\mddefault}{\updefault}$5$}}}}}
\put(0,0){\makebox(0,0)[lb]{\smash{{{\SetFigFont{8}{9.6}{\rmdefault}{\mddefault}{\updefault}southwest border}}}}}
\put(4050,1650){\makebox(0,0)[lb]{\smash{{{\SetFigFont{8}{9.6}{\rmdefault}{\mddefault}{\updefault}northeast border}}}}}
\end{picture}
} \end{center}
\caption{The northeast and southwest borders of a planar poset generated by a checker
configuration; description of a Bott-Samelson \Qflag variety as a tower
of five $\proj^1$-bundles over $Fl(4)$ \lremind{podunk}}
\label{podunk}
\end{figure}

\epoint{Strata of Bott-Samelson \Qflag varieties} \label{bsvss}
\lremind{bsvss}  Any subset $Q$ of
the quadrilaterals of a planar poset
determines a {\em 
stratum}\notation{strata of a Bott-Samelson \Qflag variety} of the Bott-Samelson \Qflag variety.  The closed stratum corresponds to
requiring the subspaces of the opposite corners of the quadrilaterals in
$Q$ of the same dimension to be the same.  The open stratum
corresponds to also requiring the spaces of the opposite corners of
the quadrilaterals {\em not} in $Q$ to be distinct.
By the construction in the proof of Lemma~\ref{bssmooth}, 
(i) the open strata give a stratification, (ii) the
closed strata are smooth, and (iii) the codimension of the stratum is
the size of the subset $Q$.
It will be convenient to depict
a stratum by placing an ``$=$'' in the quadrilaterals
of $Q$, indicating the pairs of subspaces that are required
to be equal (see for example Figure~\ref{boston2}).

\epoint{Example: planar posets generated by a set of
  checkers}\label{huntsville}\lremind{huntsville} Given a checker
configuration $\bullet$ (or $\circ$) and a positive integer $n$, we
define the planar poset $\poset_{\bullet}$ (or
$\poset_{\circ}$) as follows.\notation{$\poset_{\bullet},
  \poset_{\circ}$} Include the squares of the table where there is a
checker above (or possibly in the same square), and a checker to the
left (or in the same square); include also a ``zero element'' above
and to the left of the checkers.  (The definition of {happy} in
Sect.~\ref{rulewhite} can be rephrased as: the white checkers lie on
elements of $\poset_{\bullet}$.)  For $\mathbf{s} \in \poset$, let
$\dim \mathbf{s}$ be the number of checkers dominated by $\mathbf{s}$.

For example, if $\bullet$ is a configuration of $n$ checkers 
(as in Sect.~\ref{seattle}), then 
the southwest border of $\poset_{\bullet}$ corresponds
to $\ff_{\bu}$, and the northeast border corresponds to $\fm_{\bu}$; 
$\PF(\poset_{\bullet})$ is a fibration over $Fl(n)$ (where
$Fl(n)$ parametrizes $\ff_{\bu}$), and the fiber is a Bott-Samelson
resolution of $\Omega_{\bullet}(\ff_{\bu})$.    
Figure~\ref{podunk} describes a Bott-Samelson resolution of
the double Schubert variety corresponding to $1324$.
Similarly, the morphism $\PF(\poset_{\bullet}) \rightarrow
\overline{X}_{\bullet}$ is a resolution of singularities.
This morphism
restricts to an isomorphism of the dense open stratum of
$\PF(\poset_{\bullet})$ with $X_{\bullet}$.
 If $\bullet$ is in the specialization order, then
$\PF(\poset_{\bn})$ is (isomorphic to) a  codimension 1 stratum of
this Bott-Samelson \Qflag variety.  See Figure~\ref{boston2} for an example.

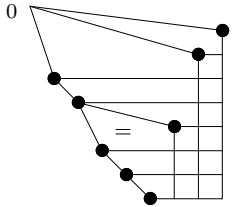
\begin{figure}[ht]
\begin{center} \setlength{\unitlength}{0.00083333in}
\begingroup\makeatletter\ifx\SetFigFont\undefined%
\gdef\SetFigFont#1#2#3#4#5{%
  \reset@font\fontsize{#1}{#2pt}%
  \fontfamily{#3}\fontseries{#4}\fontshape{#5}%
  \selectfont}%
\fi\endgroup%
{\renewcommand{\dashlinestretch}{30}
\begin{picture}(1396,1272)(0,-10)
\put(1350,1095){\blacken\ellipse{76}{76}}
\put(1350,1095){\ellipse{76}{76}}
\put(1200,945){\blacken\ellipse{76}{76}}
\put(1200,945){\ellipse{76}{76}}
\put(300,795){\blacken\ellipse{76}{76}}
\put(300,795){\ellipse{76}{76}}
\put(1050,495){\blacken\ellipse{76}{76}}
\put(1050,495){\ellipse{76}{76}}
\put(750,195){\blacken\ellipse{76}{76}}
\put(750,195){\ellipse{76}{76}}
\put(900,45){\blacken\ellipse{76}{76}}
\put(900,45){\ellipse{76}{76}}
\put(450,645){\blacken\ellipse{76}{76}}
\put(450,645){\ellipse{76}{76}}
\put(600,345){\blacken\ellipse{76}{76}}
\put(600,345){\ellipse{76}{76}}
\path(1350,1095)(1350,45)(900,45)
	(600,345)(450,645)(300,795)(1350,795)
\path(1350,945)(1200,945)(1200,45)
\path(450,645)(1350,645)
\path(450,645)(1050,495)(1350,495)
\path(600,345)(1350,345)
\path(1050,495)(1050,45)
\path(750,195)(1350,195)
\path(150,1245)(1350,1095)
\path(1200,945)(150,1245)(300,795)
\put(0,1170){\makebox(0,0)[lb]{\smash{{{\SetFigFont{8}{9.6}{\rmdefault}{\mddefault}{\updefault}$0$}}}}}
\put(675,432){\makebox(0,0)[lb]{\smash{{{\SetFigFont{8}{9.6}{\rmdefault}{\mddefault}{\updefault}$=$}}}}}
\end{picture}
} \end{center}
\caption{The poset corresponding to $\bullet$ of Figure~\ref{boston}, with the  divisorial stratum corresponding to $\bn$ marked with an ``$=$''
\lremind{boston2} 
} \label{boston2}
\end{figure}

\begin{figure}[ht]
\begin{center} \setlength{\unitlength}{0.00083333in}
\begingroup\makeatletter\ifx\SetFigFont\undefined%
\gdef\SetFigFont#1#2#3#4#5{%
  \reset@font\fontsize{#1}{#2pt}%
  \fontfamily{#3}\fontseries{#4}\fontshape{#5}%
  \selectfont}%
\fi\endgroup%
{\renewcommand{\dashlinestretch}{30}
\begin{picture}(1833,897)(0,-10)
\put(450,636){\blacken\ellipse{50}{50}}
\put(450,636){\ellipse{50}{50}}
\put(450,186){\blacken\ellipse{50}{50}}
\put(450,186){\ellipse{50}{50}}
\put(1800,186){\blacken\ellipse{50}{50}}
\put(1800,186){\ellipse{50}{50}}
\put(1800,636){\blacken\ellipse{50}{50}}
\put(1800,636){\ellipse{50}{50}}
\path(450,636)(1800,636)(1800,186)
	(450,186)(450,636)
\put(0,786){\makebox(0,0)[lb]{\smash{{{\SetFigFont{8}{9.6}{\rmdefault}{\mddefault}{\updefault}$\dim \mathbf{x} = n-2$}}}}}
\put(1425,786){\makebox(0,0)[lb]{\smash{{{\SetFigFont{8}{9.6}{\rmdefault}{\mddefault}{\updefault}$\dim \mathbf{y'} = n-1$}}}}}
\put(1425,36){\makebox(0,0)[lb]{\smash{{{\SetFigFont{8}{9.6}{\rmdefault}{\mddefault}{\updefault}$\dim \mathbf{z}=n$}}}}}
\put(0,36){\makebox(0,0)[lb]{\smash{{{\SetFigFont{8}{9.6}{\rmdefault}{\mddefault}{\updefault}$\dim \mathbf{y} = n-1$}}}}}
\end{picture}
} \end{center}
\caption{The planar poset $\pBox$ \lremind{cupertino}}
\label{cupertino}
\end{figure}

\bpoint{Bott-Samelson \Qflag varieties of morphisms of posets} 
\label{augusta} \lremind{augusta}
Suppose $\poset$
and $\cq$ are two planar posets (for the same $n$), and $p: \poset
\rightarrow \cq$ is a morphism of posets (i.e. a weakly order-preserving map
of sets, with no conditions on the planar structures).  Then let $\PF(p) = \PF(\poset \rightarrow \cq)$
\notation{$\PF(p) = \PF(\poset \rightarrow \cq)$} be the subvariety
of $\PF(\poset) \times \PF(\cq)$ such that if $p(\mathbf{s}) =
\mathbf{t}$, then the subspace corresponding to $\mathbf{s}$ is contained in the
subspace corresponding to $\mathbf{t}$.
For example, a configuration of black and white checkers $\cb$
induces a morphism of posets $p_{\cb}:  \poset_{\circ} \rightarrow
\poset_\bullet$\notation{$p_{\cb}$}, and 
$X_{\cb}$ is an open subset of $\PF(p_{\cb}:  \poset_{\circ} \rightarrow \poset_\bullet)$.

{\em Caution:} Unlike $\PF(\poset)$, the Bott-Samelson \Qflag variety $\PF( \poset \rightarrow
\cq)$ may be singular.  For example, if $\poset = \{ (x,y) : x, y \in
\{0,1\} \}$, $\cq = \{ (x,y) : x, y-1 \in \{0,1\} \}$, $p(x,y) =
(x,y+1)$, and $\dim(x,y) = x+y$, then $\PF(p)$ (with $n=3$) is the
triangle variety (parametrizing points $p_1$, $p_2$, $p_3$ and 
lines $\ell_{12}$, $\ell_{23}$, $\ell_{31}$ 
in $\proj^2$
with $p_i, p_j \in \ell_{ij}$) 
and hence singular.  Also, $X_{\cb}$ need not be
dense in $\PF(\poset_{\circ}
\rightarrow \poset_{\bullet})$;  see Section~\ref{cap} (b).

\bpoint{Relating $X_{\bullet} \cup X_{\bn}$ to the simpler variety
$\PF({\pBox})$ via $X_{\bB}$}
\label{portland}\lremind{portland}
Let ${\pBox}$ be the poset of Figure~\ref{cupertino}.
Then $\PF({\pBox})$ parametrizes
hyperplanes $V_{\y}$ and $V_{\yp}$, and a codimension 2 subspace $V_{\x} \subset
V_{\y}, V_{\yp}$.  (Of course $V_{\z} = K^n$.)
 The variety $\PF(\pBox)$ 
has one divisorial stratum $D_{\Box}$, corresponding to 
$V_{\y}=V_{\yp}$.\notation{$\pBox,\x,\y,\yp,\z,V_{\x},V_{\y},V_{\yp},V_{\z}, D_\Box$}  
This variety will play a central role in the proof.

It will be useful to describe $X_{\bullet} \cup X_{\bn}$ in
terms of $\PF(\pBox)$ and $D_{\Box}$.
We do this by way of  a variety $X_{\bB}$
which is a certain open subset of $\PF(p_{\bB}: \poset_\bullet \rightarrow \pBox)$
(defined in Sect.~\ref{firstdesc})
whose image in $\PF(\poset_\bullet)$ is  $\Xbb$,
where  $p_{\bullet \Box}$ is the
morphism of posets depicted in
Figure~\ref{bronx}.\notation{$p_{\bullet \Box}$}

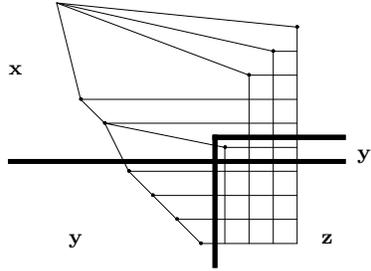
\begin{figure}[ht]
\begin{center} \setlength{\unitlength}{0.00083333in}
\begingroup\makeatletter\ifx\SetFigFont\undefined%
\gdef\SetFigFont#1#2#3#4#5{%
  \reset@font\fontsize{#1}{#2pt}%
  \fontfamily{#3}\fontseries{#4}\fontshape{#5}%
  \selectfont}%
\fi\endgroup%
{\renewcommand{\dashlinestretch}{30}
\begin{picture}(2187,1689)(0,-10)
\put(1812,1512){\blacken\ellipse{20}{20}}
\put(1812,1512){\ellipse{20}{20}}
\put(1662,1362){\blacken\ellipse{20}{20}}
\put(1662,1362){\ellipse{20}{20}}
\put(1512,1212){\blacken\ellipse{20}{20}}
\put(1512,1212){\ellipse{20}{20}}
\put(462,1062){\blacken\ellipse{20}{20}}
\put(462,1062){\ellipse{20}{20}}
\put(612,912){\blacken\ellipse{20}{20}}
\put(612,912){\ellipse{20}{20}}
\put(1362,762){\blacken\ellipse{20}{20}}
\put(1362,762){\ellipse{20}{20}}
\put(762,612){\blacken\ellipse{20}{20}}
\put(762,612){\ellipse{20}{20}}
\put(912,462){\blacken\ellipse{20}{20}}
\put(912,462){\ellipse{20}{20}}
\put(1062,312){\blacken\ellipse{20}{20}}
\put(1062,312){\ellipse{20}{20}}
\put(1212,162){\blacken\ellipse{20}{20}}
\put(1212,162){\ellipse{20}{20}}
\path(312,1662)(462,1062)(1812,1062)
\path(312,1662)(1812,1512)(1812,162)
	(1212,162)(762,612)(612,912)(462,1062)
\path(1812,1362)(1662,1362)(1662,162)
\path(312,1662)(1662,1362)
\path(312,1662)(1512,1212)
\path(1812,1212)(1512,1212)(1512,162)
\path(1812,912)(612,912)(1362,762)(1362,162)
\path(1812,762)(1362,762)
\path(762,612)(1812,612)
\path(912,462)(1812,462)
\path(1062,312)(1812,312)
\blacken\path(1287,837)(2112,837)(2112,812)
	(1287,812)(1287,837)
\path(1287,837)(2112,837)(2112,812)
	(1287,812)(1287,837)
\blacken\path(1287,837)(1312,837)(1312,12)
	(1287,12)(1287,837)
\path(1287,837)(1312,837)(1312,12)
	(1287,12)(1287,837)
\blacken\path(2112,662)(12,662)(12,687)
	(2112,687)(2112,662)
\path(2112,662)(12,662)(12,687)
	(2112,687)(2112,662)
\put(387,162){\makebox(0,0)[lb]{\smash{{{\SetFigFont{8}{9.6}{\rmdefault}{\mddefault}{\updefault}$\mathbf{y}$}}}}}
\put(1962,162){\makebox(0,0)[lb]{\smash{{{\SetFigFont{8}{9.6}{\rmdefault}{\mddefault}{\updefault}$\mathbf{z}$}}}}}
\put(12,1212){\makebox(0,0)[lb]{\smash{{{\SetFigFont{8}{9.6}{\rmdefault}{\mddefault}{\updefault}$\mathbf{x}$}}}}}
\put(2187,687){\makebox(0,0)[lb]{\smash{{{\SetFigFont{8}{9.6}{\rmdefault}{\mddefault}{\updefault}$\mathbf{y'}$}}}}}
\end{picture}
} \end{center}
\caption{A pictorial depiction of $p_{\bullet \Box}:  \poset_{\bullet} \rightarrow \pBox$;
elements of $p_{\bullet\Box}^{-1}(\x)$,
$p_{\bullet \Box}^{-1}(\y)$,
$p_{\bullet \Box}^{-1}(\yp)$,
$p_{\bullet \Box}^{-1}(\z)$ lie in regions labeled $\x$, $\y$, $\yp$, $\z$ respectively.  
\lremind{bronx}}
\label{bronx}
\end{figure}

\epoint{First description of $X_{\bB}$} \label{firstdesc}
\lremind{firstdesc}
Suppose the descending checker is at $(r,c)$.\notation{first use of
  $(r,c)$} Let \lremind{merced1--2}
\begin{equation}
\label{merced1}
S_{c-1} \subset S_c
\subset \dots \subset S_n
\end{equation} 
be the subspaces corresponding to elements of $\poset_{\bullet}$
in the bottom row of the table (part of the flag $\ff_{\bu}$).  
(Subscripts denote dimension.)
Let
\begin{equation}
\label{merced2}
T_0 \subset T_1 \subset \dots \subset
T_n
\end{equation}
be the subspaces corresponding to the northeast border of
$\poset_\bullet$ (part of the flag $\fm_{\bu}$).\notation{$S_{\bu}$, $T_{\bu}$}
Define
$X_{\bB}$\notation{$X_{\bB}$} to be the locally  closed subvariety of $\left( \Xbb \right)
\times \PF(\pBox)$ such that\lremind{carmel}
\begin{equation}
\label{carmel}V_{\x} \cap T_{r+1} = T_{r-1}, \quad V_{\yp} \cap T_{r+1} = T_r, \quad V_{\y} \cap S_c = S_{c-1}.
\end{equation}
Then the vector spaces corresponding to elements of $\poset_\bullet$
are determined by $V_{\x}$, $V_{\y}$, $V_{\yp}$, and \lremind{neworleans}
\begin{equation}
\label{neworleans}
S_c \subset S_{c+1} \subset \cdots \subset S_n, \quad
T_0 \subset T_1 \subset \dots \subset T_{r-2} \subset T_{r+1}
\subset \dots \subset T_{n-1} \subset T_n
\end{equation}
as shown in Figure~\ref{boston3}.
(This is because each element of $\poset_\bullet$ and $\poset_{\bn}$
is of the form $\inf(\mathbf{b}, \mathbf{b'})$, for some  $\mathbf{b'}$ on the northeast border,
and $\mathbf{b}$ in the bottom row of the table;
and $V_{\inf(\mathbf{b},\mathbf{b'})} = V_{\mathbf{b}} \cap V_{\mathbf{b'}}$ for any point of $X_\bullet$
and $X_{\bn}$.  Hence $V_{\inf(\mathbf{b},\mathbf{b'})}$ is determined by \eqref{merced1}
and \eqref{merced2}, and thus \eqref{neworleans} and $V_{\x}$, $V_{\y}$,
$V_{\yp}$.)

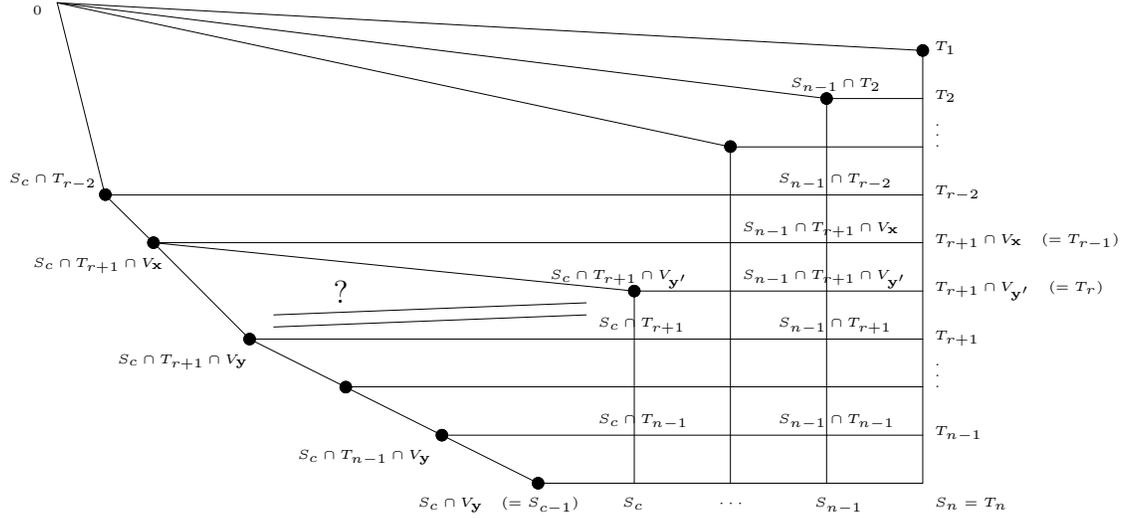
\begin{figure}[ht]
\begin{center} \setlength{\unitlength}{0.00083333in}
\begingroup\makeatletter\ifx\SetFigFont\undefined%
\gdef\SetFigFont#1#2#3#4#5{%
  \reset@font\fontsize{#1}{#2pt}%
  \fontfamily{#3}\fontseries{#4}\fontshape{#5}%
  \selectfont}%
\fi\endgroup%
{\renewcommand{\dashlinestretch}{30}
\begin{picture}(5775,3213)(0,-10)
\put(5700,2886){\blacken\ellipse{74}{74}}
\put(5700,2886){\ellipse{74}{74}}
\put(5100,2586){\blacken\ellipse{74}{74}}
\put(5100,2586){\ellipse{74}{74}}
\put(4500,2286){\blacken\ellipse{74}{74}}
\put(4500,2286){\ellipse{74}{74}}
\put(3900,1386){\blacken\ellipse{74}{74}}
\put(3900,1386){\ellipse{74}{74}}
\put(3300,186){\blacken\ellipse{74}{74}}
\put(3300,186){\ellipse{74}{74}}
\put(2700,486){\blacken\ellipse{74}{74}}
\put(2700,486){\ellipse{74}{74}}
\put(2100,786){\blacken\ellipse{74}{74}}
\put(2100,786){\ellipse{74}{74}}
\put(900,1686){\blacken\ellipse{74}{74}}
\put(900,1686){\ellipse{74}{74}}
\put(1500,1086){\blacken\ellipse{74}{74}}
\put(1500,1086){\ellipse{74}{74}}
\put(600,1986){\blacken\ellipse{74}{74}}
\put(600,1986){\ellipse{74}{74}}
\path(300,3186)(5100,2586)
\path(300,3186)(4500,2286)
\path(5700,2286)(4500,2286)(4500,186)
\path(2100,786)(5700,786)
\path(2700,486)(5700,486)
\path(300,3186)(5700,2886)(5700,186)
	(3300,186)(1500,1086)(900,1686)(600,1986)
\path(5700,2586)(5100,2586)(5100,186)
\path(300,3186)(600,1986)(5700,1986)
\path(5700,1686)(900,1686)(3900,1386)(3900,186)
\path(5700,1386)(3900,1386)
\path(1500,1086)(5700,1086)
\path(1650,1236)(3600,1311)
\path(1650,1161)(3600,1236)
\put(5775,2886){\makebox(0,0)[lb]{\smash{{{\SetFigFont{5}{6.0}{\rmdefault}{\mddefault}{\updefault}$T_1$}}}}}
\put(5775,2586){\makebox(0,0)[lb]{\smash{{{\SetFigFont{5}{6.0}{\rmdefault}{\mddefault}{\updefault}$T_2$}}}}}
\put(5775,2286){\makebox(0,0)[lb]{\smash{{{\SetFigFont{5}{6.0}{\rmdefault}{\mddefault}{\updefault}$\vdots$}}}}}
\put(5775,1986){\makebox(0,0)[lb]{\smash{{{\SetFigFont{5}{6.0}{\rmdefault}{\mddefault}{\updefault}$T_{r-2}$}}}}}
\put(5775,1686){\makebox(0,0)[lb]{\smash{{{\SetFigFont{5}{6.0}{\rmdefault}{\mddefault}{\updefault}$T_{r+1} \cap V_{\mathbf{x}} \quad (=T_{r-1})$}}}}}
\put(5775,1386){\makebox(0,0)[lb]{\smash{{{\SetFigFont{5}{6.0}{\rmdefault}{\mddefault}{\updefault}$T_{r+1} \cap V_{\mathbf{y'}} \quad (=T_r)$}}}}}
\put(5775,1086){\makebox(0,0)[lb]{\smash{{{\SetFigFont{5}{6.0}{\rmdefault}{\mddefault}{\updefault}$T_{r+1}$}}}}}
\put(5775,786){\makebox(0,0)[lb]{\smash{{{\SetFigFont{5}{6.0}{\rmdefault}{\mddefault}{\updefault}$\vdots$}}}}}
\put(5775,486){\makebox(0,0)[lb]{\smash{{{\SetFigFont{5}{6.0}{\rmdefault}{\mddefault}{\updefault}$T_{n-1}$}}}}}
\put(5775,36){\makebox(0,0)[lb]{\smash{{{\SetFigFont{5}{6.0}{\rmdefault}{\mddefault}{\updefault}$S_n=T_n$}}}}}
\put(4875,2661){\makebox(0,0)[lb]{\smash{{{\SetFigFont{5}{6.0}{\rmdefault}{\mddefault}{\updefault}$S_{n-1} \cap T_2$}}}}}
\put(4575,1761){\makebox(0,0)[lb]{\smash{{{\SetFigFont{5}{6.0}{\rmdefault}{\mddefault}{\updefault}$S_{n-1} \cap T_{r+1} \cap V_{\mathbf{x}}$}}}}}
\put(3675,1161){\makebox(0,0)[lb]{\smash{{{\SetFigFont{5}{6.0}{\rmdefault}{\mddefault}{\updefault}$S_{c} \cap T_{r+1}$}}}}}
\put(3675,561){\makebox(0,0)[lb]{\smash{{{\SetFigFont{5}{6.0}{\rmdefault}{\mddefault}{\updefault}$S_{c} \cap T_{n-1}$}}}}}
\put(5025,36){\makebox(0,0)[lb]{\smash{{{\SetFigFont{5}{6.0}{\rmdefault}{\mddefault}{\updefault}$S_{n-1}$}}}}}
\put(3825,36){\makebox(0,0)[lb]{\smash{{{\SetFigFont{5}{6.0}{\rmdefault}{\mddefault}{\updefault}$S_c$}}}}}
\put(4425,36){\makebox(0,0)[lb]{\smash{{{\SetFigFont{5}{6.0}{\rmdefault}{\mddefault}{\updefault}$\cdots$}}}}}
\put(150,1536){\makebox(0,0)[lb]{\smash{{{\SetFigFont{5}{6.0}{\rmdefault}{\mddefault}{\updefault}$S_c \cap T_{r+1} \cap V_{\mathbf{x}}$}}}}}
\put(150,3111){\makebox(0,0)[lb]{\smash{{{\SetFigFont{5}{6.0}{\rmdefault}{\mddefault}{\updefault}$0$}}}}}
\put(675,936){\makebox(0,0)[lb]{\smash{{{\SetFigFont{5}{6.0}{\rmdefault}{\mddefault}{\updefault}$S_c \cap T_{r+1} \cap V_{\mathbf{y}}$}}}}}
\put(1800,336){\makebox(0,0)[lb]{\smash{{{\SetFigFont{5}{6.0}{\rmdefault}{\mddefault}{\updefault}$S_c \cap T_{n-1} \cap V_{\mathbf{y}}$}}}}}
\put(4800,2061){\makebox(0,0)[lb]{\smash{{{\SetFigFont{5}{6.0}{\rmdefault}{\mddefault}{\updefault}$S_{n-1} \cap T_{r-2}$}}}}}
\put(4575,1461){\makebox(0,0)[lb]{\smash{{{\SetFigFont{5}{6.0}{\rmdefault}{\mddefault}{\updefault}$S_{n-1} \cap T_{r+1} \cap V_{\mathbf{y'}}$}}}}}
\put(4800,1161){\makebox(0,0)[lb]{\smash{{{\SetFigFont{5}{6.0}{\rmdefault}{\mddefault}{\updefault}$S_{n-1} \cap T_{r+1}$}}}}}
\put(4800,561){\makebox(0,0)[lb]{\smash{{{\SetFigFont{5}{6.0}{\rmdefault}{\mddefault}{\updefault}$S_{n-1} \cap T_{n-1}$}}}}}
\put(3375,1461){\makebox(0,0)[lb]{\smash{{{\SetFigFont{5}{6.0}{\rmdefault}{\mddefault}{\updefault}$S_c \cap T_{r+1} \cap V_{\mathbf{y'}}$}}}}}
\put(2025,1311){\makebox(0,0)[lb]{\smash{{{\SetFigFont{12}{14.4}{\rmdefault}{\mddefault}{\updefault}?}}}}}
\put(0,2061){\makebox(0,0)[lb]{\smash{{{\SetFigFont{5}{6.0}{\rmdefault}{\mddefault}{\updefault}$S_{c} \cap T_{r-2}$}}}}}
\put(2550,36){\makebox(0,0)[lb]{\smash{{{\SetFigFont{5}{6.0}{\rmdefault}{\mddefault}{\updefault}$S_c \cap V_{\mathbf{y}} \quad (=S_{c-1})$}}}}}
\end{picture}
} \end{center}
\caption{Pictorial description of the morphism $X_{\bB} \rightarrow 
X_{\bullet} \cup X_{\bn}$ (in the guise of the morphism $X_{\bB} \rightarrow \PF(\poset_{\bullet})$) \lremind{boston3}}
\label{boston3}
\end{figure}

\epoint{Second description of $X_{\bB}$}  
\label{secdesc} \lremind{secdesc} Equivalently, $X_{\bB}$
parametrizes
$V_{\x}$, $V_{\y}$, $V_{\yp}$,
  and
 the  partial flags $S_{\bu}$ and $T_{\bu}$ of \eqref{neworleans}
({\em not} \eqref{merced1} and \eqref{merced2}),
such that 
\begin{itemize}
\item the partial flags $S_{\bu}$ and $T_{\bu}$ are transverse,
\item $T_{r-2} \subset T_{r+1} \cap V_{\x}$, 
\item $V_{\x} \subset V_{\y}, V_{\yp}$, and
\item $V_{\x}$ is transverse to $S_c \cap T_{r+1}$.
(As $r+c \geq n+1$, $\dim S_c \cap T_{r+1} \geq 2$,  so $V_{\y}$
and $V_{\yp}$ are also transverse to $S_c \cap T_{r+1}$.)
\end{itemize}
If these conditions hold, we say
that $V_{\x}$, $V_{\y}$, $V_{\yp}$, $S_{\bu}$, $T_{\bu}$,  are in
{\em $X_{\bB}$-position}.\notation{$X_{\bB}$-position}

The projection $f_{\bullet \bn}: X_{\bB} \rightarrow
X_{\bullet} \cup X_{\bn}$ is smooth
and surjective: given an element of $X_{\bullet} \cup X_{\bn}$,
where $S_{\bu}$ and $T_{\bu}$ are the partial flags of (\ref{merced1}) and
(\ref{merced2}),
the fiber 
 corresponds to a choice of 
$(V_{\x}, V_{\y}, V_{\yp})$ satisfying (\ref{carmel}).
More precisely, (i) let $V_{\y} = T_{r-1} + S_{c-1}$, (ii)
then choose $V_{\x}$ in $V_{\y}$ such that $V_{\x} \cap T_{r+1} = 
T_{r-1}$ (or equivalently such that $V_{\x} \supset T_{r-1}$), (iii) then let $V_{\yp} = V_{\x} + T_r$.

Similarly, the  morphism $f_{\Box}: X_{\bB} \rightarrow
\PF(\pBox)$ is  smooth and
surjective; by the definition of $X_{\bB}$-position,
$f_{\Box}$ expresses $X_{\bB}$ as an open subset of
a tower of projective bundles over $\PF(\pBox)$.

\point \lremind{alewife} \label{alewife} 
Note that $f_{\bullet \bn}^{-1} X_{\bn} = f_{\Box}^{-1}
D_{\Box}$ (see Figure~\ref{boston3} --- $V_{\y} = V_{\yp}$ iff $S_c \cap T_{r+1} \cap V_{\y} = S_c \cap T_{r+1} \cap V_{\yp}$).  This
will allow us to translate questions about codimension one
degenerations in the specialization order (involving $X_{\bn} \subset
X_{\bullet} \cup X_{\bn}$) to simpler facts about codimension one
degenerations relating to $D_{\Box} \subset \PF(\pBox)$.  
The notation $f_{\bullet \bn}$ and $f_{\Box}$ will not be used hereafter.

\section{Proof of the \LR}
\label{thmpf} \lremind{thmpf}
We now prove (the final version of) the \LR, Theorem~\ref{glr3}.

\bpoint{Strategy of proof} \label{strategy} \lremind{strategy}
The strategy is as follows.  Instead of considering the ``divisor
at $\infty$'' of the closure of $X_{\cb}$ in $G(k,n) \times 
\left( \Xbb \right)$, we consider the corresponding divisor
on closures of other sets in different sets in different spaces,
shown in (\ref{newton}).  (The varieties $X_{\cbB}$ and 
$X_{\cB}$ will be defined in Sect.~\ref{providence},
and the smoothness and surjectivity of the morphisms of the
top row will be established.
For convenience, let $Z$
 be $\Cl_{\PF(\poset_\circ) \times \PF(\pBox)} X_{\cB}$.)
There is a ``divisor at $\infty$'' on each of these varieties
that behaves well with respect to pullback by these morphisms;
it corresponds to $X_{\bn} \subset \Xbb$ on the left
and $D_\Box \subset \PF(\pBox)$ on the right.
\lremind{newton}
\begin{equation}\label{newton}
\begin{array} {ccccc}
\Cl_{\PF(\poset_\circ) \times \left( \Xbb \right) }  X_{\cb} &
\stackrel {\text{sm. surj.}} \longleftarrow & 
\Cl_{\PF(\poset_\circ) \times \left( \Xbb \right) } X_{\circ \bB} &
\stackrel {\text{sm. surj.}} \longrightarrow & 
Z := \Cl_{\PF(\poset_\circ) \times \PF(\pBox)} X_{\cB} \\
\downarrow \text{\scriptsize {proper birat'l}} & & & & \downarrow \text{\scriptsize {closed imm.}} \\
\Cl_{G(k,n) \times \left( \Xbb \right)}  X_{\cb} & & & & 
\PF(\poset_\circ \rightarrow \pBox)
\end{array}
\end{equation}

We first identify the components of $\Cl_{\PF(\poset_\circ) \times
  \left( \Xbb \right) } X_{\cb}$ as follows.  The components of ``the
divisor at $\infty$'' of $\PF(\poset_\circ \rightarrow \pBox)$ of
dimension at least $\dim X_{\cB}-1$ are identified by
Theorem~\ref{detroit}.  They are denoted $D^Z_S$, where $S$ is a
certain set of ``good quadrilaterals'' in the poset $\poset_{\circ}$.
This gives a list containing the components of the ``divisor
at $\infty$'' of
$Z = \Cl_{\PF(\poset_\circ) \times \PF(\pBox)} X_{\cB}$.  Via
the smooth surjective morphisms of the top row of (\ref{newton}),
this immediately gives a list containing the components, denoted $D_S$, of 
the ``divisor at $\infty$'' of 
$\Cl_{\PF(\poset_\circ) \times \left( \Xbb \right) }  X_{\cb}$.

We say $D_S$ is {\em geometrically irrelevant}\notation{geometrically
  (ir)relevant} if its image in $G(k,n) \times X_{\bn}$ is of
smaller dimension (and {\em geometrically relevant} otherwise).  In
Section~\ref{capecod}, it is shown that all but one or two $D_S$
(corresponding to those described in the \LR) are geometrically irrelevant.  We do this
by exhibiting a one-parameter family through a general point of $D_S$
contracted by the morphism to $G(k,n) \times X_{\bn}$.

Next, in Section~\ref{amherst}, we show that in the one
or two geometrically relevant cases $D_S$ appears with multiplicity
one in the ``divisor at $\infty$'', by showing that $D^Z_S$ appears
with multiplicity one in the corresponding divisor on $Z$.

Finally, these one or two $D_S$'s map are birational to (and hence map
with degree 1 to) $X_{\cbn}$ or $X_{\cbs}$
(Proposition~\ref{birational}).  

\epoint{Remark} The bijection to puzzles
of Section~\ref{checkerpuzzle} gives a second
proof; here is a very quick sketch. 
Fix a mid-sort configuration $\cb$.
 First show that $D_X$ contains
both $X_{\cn \bn}$ and/or $X_{\cs \bn}$
(Sect.~\ref{amherst}). 
By an easy induction, $\cb$ arises in the course of a checker
game starting with some subsets $\al$ and $\be$.  Then $D_X$ can have no
more components, and these one or two must appear with multiplicity
one.  Otherwise, choose any $\ga$ such that the non-zero effective
cycle $D_X - \overline{X}_{\cbn}$,
$D_X - \overline{X}_{\cbs}$, or $D_X - \overline{X}_{\cbn}-
\overline{X}_{\cbs}$ has a non-zero coefficient of the basis element
corresponding to $\ga$.  Then the number of puzzles with inputs $\al$
and $\be$ and output $\ga$ is strictly less than the
Littlewood-Richardson coefficient $c^{\ga}_{\al \be}$, giving a
contradiction.

One advantage of the proof presented in Section~\ref{thmpf} is that in
order to generalize the \LR{} to other geometric situations (see for
example \cite{si} 
and Conjecture~\ref{bvconjecture}), one needs the geometry behind it.

\bpoint{Reduction to $\PF(\pBox)$} 
\label{providence} 
\lremind{providence}
We first reduce much of the argument to
statements involving the poset $\pBox$ (Figure~\ref{cupertino}) rather than the
more complicated $\poset_{\bullet}$
as follows.  

The composition of  $p_{\cb}: \poset_{\circ} \rightarrow
\poset_{\bullet}$ and $p_{\bullet \Box}: \poset_{\bullet} \rightarrow
\pBox$ is the morphism $\pcb: \poset_{\circ} \rightarrow
\pBox$ depicted in Figure~\ref{dc}.
Let
$X_{\cB}$\notation{$\pcb$, $X_{\cB}$} be the open subset of $\PF(\pcb:
\poset_{\circ} \rightarrow \pBox)$ such that if $\mathbf{b} \in \poset_{\circ}$
and $\mathbf{c} \in \pBox$, then $V_{\mathbf{b}} \subset V_{\mathbf{c}}$ if {\em and only if}
$\pcb(\mathbf{b}) \prec \mathbf{c}$.
(One could call this
a ``precise Bott-Samelson \Qflag variety'' of the morphism
of posets, $P\PF(p_{\cB})$.  Then $X_{\cb} = P\PF(p_{\cb})|_{X_{\bullet}}$,
$X_{\bB}= P\PF(p_{\bB})|_{\Xbb}$, and $P\PF(p) \subset \PF(p)$ is
an open immersion for
all $p$.  However, to minimize notation, we will not use this terminology.)

\begin{figure}[ht]
\begin{center} \setlength{\unitlength}{0.00083333in}
\begingroup\makeatletter\ifx\SetFigFont\undefined%
\gdef\SetFigFont#1#2#3#4#5{%
  \reset@font\fontsize{#1}{#2pt}%
  \fontfamily{#3}\fontseries{#4}\fontshape{#5}%
  \selectfont}%
\fi\endgroup%
{\renewcommand{\dashlinestretch}{30}
\begin{picture}(3825,3130)(0,-10)
\texture{88555555 55000000 555555 55000000 555555 55000000 555555 55000000 
	555555 55000000 555555 55000000 555555 55000000 555555 55000000 
	555555 55000000 555555 55000000 555555 55000000 555555 55000000 
	555555 55000000 555555 55000000 555555 55000000 555555 55000000 }
\shade\path(3375,1456)(3525,1456)(3525,1306)
	(3375,1306)(3375,1456)
\path(3375,1456)(3525,1456)(3525,1306)
	(3375,1306)(3375,1456)
\texture{0 0 0 888888 88000000 0 0 80808 
	8000000 0 0 888888 88000000 0 0 80808 
	8000000 0 0 888888 88000000 0 0 80808 
	8000000 0 0 888888 88000000 0 0 80808 }
\shade\path(3375,1156)(3525,1156)(3525,1006)
	(3375,1006)(3375,1156)
\path(3375,1156)(3525,1156)(3525,1006)
	(3375,1006)(3375,1156)
\put(3600,1006){\makebox(0,0)[lb]{\smash{{{\SetFigFont{8}{9.6}{\rmdefault}{\mddefault}{\updefault}$=$ right good quadrilaterals}}}}}
\put(3600,1306){\makebox(0,0)[lb]{\smash{{{\SetFigFont{8}{9.6}{\rmdefault}{\mddefault}{\updefault}$=$ left good quadrilaterals}}}}}
\put(1275,1456){\ellipse{50}{50}}
\put(1425,1306){\ellipse{50}{50}}
\put(1575,1006){\ellipse{50}{50}}
\put(1725,856){\ellipse{50}{50}}
\put(1875,1156){\ellipse{50}{50}}
\put(2025,1606){\ellipse{50}{50}}
\put(1125,2056){\ellipse{50}{50}}
\put(2175,1906){\ellipse{50}{50}}
\put(975,2506){\ellipse{50}{50}}
\put(2325,2206){\ellipse{50}{50}}
\put(2475,2356){\ellipse{50}{50}}
\path(1125,2056)(2475,2056)
\path(1125,2056)(2175,1906)(2475,1906)
\path(1125,2056)(2025,1606)(2475,1606)
\path(2175,1906)(2175,856)
\path(2025,1606)(2025,856)
\path(1275,1456)(2475,1456)
\path(1425,1306)(2475,1306)
\path(1425,1306)(1875,1156)(2475,1156)
\path(1575,1006)(2475,1006)
\path(1875,1156)(1875,856)
\path(3750,2206)(2475,1906)
\blacken\path(2584.939,1962.687)(2475.000,1906.000)(2598.681,1904.282)(2584.939,1962.687)
\path(3750,1981)(2475,1606)
\blacken\path(2581.659,1668.641)(2475.000,1606.000)(2598.589,1611.079)(2581.659,1668.641)
\path(679,1069)(1500,1006)
\blacken\path(1378.056,985.269)(1500.000,1006.000)(1382.647,1045.093)(1378.056,985.269)
\path(675,2806)(975,2506)(1125,2056)
	(1275,1456)(1425,1306)(1575,1006)
	(1725,856)(2475,856)(2475,2356)
	(975,2506)(2325,2206)(2325,856)
\dashline{60.000}(1575,931)(2550,931)(2550,781)
	(1575,781)(1575,931)
\dashline{60.000}(2550,2431)(2550,1756)(2025,1756)
	(525,2656)(525,2956)(825,2956)(2550,2431)
\path(375,1906)(1125,2056)
\blacken\path(1013.214,2003.049)(1125.000,2056.000)(1001.447,2061.883)(1013.214,2003.049)
\path(2400,181)(2100,781)
\blacken\path(2180.498,687.085)(2100.000,781.000)(2126.833,660.252)(2180.498,687.085)
\path(3150,3031)(2550,2431)
\blacken\path(2613.640,2537.066)(2550.000,2431.000)(2656.066,2494.640)(2613.640,2537.066)
\shade\path(2025,1606)(2475,1606)(2475,1006)
	(2025,1006)(2025,1606)
\path(2025,1606)(2475,1606)(2475,1006)
	(2025,1006)(2025,1606)
\path(2175,1606)(2175,1006)
\path(2175,1606)(2175,1006)
\path(2325,1606)(2325,1006)
\path(2325,1606)(2325,1006)
\path(2025,1456)(2475,1456)
\path(2025,1456)(2475,1456)
\path(2025,1306)(2475,1306)
\path(2025,1306)(2475,1306)
\path(2025,1156)(2475,1156)
\path(2025,1156)(2475,1156)
\texture{88555555 55000000 555555 55000000 555555 55000000 555555 55000000 
	555555 55000000 555555 55000000 555555 55000000 555555 55000000 
	555555 55000000 555555 55000000 555555 55000000 555555 55000000 
	555555 55000000 555555 55000000 555555 55000000 555555 55000000 }
\shade\path(1125,2056)(1275,1456)(1425,1306)
	(2025,1306)(2025,1606)(1125,2056)
\path(1125,2056)(1275,1456)(1425,1306)
	(2025,1306)(2025,1606)(1125,2056)
\blacken\path(675,1531)(2775,1531)(2775,1506)
	(675,1506)(675,1531)
\path(675,1531)(2775,1531)(2775,1506)
	(675,1506)(675,1531)
\blacken\path(1650,556)(1675,556)(1675,1681)
	(1650,1681)(1650,556)
\path(1650,556)(1675,556)(1675,1681)
	(1650,1681)(1650,556)
\blacken\path(1650,1681)(2775,1681)(2775,1656)
	(1650,1656)(1650,1681)
\path(1650,1681)(2775,1681)(2775,1656)
	(1650,1656)(1650,1681)
\path(1275,1456)(2025,1456)
\put(600,706){\makebox(0,0)[lb]{\smash{{{\SetFigFont{8}{9.6}{\rmdefault}{\mddefault}{\updefault}$\mathbf{y}$}}}}}
\put(2925,1606){\makebox(0,0)[lb]{\smash{{{\SetFigFont{8}{9.6}{\rmdefault}{\mddefault}{\updefault}$\mathbf{y'}$}}}}}
\put(2925,706){\makebox(0,0)[lb]{\smash{{{\SetFigFont{8}{9.6}{\rmdefault}{\mddefault}{\updefault}$\mathbf{z}$}}}}}
\put(3825,2206){\makebox(0,0)[lb]{\smash{{{\SetFigFont{8}{9.6}{\rmdefault}{\mddefault}{\updefault}$\mathbf{a''}$}}}}}
\put(3825,1906){\makebox(0,0)[lb]{\smash{{{\SetFigFont{8}{9.6}{\rmdefault}{\mddefault}{\updefault}$\mathbf{a'}$}}}}}
\put(0,2206){\makebox(0,0)[lb]{\smash{{{\SetFigFont{8}{9.6}{\rmdefault}{\mddefault}{\updefault}$\mathbf{x}$}}}}}
\put(0,1756){\makebox(0,0)[lb]{\smash{{{\SetFigFont{8}{9.6}{\rmdefault}{\mddefault}{\updefault}$\inf(\mathbf{a},\mathbf{a'})$}}}}}
\put(3150,3031){\makebox(0,0)[lb]{\smash{{{\SetFigFont{8}{9.6}{\rmdefault}{\mddefault}{\updefault}Remove these to obtain $\mathcal{P}'_{\circ}$}}}}}
\put(2400,31){\makebox(0,0)[lb]{\smash{{{\SetFigFont{8}{9.6}{\rmdefault}{\mddefault}{\updefault}Remove these to obtain $\mathcal{P}'_{\circ}$}}}}}
\put(375,1081){\makebox(0,0)[lb]{\smash{{{\SetFigFont{8}{9.6}{\rmdefault}{\mddefault}{\updefault}$\mathbf{a}$}}}}}
\put(825,2806){\makebox(0,0)[lb]{\smash{{{\SetFigFont{8}{9.6}{\rmdefault}{\mddefault}{\updefault}$0$}}}}}
\end{picture}
} \end{center}
\caption{$p_{\cB}:  \poset_{\circ} \rightarrow \pBox$, $\mathbf{a}$, $\mathbf{a'}$, ${\mathbf{a''}}$, $\poset'_{\circ}$
\lremind{dc}}
\label{dc}
\end{figure}

\begin{figure}
\begin{center} \setlength{\unitlength}{0.00083333in}
\begingroup\makeatletter\ifx\SetFigFont\undefined%
\gdef\SetFigFont#1#2#3#4#5{%
  \reset@font\fontsize{#1}{#2pt}%
  \fontfamily{#3}\fontseries{#4}\fontshape{#5}%
  \selectfont}%
\fi\endgroup%
{\renewcommand{\dashlinestretch}{30}
\begin{picture}(1452,1935)(0,-10)
\put(1200,1611){\makebox(0,0)[lb]{\smash{{{\SetFigFont{5}{6.0}{\rmdefault}{\mddefault}{\updefault}smooth}}}}}
\put(1200,1536){\makebox(0,0)[lb]{\smash{{{\SetFigFont{5}{6.0}{\rmdefault}{\mddefault}{\updefault}surj.}}}}}
\put(225,1611){\makebox(0,0)[lb]{\smash{{{\SetFigFont{5}{6.0}{\rmdefault}{\mddefault}{\updefault}smooth}}}}}
\put(225,1536){\makebox(0,0)[lb]{\smash{{{\SetFigFont{5}{6.0}{\rmdefault}{\mddefault}{\updefault}surj.}}}}}
\path(787,1761)(787,1086)
\blacken\path(772.000,1146.000)(787.000,1086.000)(802.000,1146.000)(772.000,1146.000)
\path(712,1761)(187,1386)
\blacken\path(227.105,1433.080)(187.000,1386.000)(244.543,1408.668)(227.105,1433.080)
\path(862,1761)(1387,1386)
\blacken\path(1329.457,1408.668)(1387.000,1386.000)(1346.895,1433.080)(1329.457,1408.668)
\path(1425,1161)(1425,186)
\blacken\path(1410.000,246.000)(1425.000,186.000)(1440.000,246.000)(1410.000,246.000)
\path(150,1161)(150,186)
\blacken\path(135.000,246.000)(150.000,186.000)(165.000,246.000)(135.000,246.000)
\path(750,861)(225,186)
\blacken\path(249.996,242.570)(225.000,186.000)(273.677,224.152)(249.996,242.570)
\path(225,1161)(750,486)
\blacken\path(701.323,524.152)(750.000,486.000)(725.004,542.570)(701.323,524.152)
\path(1350,1161)(825,486)
\blacken\path(849.996,542.570)(825.000,486.000)(873.677,524.152)(849.996,542.570)
\path(825,861)(1350,186)
\blacken\path(1301.323,224.152)(1350.000,186.000)(1325.004,242.570)(1301.323,224.152)
\put(112,1236){\makebox(0,0)[lb]{\smash{{{\SetFigFont{8}{9.6}{\rmdefault}{\mddefault}{\updefault}$X_{\circ \bullet}$}}}}}
\put(1275,1236){\makebox(0,0)[lb]{\smash{{{\SetFigFont{8}{9.6}{\rmdefault}{\mddefault}{\updefault}$X_{\circ \Box}$}}}}}
\put(675,1836){\makebox(0,0)[lb]{\smash{{{\SetFigFont{8}{9.6}{\rmdefault}{\mddefault}{\updefault}$X_{\circ \bullet \Box}$}}}}}
\put(675,936){\makebox(0,0)[lb]{\smash{{{\SetFigFont{8}{9.6}{\rmdefault}{\mddefault}{\updefault}$X_{\bullet \Box}$}}}}}
\put(637,336){\makebox(0,0)[lb]{\smash{{{\SetFigFont{8}{9.6}{\rmdefault}{\mddefault}{\updefault}$\PF(\mathcal{P}_{\circ})$}}}}}
\put(1237,36){\makebox(0,0)[lb]{\smash{{{\SetFigFont{8}{9.6}{\rmdefault}{\mddefault}{\updefault}$\PF(\mathcal{P}_{\Box})$}}}}}
\put(262,523){\makebox(0,0)[lb]{\smash{{{\SetFigFont{5}{6.0}{\rmdefault}{\mddefault}{\updefault}smooth}}}}}
\put(262,411){\makebox(0,0)[lb]{\smash{{{\SetFigFont{5}{6.0}{\rmdefault}{\mddefault}{\updefault}surj.}}}}}
\put(1012,523){\makebox(0,0)[lb]{\smash{{{\SetFigFont{5}{6.0}{\rmdefault}{\mddefault}{\updefault}smooth}}}}}
\put(0,36){\makebox(0,0)[lb]{\smash{{{\SetFigFont{8}{9.6}{\rmdefault}{\mddefault}{\updefault}$X_{\bullet} \cup X_{\bn}$}}}}}
\put(1012,448){\makebox(0,0)[lb]{\smash{{{\SetFigFont{5}{6.0}{\rmdefault}{\mddefault}{\updefault}surj.}}}}}
\end{picture}
} \end{center}
\caption{The definition of $X_{\cbB}$  
\lremind{backbay}}
\label{backbay}
\end{figure}

Define $X_{\cbB}$ as a ``triple fibered product'' as in
Figure~\ref{backbay}, i.e. as the indirect limit of the lower part of 
the figure. Equivalently, 
$X_{\cbB}$ is the open subset of $X_{\circ \bullet}
\times_{\PF(\poset_{\circ})} X_{\cB}$ whose image in
$\PF(\poset_{\bullet}) \times \PF(\pBox)$ is in $X_{\bB}$. 
Note that
$X_{\cbB}$ is an open subset of both $X_{\cb} \times_{\Xbb} X_{\bB}$
and $X_{\cB} \times_{\PF(\pBox)} X_{\bB}$ (although not in general
equal to either), so the projections from $X_{\cbB}$ to $X_{\cb}$ and
$X_{\cB}$ are both smooth (by base change from
$X_{\bB} \rightarrow \Xbb$ and $X_{\bB} \rightarrow \PF(\pBox)$).  In addition, $X_{\cbB} \rightarrow
X_{\cb}$ is surjective; the fiber corresponds to a certain (non-empty,
open) choice of $V_{\x}$, $V_{\y}$, $V_{\yp}$.  
Similarly, $X_{\cbB} \rightarrow X_{\cB}$ is surjective.  This information is summarized in the
top and bottom rows of (\ref{miami}).\lremind{miami} 
\begin{equation}
\label{miami}
\begin{CD}
X_{\cb} 
@<\text{smooth}<\text{surjective}<
X_{\cbB}
@>\text{smooth}>\text{surjective}>
X_{\cB}
\\
@V\text{open}V\text{imm.}V
@V\text{open}V\text{imm.}V
@V\text{open}V\text{imm.}V
\\
\Cl_{\PF(\poset_{\circ}) \times \left( \Xbb \right)} X_{\cb} 
@<\text{smooth}<\text{surjective}<
\Cl_{\PF(\poset_{\circ}) \times X_{\bB}} X_{\cbB}
@>\text{smooth}>\text{surjective}>
Z = \Cl_{\PF(\poset_{\circ}) \times \PF(\pBox)} X_{\cB}
\\
@V\text{proper}VV
@V\text{proper}VV
@V\text{proper}VV
\\
\Xbb
@<\text{smooth}<\text{surjective}<
X_{\bB}
@>\text{smooth}>\text{surjective}>
\PF(\pBox)
\end{CD}
\end{equation}
As remarked earlier, 
$\Cl_{\PF(\poset_{\circ}) \times \PF(\pBox)} X_{\cB}$ 
is denoted by $Z$.\notation{$Z$}
The
morphism $$\Cl_{\PF(\poset_{\circ}) \times X_{\bB}} X_{\cbB}
\rightarrow \Cl_{\PF(\poset_{\circ}) \times \left( \Xbb \right)} X_{\cb}$$ is smooth
and surjective for the same reason that $X_{\cbB} \rightarrow
X_{\cb}$ and $X_{\bB} \rightarrow \Xbb$ were: the fiber
corresponds to choosing $V_{\x}$, $V_{\y}$, $V_{\yp}$, satisfying the requirements of
$X_{\bB}$-position, except the transversality conditions with the
spaces corresponding to elements of $\poset_\bullet$, so this morphism is a
tower of projective bundles.  

\point \label{provincetown} \lremind{provincetown}
The morphism $\Cl_{\PF(\poset_{\circ})
  \times X_{\bB}} X_{\cbB} \rightarrow Z$ 
is smooth for the same reason that
$X_{\cbB} \rightarrow X_{\cB}$  and $X_{\bB} \rightarrow
\PF(\pBox)$ were.
More precisely, this morphism  is constructed
by  choosing the $S_{\bu}$'s and $T_{\bu}$'s in (\ref{neworleans}):
\begin{enumerate}
\item[(i)] inductively choosing $S_i$ (starting with $i=n$ and 
decrementing) so that $S_i$ contains (the subspace corresponding to)
the maximum element of
$\poset_{\circ}$ in column up to $i$, then
\item[(ii)] inductively choosing $T_i$ (starting with $i=n$ and
decrementing) so that $T_i$ contains 
(the subspace corresponding to) the maximum element of
$\poset_{\circ}$ in row up to $i$, then
\item[(iii)] discarding the closed subset not in $X_{\bB}$-position,
and defining $T_{r-1}$, $T_r$, $S_{c-1}$ by (\ref{carmel}).
\end{enumerate}
Note that (i) involves only the southwest border of $\poset_{\circ}$
(and not the element ${\mathbf{a}}$ shown in Figure~\ref{dc}),
and (ii) involves only the northeast border
(and not the elements ${\mathbf{a'}}$ and
${\mathbf{a''}}$ shown in Figure~\ref{dc}).

\bpoint{Proof of the \LR{} in the cases where $\pcb^{-1}(\y)$ or $\pcb^{-1}(\yp)$ is empty} This
corresponds to the cases where there are no white checkers in the
critical diagonal or the critical row, respectively (five of the nine
cases of Table~\ref{keywest}).  

In the case $\pcb^{-1}(\yp) = \{ \}$, the family $Z \rightarrow
\PF(\pBox)$ is pulled back from a family over $\PF(\pBox - \{ \yp \})$, the partial flag
parametrizing $V_{\x} \subset V_{\y} \subset K^n$.  In particular, the pullback
of $D_{\Box} = \{ V_{\y} = V_{\yp} \}$ to $Z$ consists of one component, appearing with
multiplicity 1.  The corresponding divisor on $\Cl_{\PF(\poset_{\circ})
  \times \left( \Xbb \right) } X_{\cb}$ also appears with multiplicity 1 by
the top row of (\ref{newton}) or the middle row of (\ref{miami}), and it is isomorphic to its image 
$\overline{X}_{\cbn}$
in
$\overline{X}_{\cb}$.

The same argument holds
for the case $\pcb^{-1}(\y) = \{ \}$, with the
roles of $\y$ and $\yp$ switched. 

\bpoint{Proof in the remaining cases} \label{pit} \lremind{pit}
For the rest of Section~\ref{thmpf}, we assume that both
$\pcb^{-1}(\y)$ (the critical diagonal) and $\pcb^{-1}(\yp)$ 
(the critical row) are non-empty.  The white
checkers are in mid-sort, by Lemma~\ref{midsort}.
Let $\mathbf{a}$ (resp. $\mathbf{a'}$, ${\mathbf{a''}}$) be the maximum\notation{$\mathbf{a}$,$\mathbf{a'}$, ${\mathbf{a''}}$}
of $\pcb^{-1}(\y)$ (resp. $\pcb^{-1}(\yp)$, $\pcb^{-1}(\x)$),
so $V_{\mathbf{a}}$, $V_{\mathbf{a'}}$, $V_{{\mathbf{a''}}}$ are the corresponding
subspaces.  (See Figure~\ref{dc}.)

In lieu of studying the components of the preimage of $D_{\Box}$ on
$Z$, we find the components of the pullback of $D_{\Box}$ to the (a priori
larger) variety $\PF(\poset_{\circ} \rightarrow \poset_{\Box})$, and
show that they all have dimension at most $\dim X_{\cB}-1$.  

The pullback of $D_{\Box}$ to $\PF(\poset_{\circ} \rightarrow \pBox)$
is best described in terms of the stratification of $\PF(\poset_{\circ})$.
Define the {\em good quadrilaterals}
of $\poset_{\circ}$ to be those quadrilaterals 
whose minimal element does not dominate $a$, and either (a)
whose minimal element dominates the white checker  $\min( p^{-1}_{\cB}(\yp))$ in the critical row
$\pcb^{-1} (\yp)$, or (b) with two vertices dominating 
$\min( p^{-1}_{\cB}(\yp))$,
and two vertices in $\pcb^{-1}( \{ \x, \y \})$.  
Call those in case (a) {\em right good quadrilaterals} and those
in case (b)
{\em left good quadrilaterals}; see
Figure~\ref{dc}.  Note that the good quadrilaterals appear in
columns,\notation{left and right good quadrilaterals, columns} and that, if there is a blocker,  there is no column to the left of
the white checker in the critical row $\min(\pcb^{-1} (\yp))$ (i.e. no 
left good quadrilaterals).

The following theorem describes the pullback of $D_{\Box}$
to $\PF(\poset_{\circ} \rightarrow \pBox)$.
\tpoint{Theorem} \label{detroit} \lremind{detroit} {\em The components
  of the pullback of $D_{\Box}$ to $\PF(\poset_{\circ} \rightarrow
  \pBox)$ of dimension at least $\dim X_{\cB} -1$ correspond to sets of good quadrilaterals, with at most
  one quadrilateral in each column, as follows.  If $S$ is such a set
  of good quadrilaterals, the corresponding component is the closure
  of the pullback to $\PF(\poset_{\circ} \rightarrow \pBox)$ of the open
  stratum on $\PF(\poset_{\circ})$ corresponding to $S$.  Its 
dimension is exactly 
$\dim X_{\cB} -1$. }

Let $D^Z_S$ be the component corresponding to $S$.  Let $D_S$
be the corresponding subvariety of $\PF(\poset_{\circ}) \times X_{\bn}$
(obtained via the top row of (\ref{newton}) or the middle row of (\ref{miami})).\notation{$D^Z_S, D_S$}

\bpf  
We note first that 
$X_{\cB}$
can be constructed 
by starting with the dense open stratum of  $\PF(\poset_{\circ})$, then choosing
$V_{\y}$ containing $V_{\mathbf{a}}$ and $V_{{\mathbf{a''}}}$ (giving a dimensional
contribution of 
$n - (\dim \mathbf{a} + \dim {\mathbf{a''}} - \dim \inf(\mathbf{a},\mathbf{a'}))$,
as $V_{\mathbf{a}} \cap V_{{\mathbf{a''}}} = V_{\inf(\mathbf{a},\mathbf{a'})}$ for an element of 
the dense open stratum), then choosing $V_{\yp}$ containing $V_{\mathbf{a'}}$ 
(giving a dimensional contribution of  $n- \dim \mathbf{a'}$).  
Then $V_{\x} = V_{\y} \cap V_{\yp}$ is determined.
Hence \lremind{knoxville}
\begin{equation} \label{knoxville}
\dim X_{\cB} = 2n - \dim \mathbf{a} - \dim \mathbf{a'} - \dim {\mathbf{a''}} + \dim \inf(\mathbf{a},\mathbf{a'}).
\end{equation}

Next, consider an irreducible component of the pullback of $D_{\Box}$
to $\PF(\poset_{\circ} \rightarrow \poset_{\Box})$ of dimension at least $\dim X_{\cB}-1$.  We will show
that its dimension is precisely this, and that the component is of the
form described in Theorem~\ref{detroit}.

Let $\ell=\dim (V_{\mathbf{a}} \cap V_{\mathbf{a'}}) -
\dim \inf(\mathbf{a},\mathbf{a'})$ for   
 a general point of this component.
(Note that $\ell \geq 0$, as $V_{\mathbf{a}} \cap V_{\mathbf{a'}}$ contains
$V_{\inf(\mathbf{a},\mathbf{a'})}$.)
Let $Q_\ell \subset \PF(\poset_{\circ})$ be the
  locus where $\dim (V \cap V') = \ell+ \dim\inf(\mathbf{a},\mathbf{a'})$.\notation{$Q_\ell$}
Then a dense open subset
$U$ of the component with a morphism to $Q_\ell$, and $U$ is an open set of 
the fibration over $Q_\ell$ parametrizing choices of $V_{\y}=V_{\yp}$
containing $V_{\mathbf{a}}$ and $V_{\mathbf{a'}}$ (giving a dimensional contribution of $n-(\dim
V_{\mathbf{a}} + \dim V_{\mathbf{a'}} - \dim (V_{\mathbf{a}} \cap V_{\mathbf{a'}}))$), and $V_{\x}$ contained 
in $V_{\y}=V_{\yp}$ containing $V_{{\mathbf{a''}}}$
(giving a dimensional contribution of $n-1- \dim V_{{\mathbf{a''}}}$ once $V_{\y} = V_{\yp}$ is
chosen).  Thus by Proposition~\ref{oakpark} below,
\begin{eqnarray*}
\dim U &\leq& \dim Q_\ell + (n- \dim \mathbf{a} - \dim \mathbf{a'} + (\ell+ \dim \inf(\mathbf{a},\mathbf{a'}))) + 
(n-1- \dim {\mathbf{a''}})\\
&=& \left(\dim Q_\ell + \ell \right)  + 2n - \dim \mathbf{a} - \dim \mathbf{a'} - \dim {\mathbf{a''}} + \dim \inf(\mathbf{a},\mathbf{a'})-1 \\
&\leq& \dim X_{\cB} - 1  \quad \quad \quad \text{(by (\ref{knoxville}))},
\end{eqnarray*}
and if equality holds, then the the component
is of the form desired (also by Proposition~\ref{oakpark}).
\epf

\tpoint{Proposition} \label{oakpark} \lremind{oakpark} {\em The
  irreducible components of $Q_\ell$ have codimension at least $\ell$
  in $\PF(\poset_{\circ})$, and those components of codimension exactly
  $\ell$ are strata of the form described in
  Theorem~\ref{detroit}.}

\bpf 
Suppose $\poset'_{\circ}$\notation{$\poset'_{\circ}$} is defined by removing
from $\poset_{\circ}$ all rows below $\mathbf{a}$, and all elements of
$p_{\cB}^{-1}(\x)$ except $\inf(\mathbf{a},\mathbf{a'})$ (see Figure~\ref{dc}).  As in the proof of Lemma~\ref{bssmooth},
the natural morphism
$\PF(\poset_{\circ}) \rightarrow \PF(\poset'_{\circ})$
is a tower of projective bundles , and the
components of $Q_\ell$ on $\PF(\poset_{\circ})$
are precisely the pullback of components of the analogous $Q_\ell$ on
$\PF(\poset'_{\circ})$.
Hence the Proposition follows from the analogous result for $\PF(\poset_{\circ}')$, Proposition~\ref{oakpark2}. \epf

More precisely, it suffices to prove the result for $\poset'_{\circ}$
of the following form.  Fix $r_1 < \dots < r_y$ and $c_1 < \dots <
c_x$ where $x, y >1$, and $r_1 < R_1 < \dots <
R_z < r_y$ and $c_2 > C_1 > \dots > C_z > c_1$ ($(R_i,C_i)$ will be
location of white checkers not contained in the grid $\{ (r_i,
c_j)\}$, for example a blocker).  Then $\poset'_\circ$ is the poset
generated by the set $\{ (r_1,c_i) \}_{1 \leq i \leq x } \cup \{
(r_i,c_1) \}_{1 \leq i \leq y } \cup \{ (R_i,C_i) \}_{1 \leq i \leq
  z}$ (see Sect.~\ref{huntsville}, and Figure~\ref{tuscaloosa} for
an example). 
Note that the critical row corresponds to $\{ (r_1, c_2), \dots,
(r_1,c_x) \}$, the critical diagonal corresponds to $\{ (r_2, c_1),
\dots, (r_y,c_1) \}$, $\mathbf{a}=(r_y,c_1)$, $\mathbf{a'} = (r_1,c_x)$, and $\inf(\mathbf{a},\mathbf{a'}) =
(r_1,c_1)$.

\begin{figure}[ht]
\begin{center} \setlength{\unitlength}{0.00083333in}
\begingroup\makeatletter\ifx\SetFigFont\undefined%
\gdef\SetFigFont#1#2#3#4#5{%
  \reset@font\fontsize{#1}{#2pt}%
  \fontfamily{#3}\fontseries{#4}\fontshape{#5}%
  \selectfont}%
\fi\endgroup%
{\renewcommand{\dashlinestretch}{30}
\begin{picture}(3525,2799)(0,-10)
\texture{88555555 55000000 555555 55000000 555555 55000000 555555 55000000 
	555555 55000000 555555 55000000 555555 55000000 555555 55000000 
	555555 55000000 555555 55000000 555555 55000000 555555 55000000 
	555555 55000000 555555 55000000 555555 55000000 555555 55000000 }
\shade\path(3300,1200)(3450,1200)(3450,1050)
	(3300,1050)(3300,1200)
\path(3300,1200)(3450,1200)(3450,1050)
	(3300,1050)(3300,1200)
\texture{0 0 0 888888 88000000 0 0 80808 
	8000000 0 0 888888 88000000 0 0 80808 
	8000000 0 0 888888 88000000 0 0 80808 
	8000000 0 0 888888 88000000 0 0 80808 }
\shade\path(3300,900)(3450,900)(3450,750)
	(3300,750)(3300,900)
\path(3300,900)(3450,900)(3450,750)
	(3300,750)(3300,900)
\put(3525,750){\makebox(0,0)[lb]{\smash{{{\SetFigFont{8}{9.6}{\rmdefault}{\mddefault}{\updefault}$=$ right good quadrilaterals}}}}}
\put(3525,1050){\makebox(0,0)[lb]{\smash{{{\SetFigFont{8}{9.6}{\rmdefault}{\mddefault}{\updefault}$=$ left good quadrilaterals}}}}}
\path(900,2400)(2700,2400)(2700,300)
	(900,300)(900,2400)
\texture{88555555 55000000 555555 55000000 555555 55000000 555555 55000000 
	555555 55000000 555555 55000000 555555 55000000 555555 55000000 
	555555 55000000 555555 55000000 555555 55000000 555555 55000000 
	555555 55000000 555555 55000000 555555 55000000 555555 55000000 }
\shade\path(900,2400)(2100,2400)(2100,1800)
	(900,1800)(900,2400)
\path(900,2400)(2100,2400)(2100,1800)
	(900,1800)(900,2400)
\texture{0 0 0 888888 88000000 0 0 80808 
	8000000 0 0 888888 88000000 0 0 80808 
	8000000 0 0 888888 88000000 0 0 80808 
	8000000 0 0 888888 88000000 0 0 80808 }
\shade\path(2100,2400)(2700,2400)(2700,300)
	(2100,300)(2100,2400)
\path(2100,2400)(2700,2400)(2700,300)
	(2100,300)(2100,2400)
\texture{55aaaaaa aa555555 55aaaaaa aa545454 54aaaaaa aa555555 55aaaaaa aa444444 
	44aaaaaa aa555555 55aaaaaa aa445444 54aaaaaa aa555555 55aaaaaa aa444444 
	44aaaaaa aa555555 55aaaaaa aa545454 54aaaaaa aa555555 55aaaaaa aa444444 
	44aaaaaa aa555555 55aaaaaa aa445444 54aaaaaa aa555555 55aaaaaa aa444444 }
\path(900,2100)(2700,2100)
\path(900,2100)(2700,2100)
\path(2100,1800)(2700,1800)
\path(2100,1800)(2700,1800)
\path(900,1800)(1800,1500)(2700,1500)
\path(900,1800)(1500,1200)(2700,1200)
\path(900,900)(2700,900)
\path(900,900)(1200,600)(2700,600)
\path(1800,1500)(1800,300)
\path(1500,1200)(1500,300)
\path(1200,600)(1200,300)
\path(2400,2400)(2400,300)
\path(2400,2400)(2400,300)
\put(2325,2475){\makebox(0,0)[lb]{\smash{{{\SetFigFont{5}{6.0}{\rmdefault}{\mddefault}{\updefault}$\cdots$}}}}}
\put(2625,2475){\makebox(0,0)[lb]{\smash{{{\SetFigFont{5}{6.0}{\rmdefault}{\mddefault}{\updefault}$(r_1,c_x)=\mathbf{a'}$}}}}}
\put(675,1350){\makebox(0,0)[lb]{\smash{{{\SetFigFont{5}{6.0}{\rmdefault}{\mddefault}{\updefault}$\vdots$}}}}}
\put(1650,2700){\makebox(0,0)[lb]{\smash{{{\SetFigFont{8}{9.6}{\rmdefault}{\mddefault}{\updefault}top}}}}}
\put(1500,0){\makebox(0,0)[lb]{\smash{{{\SetFigFont{8}{9.6}{\rmdefault}{\mddefault}{\updefault}bottom}}}}}
\put(1875,2475){\makebox(0,0)[lb]{\smash{{{\SetFigFont{5}{6.0}{\rmdefault}{\mddefault}{\updefault}$(r_1,c_2)$}}}}}
\put(450,2100){\makebox(0,0)[lb]{\smash{{{\SetFigFont{5}{6.0}{\rmdefault}{\mddefault}{\updefault}$(r_2,c_1)$}}}}}
\put(0,1500){\makebox(0,0)[lb]{\smash{{{\SetFigFont{8}{9.6}{\rmdefault}{\mddefault}{\updefault}left}}}}}
\put(1650,1575){\makebox(0,0)[lb]{\smash{{{\SetFigFont{5}{6.0}{\rmdefault}{\mddefault}{\updefault}$(R_1,C_1)$}}}}}
\put(975,675){\makebox(0,0)[lb]{\smash{{{\SetFigFont{5}{6.0}{\rmdefault}{\mddefault}{\updefault}$(R_z,C_z)$}}}}}
\put(450,2475){\makebox(0,0)[lb]{\smash{{{\SetFigFont{5}{6.0}{\rmdefault}{\mddefault}{\updefault}$(r_1,c_1) = \inf(\mathbf{a},\mathbf{a'})$}}}}}
\put(450,225){\makebox(0,0)[lb]{\smash{{{\SetFigFont{5}{6.0}{\rmdefault}{\mddefault}{\updefault}$(r_y,c_1)=\mathbf{a}$}}}}}
\put(2700,225){\makebox(0,0)[lb]{\smash{{{\SetFigFont{5}{6.0}{\rmdefault}{\mddefault}{\updefault}$(r_y,c_x)$}}}}}
\put(2850,1500){\makebox(0,0)[lb]{\smash{{{\SetFigFont{8}{9.6}{\rmdefault}{\mddefault}{\updefault}right}}}}}
\end{picture}
} \end{center}
\caption{An example of $\poset'_\circ$
\lremind{tuscaloosa}}
\label{tuscaloosa}
\end{figure}

\tpoint{Proposition} \label{oakpark2} \lremind{oakpark2} {\em The
  irreducible components of $Q_\ell$ on $\PF(\poset'_\circ)$ (described in the previous paragraph) have codimension at least $\ell$
  in $\PF(\poset'_{\circ})$, and those components of codimension exactly
  $\ell$ are strata of the form described in
  Theorem~\ref{detroit} (corresponding
to sets of good quadrilaterals, at most one per column).}

\bpf
For convenience,
we assume $\dim \inf(\mathbf{a},\mathbf{a'})=0$.  (Geometrically,
this corresponds to considering instead the
quotient of all subspaces by $V_{\inf(\mathbf{a},\mathbf{a'})}$.)

Fix an irreducible component of $Q_\ell$ of 
codimension at most $\ell$, and choose
a general point of this component; this corresponds
to some configuration of subspaces $\{ V_{\mathbf{c}} \subset K^m : \mathbf{c} \in \poset'_{\circ} \}$.
Label each element $\mathbf{c}$ of $\poset'_{\circ}$ with $\dim (V_{\mathbf{c}} \cap V_{\mathbf{a}})$.  Hence the label on any vertex dominating $\mathbf{a}$ is
$\dim \mathbf{a}$, and the label on $\mathbf{a'}$ is $\ell$.  
The label on any $\mathbf{c} \in \poset'_{\circ}$ is at most $\dim \mathbf{c}$.

This point lies in some open stratum of 
$\PF(\poset'_{\circ})$; mark the quadrilaterals
corresponding to that stratum with ``$=$''.  Then each 
labeled quadrilateral looks like one of the examples in Figure~\ref{allentown}.  

We argue by induction on the number of good quadrilaterals.  
The base case is
given in {\em (c)} and 
{\em (d)} below, and the inductive step is given in 
{\em (a)} and
{\em (b)}.

\begin{figure}[ht]
\begin{center} \setlength{\unitlength}{0.00083333in}
\begingroup\makeatletter\ifx\SetFigFont\undefined%
\gdef\SetFigFont#1#2#3#4#5{%
  \reset@font\fontsize{#1}{#2pt}%
  \fontfamily{#3}\fontseries{#4}\fontshape{#5}%
  \selectfont}%
\fi\endgroup%
{\renewcommand{\dashlinestretch}{30}
\begin{picture}(5937,2005)(0,-10)
\path(375,481)(600,556)
\path(450,406)(675,481)
\path(1575,481)(1800,556)
\path(1650,406)(1875,481)
\path(5325,481)(5550,556)
\path(5400,406)(5625,481)
\path(3375,481)(3600,556)
\path(3450,406)(3675,481)
\path(3075,781)(3225,331)(3975,181)
	(3825,631)(3075,781)
\path(75,1831)(225,1381)(975,1231)
	(825,1681)(75,1831)
\path(1275,1831)(1425,1381)(2175,1231)
	(2025,1681)(1275,1831)
\path(2475,1831)(2625,1381)(3375,1231)
	(3225,1681)(2475,1831)
\path(3675,1831)(3825,1381)(4575,1231)
	(4425,1681)(3675,1831)
\path(75,781)(225,331)(975,181)
	(825,631)(75,781)
\path(1275,781)(1425,331)(2175,181)
	(2025,631)(1275,781)
\path(4875,1831)(5025,1381)(5775,1231)
	(5625,1681)(4875,1831)
\path(5025,781)(5175,331)(5925,181)
	(5775,631)(5025,781)
\put(0,1906){\makebox(0,0)[lb]{\smash{{{\SetFigFont{8}{9.6}{\rmdefault}{\mddefault}{\updefault}$m$}}}}}
\put(0,856){\makebox(0,0)[lb]{\smash{{{\SetFigFont{8}{9.6}{\rmdefault}{\mddefault}{\updefault}$m$}}}}}
\put(75,1231){\makebox(0,0)[lb]{\smash{{{\SetFigFont{8}{9.6}{\rmdefault}{\mddefault}{\updefault}$m$}}}}}
\put(825,1756){\makebox(0,0)[lb]{\smash{{{\SetFigFont{8}{9.6}{\rmdefault}{\mddefault}{\updefault}$m$}}}}}
\put(975,1081){\makebox(0,0)[lb]{\smash{{{\SetFigFont{8}{9.6}{\rmdefault}{\mddefault}{\updefault}$m$}}}}}
\put(975,31){\makebox(0,0)[lb]{\smash{{{\SetFigFont{8}{9.6}{\rmdefault}{\mddefault}{\updefault}$m$}}}}}
\put(75,181){\makebox(0,0)[lb]{\smash{{{\SetFigFont{8}{9.6}{\rmdefault}{\mddefault}{\updefault}$m$}}}}}
\put(825,706){\makebox(0,0)[lb]{\smash{{{\SetFigFont{8}{9.6}{\rmdefault}{\mddefault}{\updefault}$m$}}}}}
\put(1200,1906){\makebox(0,0)[lb]{\smash{{{\SetFigFont{8}{9.6}{\rmdefault}{\mddefault}{\updefault}$m$}}}}}
\put(2025,1756){\makebox(0,0)[lb]{\smash{{{\SetFigFont{8}{9.6}{\rmdefault}{\mddefault}{\updefault}$m$}}}}}
\put(1275,1231){\makebox(0,0)[lb]{\smash{{{\SetFigFont{8}{9.6}{\rmdefault}{\mddefault}{\updefault}$m$}}}}}
\put(1200,856){\makebox(0,0)[lb]{\smash{{{\SetFigFont{8}{9.6}{\rmdefault}{\mddefault}{\updefault}$m$}}}}}
\put(2025,706){\makebox(0,0)[lb]{\smash{{{\SetFigFont{8}{9.6}{\rmdefault}{\mddefault}{\updefault}$m$}}}}}
\put(1275,181){\makebox(0,0)[lb]{\smash{{{\SetFigFont{8}{9.6}{\rmdefault}{\mddefault}{\updefault}$m$}}}}}
\put(2025,1081){\makebox(0,0)[lb]{\smash{{{\SetFigFont{8}{9.6}{\rmdefault}{\mddefault}{\updefault}$m+1$}}}}}
\put(2025,31){\makebox(0,0)[lb]{\smash{{{\SetFigFont{8}{9.6}{\rmdefault}{\mddefault}{\updefault}$m+1$}}}}}
\put(3825,31){\makebox(0,0)[lb]{\smash{{{\SetFigFont{8}{9.6}{\rmdefault}{\mddefault}{\updefault}$m+1$}}}}}
\put(3225,1081){\makebox(0,0)[lb]{\smash{{{\SetFigFont{8}{9.6}{\rmdefault}{\mddefault}{\updefault}$m+1$}}}}}
\put(4425,1081){\makebox(0,0)[lb]{\smash{{{\SetFigFont{8}{9.6}{\rmdefault}{\mddefault}{\updefault}$m+1$}}}}}
\put(2400,1906){\makebox(0,0)[lb]{\smash{{{\SetFigFont{8}{9.6}{\rmdefault}{\mddefault}{\updefault}$m$}}}}}
\put(3600,1906){\makebox(0,0)[lb]{\smash{{{\SetFigFont{8}{9.6}{\rmdefault}{\mddefault}{\updefault}$m$}}}}}
\put(3000,856){\makebox(0,0)[lb]{\smash{{{\SetFigFont{8}{9.6}{\rmdefault}{\mddefault}{\updefault}$m$}}}}}
\put(3150,1756){\makebox(0,0)[lb]{\smash{{{\SetFigFont{8}{9.6}{\rmdefault}{\mddefault}{\updefault}$m$}}}}}
\put(3675,1231){\makebox(0,0)[lb]{\smash{{{\SetFigFont{8}{9.6}{\rmdefault}{\mddefault}{\updefault}$m$}}}}}
\put(2400,1231){\makebox(0,0)[lb]{\smash{{{\SetFigFont{8}{9.6}{\rmdefault}{\mddefault}{\updefault}$m+1$}}}}}
\put(4275,1756){\makebox(0,0)[lb]{\smash{{{\SetFigFont{8}{9.6}{\rmdefault}{\mddefault}{\updefault}$m+1$}}}}}
\put(3675,706){\makebox(0,0)[lb]{\smash{{{\SetFigFont{8}{9.6}{\rmdefault}{\mddefault}{\updefault}$m+1$}}}}}
\put(2925,181){\makebox(0,0)[lb]{\smash{{{\SetFigFont{8}{9.6}{\rmdefault}{\mddefault}{\updefault}$m+1$}}}}}
\put(4800,1906){\makebox(0,0)[lb]{\smash{{{\SetFigFont{8}{9.6}{\rmdefault}{\mddefault}{\updefault}$m$}}}}}
\put(4950,856){\makebox(0,0)[lb]{\smash{{{\SetFigFont{8}{9.6}{\rmdefault}{\mddefault}{\updefault}$m$}}}}}
\put(5475,1756){\makebox(0,0)[lb]{\smash{{{\SetFigFont{8}{9.6}{\rmdefault}{\mddefault}{\updefault}$m+1$}}}}}
\put(5625,706){\makebox(0,0)[lb]{\smash{{{\SetFigFont{8}{9.6}{\rmdefault}{\mddefault}{\updefault}$m+1$}}}}}
\put(4950,181){\makebox(0,0)[lb]{\smash{{{\SetFigFont{8}{9.6}{\rmdefault}{\mddefault}{\updefault}$m+1$}}}}}
\put(4800,1231){\makebox(0,0)[lb]{\smash{{{\SetFigFont{8}{9.6}{\rmdefault}{\mddefault}{\updefault}$m+1$}}}}}
\put(5625,1081){\makebox(0,0)[lb]{\smash{{{\SetFigFont{8}{9.6}{\rmdefault}{\mddefault}{\updefault}$m+2$}}}}}
\put(5775,31){\makebox(0,0)[lb]{\smash{{{\SetFigFont{8}{9.6}{\rmdefault}{\mddefault}{\updefault}$m+2$}}}}}
\end{picture}
} \end{center}
\caption{Possible quadrilaterals in proof of Proposition~\ref{oakpark2}
(where $\mathbf{c}$ is labeled with $\dim(V_{\mathbf{c}} \cap V_{\mathbf{a}})$)
\lremind{allentown}}
\label{allentown}
\end{figure}
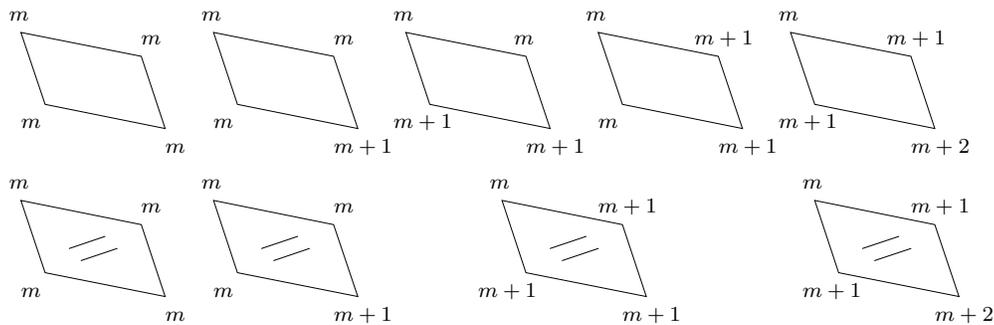

{\em (a)}
First, if $x>2$, consider the
right-most column of good quadrilaterals (see Figure~\ref{atlanticcity}).  The
two vertices in the bottom row are labeled $\dim \mathbf{a}$, and the vertex
in the upper right is labeled $\ell$.  The vertex in the upper left of the column
is labeled either $\ell-1$ or $\ell$.  If it is labeled $\ell$, then
the resulting poset with the right-most column of vertices removed
also satisfies the hypothesis of the Proposition, and hence 
by the inductive hypothesis there
must be at least $\ell$ equal signs further to the left, of
the desired form (all good quadrilaterals, at most one per column).  
Hence the Proposition holds in this
case.

\begin{figure}[ht]
\begin{center} \setlength{\unitlength}{0.00083333in}
\begingroup\makeatletter\ifx\SetFigFont\undefined%
\gdef\SetFigFont#1#2#3#4#5{%
  \reset@font\fontsize{#1}{#2pt}%
  \fontfamily{#3}\fontseries{#4}\fontshape{#5}%
  \selectfont}%
\fi\endgroup%
{\renewcommand{\dashlinestretch}{30}
\begin{picture}(1512,1705)(0,-10)
\path(300,1456)(1500,1456)(1500,256)
	(300,256)(300,1456)
\path(300,1156)(1500,1156)
\path(300,856)(1500,856)
\path(300,556)(1500,556)
\put(0,1606){\makebox(0,0)[lb]{\smash{{{\SetFigFont{8}{9.6}{\rmdefault}{\mddefault}{\updefault}$\ell-1$ or $\ell$}}}}}
\put(0,31){\makebox(0,0)[lb]{\smash{{{\SetFigFont{8}{9.6}{\rmdefault}{\mddefault}{\updefault}$\dim \mathbf{a}$}}}}}
\put(1500,1606){\makebox(0,0)[lb]{\smash{{{\SetFigFont{8}{9.6}{\rmdefault}{\mddefault}{\updefault}$\ell$}}}}}
\put(1275,31){\makebox(0,0)[lb]{\smash{{{\SetFigFont{8}{9.6}{\rmdefault}{\mddefault}{\updefault}$\dim \mathbf{a}$}}}}}
\end{picture}
} \end{center}
\caption{The right-most column of quadrilaterals (in case (a) of the proof of Proposition~\ref{oakpark2}) \lremind{atlanticcity}}
\label{atlanticcity}
\end{figure}

If otherwise the vertex in the upper left corner of the right-most
column of quadrilaterals  is labeled $\ell-1$, then by the inductive
hypothesis there must be at least $\ell-1$ equal signs further to the
left, of the desired form.  Furthermore,
by inspection of  Figure~\ref{allentown}, there must be another equal sign in 
the right-most column as well, and again the Proposition holds.  
See Figure~\ref{atlanticcity2} for an example.

\begin{figure}[ht]
\begin{center} \setlength{\unitlength}{0.00083333in}
\begingroup\makeatletter\ifx\SetFigFont\undefined%
\gdef\SetFigFont#1#2#3#4#5{%
  \reset@font\fontsize{#1}{#2pt}%
  \fontfamily{#3}\fontseries{#4}\fontshape{#5}%
  \selectfont}%
\fi\endgroup%
{\renewcommand{\dashlinestretch}{30}
\begin{picture}(1425,1933)(0,-10)
\path(150,1606)(1350,1606)
\path(150,1306)(1350,1306)
\path(150,1006)(1350,1006)
\path(150,1906)(1350,1906)(1350,106)
	(150,106)(150,1906)
\path(150,706)(1350,706)
\path(150,406)(1350,406)
\path(600,1081)(900,1156)
\path(600,1156)(900,1231)
\put(0,1831){\makebox(0,0)[lb]{\smash{{{\SetFigFont{8}{9.6}{\rmdefault}{\mddefault}{\updefault}$1$}}}}}
\put(0,1531){\makebox(0,0)[lb]{\smash{{{\SetFigFont{8}{9.6}{\rmdefault}{\mddefault}{\updefault}$2$}}}}}
\put(0,1231){\makebox(0,0)[lb]{\smash{{{\SetFigFont{8}{9.6}{\rmdefault}{\mddefault}{\updefault}$2$}}}}}
\put(0,931){\makebox(0,0)[lb]{\smash{{{\SetFigFont{8}{9.6}{\rmdefault}{\mddefault}{\updefault}$3$}}}}}
\put(0,631){\makebox(0,0)[lb]{\smash{{{\SetFigFont{8}{9.6}{\rmdefault}{\mddefault}{\updefault}$4$}}}}}
\put(0,331){\makebox(0,0)[lb]{\smash{{{\SetFigFont{8}{9.6}{\rmdefault}{\mddefault}{\updefault}$5$}}}}}
\put(0,31){\makebox(0,0)[lb]{\smash{{{\SetFigFont{8}{9.6}{\rmdefault}{\mddefault}{\updefault}$5$}}}}}
\put(1425,1831){\makebox(0,0)[lb]{\smash{{{\SetFigFont{8}{9.6}{\rmdefault}{\mddefault}{\updefault}$2$}}}}}
\put(1425,1531){\makebox(0,0)[lb]{\smash{{{\SetFigFont{8}{9.6}{\rmdefault}{\mddefault}{\updefault}$3$}}}}}
\put(1425,1231){\makebox(0,0)[lb]{\smash{{{\SetFigFont{8}{9.6}{\rmdefault}{\mddefault}{\updefault}$3$}}}}}
\put(1425,931){\makebox(0,0)[lb]{\smash{{{\SetFigFont{8}{9.6}{\rmdefault}{\mddefault}{\updefault}$3$}}}}}
\put(1425,631){\makebox(0,0)[lb]{\smash{{{\SetFigFont{8}{9.6}{\rmdefault}{\mddefault}{\updefault}$4$}}}}}
\put(1425,331){\makebox(0,0)[lb]{\smash{{{\SetFigFont{8}{9.6}{\rmdefault}{\mddefault}{\updefault}$5$}}}}}
\put(1425,31){\makebox(0,0)[lb]{\smash{{{\SetFigFont{8}{9.6}{\rmdefault}{\mddefault}{\updefault}$5$}}}}}
\end{picture}
} \end{center}
\caption{Sample labeling of the right-most column of quadrilaterals (in case (a) of the proof of Proposition~\ref{oakpark2})
\lremind{atlanticcity2}}
\label{atlanticcity2}
\end{figure}

{\em (b)} Next suppose $x=2$ (so there are no right good quadrilaterals),
and there is at least one left good
quadrilateral (hence there is no blocker).  
Let $\mathbf{w}=(r_1,c_2)$ (the rightmost element of the
top row; this corresponds to the white
checker in the critical row of $\circ$).  The label on $\mathbf{w}$ must be 0 or
1,  as $\dim \mathbf{w} =1$.  If $\mathbf{w}$ has label 0, then the Proposition 
follows 
immediately.  If $\mathbf{w}$ has label 1, then the top left quadrilateral
 must be of one of the  forms shown in
Figure~\ref{jerseycity}.  In each case, we can remove
the top two vertices of the good quadrilateral, and in all but
the first case subtract 1 from each
of the other vertices (as shown in Figure~\ref{jerseycity}),
and reduce to the case with one fewer good quadrilaterals.    
Hence by the inductive hypothesis, the last case cannot happen
(there must be a second ``$=$'' in a good quadrilateral further below,
yet we are allowed only $\ell=1$).  The other cases
may occur, and by our inductive hypothesis the Proposition holds.

\begin{figure}[ht]
\begin{center} \setlength{\unitlength}{0.00083333in}
\begingroup\makeatletter\ifx\SetFigFont\undefined%
\gdef\SetFigFont#1#2#3#4#5{%
  \reset@font\fontsize{#1}{#2pt}%
  \fontfamily{#3}\fontseries{#4}\fontshape{#5}%
  \selectfont}%
\fi\endgroup%
{\renewcommand{\dashlinestretch}{30}
\begin{picture}(4737,2155)(0,-10)
\path(75,1981)(675,1981)(675,1381)
	(75,1381)(75,1981)
\path(1275,1381)(1875,1381)
\path(75,781)(675,781)(675,181)
	(75,181)(75,781)
\path(1275,181)(1875,181)
\path(2925,1981)(3525,1981)(3525,1381)
	(2925,1381)(2925,1981)
\path(4125,1381)(4725,1381)
\path(4125,181)(4725,181)
\path(150,256)(600,631)
\path(150,331)(600,706)
\path(3000,256)(3450,631)
\path(3000,331)(3450,706)
\path(2925,781)(3525,781)(3525,181)
	(2925,181)(2925,781)
\put(1200,1231){\makebox(0,0)[lb]{\smash{{{\SetFigFont{8}{9.6}{\rmdefault}{\mddefault}{\updefault}$0$}}}}}
\put(675,1231){\makebox(0,0)[lb]{\smash{{{\SetFigFont{8}{9.6}{\rmdefault}{\mddefault}{\updefault}$1$}}}}}
\put(1875,1231){\makebox(0,0)[lb]{\smash{{{\SetFigFont{8}{9.6}{\rmdefault}{\mddefault}{\updefault}$1$}}}}}
\put(0,1231){\makebox(0,0)[lb]{\smash{{{\SetFigFont{8}{9.6}{\rmdefault}{\mddefault}{\updefault}$0$}}}}}
\put(675,2056){\makebox(0,0)[lb]{\smash{{{\SetFigFont{8}{9.6}{\rmdefault}{\mddefault}{\updefault}$1$}}}}}
\put(0,2056){\makebox(0,0)[lb]{\smash{{{\SetFigFont{8}{9.6}{\rmdefault}{\mddefault}{\updefault}$0$}}}}}
\put(900,1606){\makebox(0,0)[lb]{\smash{{{\SetFigFont{8}{9.6}{\rmdefault}{\mddefault}{\updefault}$\Longrightarrow$}}}}}
\put(1200,31){\makebox(0,0)[lb]{\smash{{{\SetFigFont{8}{9.6}{\rmdefault}{\mddefault}{\updefault}$0$}}}}}
\put(675,31){\makebox(0,0)[lb]{\smash{{{\SetFigFont{8}{9.6}{\rmdefault}{\mddefault}{\updefault}$1$}}}}}
\put(1875,31){\makebox(0,0)[lb]{\smash{{{\SetFigFont{8}{9.6}{\rmdefault}{\mddefault}{\updefault}$0$}}}}}
\put(0,31){\makebox(0,0)[lb]{\smash{{{\SetFigFont{8}{9.6}{\rmdefault}{\mddefault}{\updefault}$1$}}}}}
\put(675,856){\makebox(0,0)[lb]{\smash{{{\SetFigFont{8}{9.6}{\rmdefault}{\mddefault}{\updefault}$1$}}}}}
\put(0,856){\makebox(0,0)[lb]{\smash{{{\SetFigFont{8}{9.6}{\rmdefault}{\mddefault}{\updefault}$0$}}}}}
\put(900,406){\makebox(0,0)[lb]{\smash{{{\SetFigFont{8}{9.6}{\rmdefault}{\mddefault}{\updefault}$\Longrightarrow$}}}}}
\put(4050,1231){\makebox(0,0)[lb]{\smash{{{\SetFigFont{8}{9.6}{\rmdefault}{\mddefault}{\updefault}$0$}}}}}
\put(3525,1231){\makebox(0,0)[lb]{\smash{{{\SetFigFont{8}{9.6}{\rmdefault}{\mddefault}{\updefault}$2$}}}}}
\put(4725,1231){\makebox(0,0)[lb]{\smash{{{\SetFigFont{8}{9.6}{\rmdefault}{\mddefault}{\updefault}$1$}}}}}
\put(2850,1231){\makebox(0,0)[lb]{\smash{{{\SetFigFont{8}{9.6}{\rmdefault}{\mddefault}{\updefault}$1$}}}}}
\put(3525,2056){\makebox(0,0)[lb]{\smash{{{\SetFigFont{8}{9.6}{\rmdefault}{\mddefault}{\updefault}$1$}}}}}
\put(2850,2056){\makebox(0,0)[lb]{\smash{{{\SetFigFont{8}{9.6}{\rmdefault}{\mddefault}{\updefault}$0$}}}}}
\put(3750,1606){\makebox(0,0)[lb]{\smash{{{\SetFigFont{8}{9.6}{\rmdefault}{\mddefault}{\updefault}$\Longrightarrow$}}}}}
\put(4050,31){\makebox(0,0)[lb]{\smash{{{\SetFigFont{8}{9.6}{\rmdefault}{\mddefault}{\updefault}$0$}}}}}
\put(3525,31){\makebox(0,0)[lb]{\smash{{{\SetFigFont{8}{9.6}{\rmdefault}{\mddefault}{\updefault}$2$}}}}}
\put(4725,31){\makebox(0,0)[lb]{\smash{{{\SetFigFont{8}{9.6}{\rmdefault}{\mddefault}{\updefault}$1$}}}}}
\put(2850,31){\makebox(0,0)[lb]{\smash{{{\SetFigFont{8}{9.6}{\rmdefault}{\mddefault}{\updefault}$1$}}}}}
\put(3525,856){\makebox(0,0)[lb]{\smash{{{\SetFigFont{8}{9.6}{\rmdefault}{\mddefault}{\updefault}$1$}}}}}
\put(2850,856){\makebox(0,0)[lb]{\smash{{{\SetFigFont{8}{9.6}{\rmdefault}{\mddefault}{\updefault}$0$}}}}}
\put(3750,406){\makebox(0,0)[lb]{\smash{{{\SetFigFont{8}{9.6}{\rmdefault}{\mddefault}{\updefault}$\Longrightarrow$}}}}}
\put(750,1906){\makebox(0,0)[lb]{\smash{{{\SetFigFont{8}{9.6}{\rmdefault}{\mddefault}{\updefault}$\mathbf{w}$}}}}}
\put(750,706){\makebox(0,0)[lb]{\smash{{{\SetFigFont{8}{9.6}{\rmdefault}{\mddefault}{\updefault}$\mathbf{w}$}}}}}
\put(3600,706){\makebox(0,0)[lb]{\smash{{{\SetFigFont{8}{9.6}{\rmdefault}{\mddefault}{\updefault}$\mathbf{w}$}}}}}
\put(3600,1906){\makebox(0,0)[lb]{\smash{{{\SetFigFont{8}{9.6}{\rmdefault}{\mddefault}{\updefault}$\mathbf{w}$}}}}}
\end{picture}
} \end{center}
\caption{The top (left) good quadrilateral  (in case (b) of the proof of Proposition~\ref{oakpark2}) \lremind{jerseycity}
}
\label{jerseycity}
\end{figure}
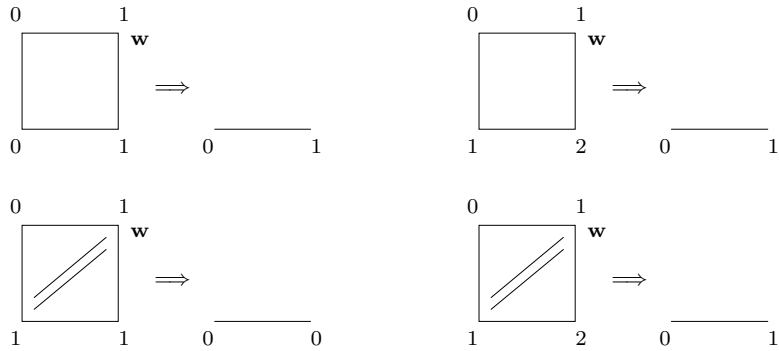

{\em (c)} Next, if there are no good quadrilaterals, and no blockers,
then the labeled poset is the one shown in Figure~\ref{lawrence}, and the
result is immediate.

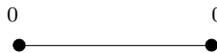
\begin{figure}[ht]
\begin{center} \setlength{\unitlength}{0.00083333in}
\begingroup\makeatletter\ifx\SetFigFont\undefined%
\gdef\SetFigFont#1#2#3#4#5{%
  \reset@font\fontsize{#1}{#2pt}%
  \fontfamily{#3}\fontseries{#4}\fontshape{#5}%
  \selectfont}%
\fi\endgroup%
{\renewcommand{\dashlinestretch}{30}
\begin{picture}(1320,294)(0,-10)
\put(75,45){\blacken\ellipse{74}{74}}
\put(75,45){\ellipse{74}{74}}
\put(1275,45){\blacken\ellipse{74}{74}}
\put(1275,45){\ellipse{74}{74}}
\path(75,45)(1275,45)
\put(0,195){\makebox(0,0)[lb]{\smash{{{\SetFigFont{8}{9.6}{\rmdefault}{\mddefault}{\updefault}$0$}}}}}
\put(1275,195){\makebox(0,0)[lb]{\smash{{{\SetFigFont{8}{9.6}{\rmdefault}{\mddefault}{\updefault}$0$}}}}}
\end{picture}
} \end{center}
\caption{The trivial case (i)(c) of the proof of Proposition~\ref{oakpark2} \lremind{lawrence}}
\label{lawrence}
\end{figure}

{\em (d)} Finally, suppose there are no good quadrilaterals (implying
$x=2$), and there is a blocker.  As in {\em (b)},
$\ell = 0$ or $1$.  If $\ell=
0$,  the Proposition is immediate.

Assume now that $\ell =1$; we will obtain a contradiction. Re-label the
poset by labeling element $\mathbf{c}$ with the dimension of the intersection
of its corresponding vector space with $V_{\mathbf{a'}}$ (rather than $V_{\mathbf{a}}$).
Hence all labels are $0$ or $1$, and the elements on the right side
and the bottom of the graph are labeled $1$.  Any quadrilateral with
minimal element labeled $0$ and other three elements labeled 1 must be
marked with an ``$=$'' (see Figure~\ref{pittsburgh3}).  By hypothesis,
as the component of $Q_1$ has codimension at most 1 in
$\PF(\poset'_{\circ})$, there can be at most one such quadrilateral.
As the subset of elements labeled 0 has a maximum, which doesn't lie
on the right side or the bottom of the graph, there is at least one
such quadrilateral, so there is exactly one.  The component of $Q_1$
lies in the (irreducible) divisorial stratum corresponding to this
quadrilateral, so these two loci are equal.  Thus a general point of
this divisorial stratum lies in $Q_1$.

\begin{figure}[ht]
\begin{center} \setlength{\unitlength}{0.00083333in}
\begingroup\makeatletter\ifx\SetFigFont\undefined%
\gdef\SetFigFont#1#2#3#4#5{%
  \reset@font\fontsize{#1}{#2pt}%
  \fontfamily{#3}\fontseries{#4}\fontshape{#5}%
  \selectfont}%
\fi\endgroup%
{\renewcommand{\dashlinestretch}{30}
\begin{picture}(1200,805)(0,-10)
\path(225,706)(1125,706)(1125,106)
	(225,106)(225,706)
\path(300,256)(1050,631)
\path(1050,556)(300,181)
\put(1200,706){\makebox(0,0)[lb]{\smash{{{\SetFigFont{8}{9.6}{\rmdefault}{\mddefault}{\updefault}$1$}}}}}
\put(0,706){\makebox(0,0)[lb]{\smash{{{\SetFigFont{8}{9.6}{\rmdefault}{\mddefault}{\updefault}$0$}}}}}
\put(0,31){\makebox(0,0)[lb]{\smash{{{\SetFigFont{8}{9.6}{\rmdefault}{\mddefault}{\updefault}$1$}}}}}
\put(1200,31){\makebox(0,0)[lb]{\smash{{{\SetFigFont{8}{9.6}{\rmdefault}{\mddefault}{\updefault}$1$}}}}}
\end{picture}
} \end{center}
\caption{ \lremind{pittsburgh3}}
\label{pittsburgh3}
\end{figure}

We construct the divisorial stratum as a tower of $\proj^1$-bundles
over the partial flag variety parametrizing the subspaces
corresponding to the vertices on the northeast border (as in the proof
of Lemma~\ref{bssmooth}, where the southwest border was used).
There can be no quadrilateral labeled as in
Figure~\ref{pittsburgh} (with no ``$=$''): $V_{\mathbf{e}}$ is a general
subspace (of $\dim \mathbf{e}$) containing $V_{\mathbf{c}}$ and contained in
$V_{\mathbf{f}}$; as $\dim V_{\mathbf{d}} \cap V_{\mathbf{a'}} = 0$, we have $\dim
V_{\mathbf{e}} \cap V_{\mathbf{a'}}=0$ as well, so the label on the lower left
vertex must be $0$, not $1$.

\begin{figure}[ht]
\begin{center} \setlength{\unitlength}{0.00083333in}
\begingroup\makeatletter\ifx\SetFigFont\undefined%
\gdef\SetFigFont#1#2#3#4#5{%
  \reset@font\fontsize{#1}{#2pt}%
  \fontfamily{#3}\fontseries{#4}\fontshape{#5}%
  \selectfont}%
\fi\endgroup%
{\renewcommand{\dashlinestretch}{30}
\begin{picture}(1200,1030)(0,-10)
\path(225,856)(1125,856)(1125,256)
	(225,256)(225,856)
\put(1200,781){\makebox(0,0)[lb]{\smash{{{\SetFigFont{8}{9.6}{\rmdefault}{\mddefault}{\updefault}$\mathbf{d}$}}}}}
\put(225,31){\makebox(0,0)[lb]{\smash{{{\SetFigFont{8}{9.6}{\rmdefault}{\mddefault}{\updefault}$1$}}}}}
\put(1050,31){\makebox(0,0)[lb]{\smash{{{\SetFigFont{8}{9.6}{\rmdefault}{\mddefault}{\updefault}$1$}}}}}
\put(225,931){\makebox(0,0)[lb]{\smash{{{\SetFigFont{8}{9.6}{\rmdefault}{\mddefault}{\updefault}$0$}}}}}
\put(1050,931){\makebox(0,0)[lb]{\smash{{{\SetFigFont{8}{9.6}{\rmdefault}{\mddefault}{\updefault}$0$}}}}}
\put(0,781){\makebox(0,0)[lb]{\smash{{{\SetFigFont{8}{9.6}{\rmdefault}{\mddefault}{\updefault}$\mathbf{c}$}}}}}
\put(1200,181){\makebox(0,0)[lb]{\smash{{{\SetFigFont{8}{9.6}{\rmdefault}{\mddefault}{\updefault}$\mathbf{f}$}}}}}
\put(0,181){\makebox(0,0)[lb]{\smash{{{\SetFigFont{8}{9.6}{\rmdefault}{\mddefault}{\updefault}$\mathbf{e}$}}}}}
\end{picture}
} \end{center}
\caption{ \lremind{pittsburgh}}
\label{pittsburgh}
\end{figure}

Now consider the part of the graph above row $r_2$, as in
Figure~\ref{pittsburgh2}.  As there are 
no quadrilaterals marked as in
Figure~\ref{pittsburgh}, all unlabeled vertices in Figure~\ref{pittsburgh2} must have label 1 (work inductively bottom to top, and left to right in each row), contradicting the fact that there is 
precisely one quadrilateral marked as in
Figure~\ref{pittsburgh3}.\epf

\begin{figure}[ht]
\begin{center} \setlength{\unitlength}{0.00083333in}
\begingroup\makeatletter\ifx\SetFigFont\undefined%
\gdef\SetFigFont#1#2#3#4#5{%
  \reset@font\fontsize{#1}{#2pt}%
  \fontfamily{#3}\fontseries{#4}\fontshape{#5}%
  \selectfont}%
\fi\endgroup%
{\renewcommand{\dashlinestretch}{30}
\begin{picture}(1500,1630)(0,-10)
\path(1425,1156)(1425,256)(525,256)
	(225,1456)(825,556)(1425,556)
\path(825,256)(825,556)
\path(1125,256)(1125,856)(1425,856)
\path(1425,1156)(225,1456)(1125,856)
\put(0,1531){\makebox(0,0)[lb]{\smash{{{\SetFigFont{8}{9.6}{\rmdefault}{\mddefault}{\updefault}$0$}}}}}
\put(1500,1081){\makebox(0,0)[lb]{\smash{{{\SetFigFont{8}{9.6}{\rmdefault}{\mddefault}{\updefault}$1$}}}}}
\put(1500,781){\makebox(0,0)[lb]{\smash{{{\SetFigFont{8}{9.6}{\rmdefault}{\mddefault}{\updefault}$1$}}}}}
\put(1500,481){\makebox(0,0)[lb]{\smash{{{\SetFigFont{8}{9.6}{\rmdefault}{\mddefault}{\updefault}$1$}}}}}
\put(1500,31){\makebox(0,0)[lb]{\smash{{{\SetFigFont{8}{9.6}{\rmdefault}{\mddefault}{\updefault}$1$}}}}}
\put(750,31){\makebox(0,0)[lb]{\smash{{{\SetFigFont{8}{9.6}{\rmdefault}{\mddefault}{\updefault}$1$}}}}}
\put(450,31){\makebox(0,0)[lb]{\smash{{{\SetFigFont{8}{9.6}{\rmdefault}{\mddefault}{\updefault}$1$}}}}}
\put(1050,31){\makebox(0,0)[lb]{\smash{{{\SetFigFont{8}{9.6}{\rmdefault}{\mddefault}{\updefault}$1$}}}}}
\end{picture}
} \end{center}
\caption{\lremind{pittsburgh2}}
\label{pittsburgh2}
\end{figure}

\epoint{Remark}  The irreducible components of $Q_\ell$ correspond
to sets of good quadrilaterals, at most one per column,
such that {\em each is in a higher row than any to its left}.  
This describes the components of the pullback of $D_{\Box}$ in
Theorem~\ref{detroit} completely.  To prove this, one must slightly
refine case (a) of the proof of Proposition~\ref{oakpark2}.  However,
this fact will not be needed.

\bpoint{Geometric irrelevance of all but one or two divisors} 
\label{capecod} \lremind{capecod}
We
next show that all but one or two of the $D_S \subset \PF(\poset_{\circ})
\times X_{\bn}$  (defined immediately after the statement of
Theorem~\ref{detroit})
are geometrically irrelevant.  
Note that 
$D_{\{ \}}$ corresponds to the entry
``stay'' in Table~\ref{keywest}, and
$D_{ \{ \text{top left good quad.}\}}$ corresponds to ``swap''.

Proposition~\ref{tacoma} shows that all other
components are geometrically irrelevant, and Proposition~\ref{redmond}
shows
that $D_{\{ \}}$ is geometrically irrelevant as well when required to
be by Table~\ref{keywest}.

The following example may motivate the argument.  Suppose
$\cb$ is the configuration of black and white checkers shown in
Figure~\ref{lubbock1}, corresponding to planes in $\proj^3$.  Then a
general point of $X_{\cb}$ is depicted in Figure~\ref{lubbock2}.  The
four points parametrized by black checkers are denoted $b_1$, $b_2$,
$b_3$, $b_4$, and the three points parametrized by the white checkers
are denoted $p_1$, $p_2$, $p_3$.  The other subspaces parametrized by
$\PF(\poset_{\circ})$ are two lines $l_{12}$ and $l_{23}$ (containing
$\{p_1, p_2 \}$ and $\{ p_2, p_3 \}$ respectively), and a plane
$f_{123}$.  For
a fixed point of $X_{\bullet}$, this locus clearly has dimension 2,
(generically) one for the choice of $p_2$ on the line $b_2b_3$ and one for the choice of
$p_3$ on the line $b_3 b_4$.

In the degeneration corresponding to $X_{\bn}$, $b_2$ tends to $b_1$
(and the line of approach $\ell$ is remembered).
Note that there are two good quadrilaterals in $\poset'_{\circ}$, one
``left'' and one ``right''; denote them $L$ and $R$.
There are then three possible degenerations of 
Figure~\ref{lubbock2} (i.e. three components of $D'_X$),
depicted in Figure~\ref{lubbock3}.

\begin{figure}[ht]
\begin{center}\setlength{\unitlength}{0.00083333in}
\begingroup\makeatletter\ifx\SetFigFont\undefined%
\gdef\SetFigFont#1#2#3#4#5{%
  \reset@font\fontsize{#1}{#2pt}%
  \fontfamily{#3}\fontseries{#4}\fontshape{#5}%
  \selectfont}%
\fi\endgroup%
{\renewcommand{\dashlinestretch}{30}
\begin{picture}(1224,1239)(0,-10)
\put(1062,837){\ellipse{74}{74}}
\put(1062,1137){\blacken\ellipse{74}{74}}
\put(1062,1137){\ellipse{74}{74}}
\put(762,837){\blacken\ellipse{74}{74}}
\put(762,837){\ellipse{74}{74}}
\put(462,537){\blacken\ellipse{74}{74}}
\put(462,537){\ellipse{74}{74}}
\put(87,237){\blacken\ellipse{74}{74}}
\put(87,237){\ellipse{74}{74}}
\put(237,237){\ellipse{74}{74}}
\put(762,537){\ellipse{74}{74}}
\path(312,1212)(312,12)
\path(612,1212)(612,12)
\path(912,1212)(912,12)
\path(12,912)(1212,912)
\path(12,312)(1212,312)
\path(12,612)(1212,612)
\path(12,1212)(1212,1212)(1212,12)
	(12,12)(12,1212)
\put(987,387){\makebox(0,0)[lb]{\smash{{{\SetFigFont{5}{6.0}{\rmdefault}{\mddefault}{\updefault}$\ell_{23}$}}}}}
\put(162,87){\makebox(0,0)[lb]{\smash{{{\SetFigFont{5}{6.0}{\rmdefault}{\mddefault}{\updefault}$p_1$}}}}}
\put(987,87){\makebox(0,0)[lb]{\smash{{{\SetFigFont{5}{6.0}{\rmdefault}{\mddefault}{\updefault}$f_{123}$}}}}}
\put(687,387){\makebox(0,0)[lb]{\smash{{{\SetFigFont{5}{6.0}{\rmdefault}{\mddefault}{\updefault}$p_2$}}}}}
\put(387,387){\makebox(0,0)[lb]{\smash{{{\SetFigFont{5}{6.0}{\rmdefault}{\mddefault}{\updefault}$b_2$}}}}}
\put(687,687){\makebox(0,0)[lb]{\smash{{{\SetFigFont{5}{6.0}{\rmdefault}{\mddefault}{\updefault}$b_3$}}}}}
\put(987,687){\makebox(0,0)[lb]{\smash{{{\SetFigFont{5}{6.0}{\rmdefault}{\mddefault}{\updefault}$p_3$}}}}}
\put(987,987){\makebox(0,0)[lb]{\smash{{{\SetFigFont{5}{6.0}{\rmdefault}{\mddefault}{\updefault}$b_4$}}}}}
\put(687,87){\makebox(0,0)[lb]{\smash{{{\SetFigFont{5}{6.0}{\rmdefault}{\mddefault}{\updefault}$\ell_{23}$}}}}}
\put(27,87){\makebox(0,0)[lb]{\smash{{{\SetFigFont{5}{6.0}{\rmdefault}{\mddefault}{\updefault}$b_1$}}}}}
\end{picture}
}\end{center}
\caption{\lremind{lubbock1}
}
\label{lubbock1}
\end{figure}

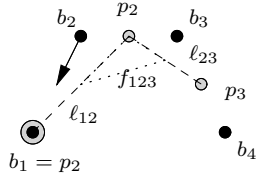
\begin{figure}[ht]
\begin{center}\setlength{\unitlength}{0.00083333in}
\begingroup\makeatletter\ifx\SetFigFont\undefined%
\gdef\SetFigFont#1#2#3#4#5{%
  \reset@font\fontsize{#1}{#2pt}%
  \fontfamily{#3}\fontseries{#4}\fontshape{#5}%
  \selectfont}%
\fi\endgroup%
{\renewcommand{\dashlinestretch}{30}
\begin{picture}(1425,1122)(0,-10)
\texture{0 115111 51000000 444444 44000000 151515 15000000 444444 
	44000000 511151 11000000 444444 44000000 151515 15000000 444444 
	44000000 115111 51000000 444444 44000000 151515 15000000 444444 
	44000000 511151 11000000 444444 44000000 151515 15000000 444444 }
\put(150,261){\shade\ellipse{150}{150}}
\put(150,261){\ellipse{150}{150}}
\put(750,861){\shade\ellipse{74}{74}}
\put(750,861){\ellipse{74}{74}}
\put(150,261){\blacken\ellipse{74}{74}}
\put(150,261){\ellipse{74}{74}}
\put(450,861){\blacken\ellipse{74}{74}}
\put(450,861){\ellipse{74}{74}}
\put(1050,861){\blacken\ellipse{74}{74}}
\put(1050,861){\ellipse{74}{74}}
\put(1350,261){\blacken\ellipse{74}{74}}
\put(1350,261){\ellipse{74}{74}}
\put(1200,561){\shade\ellipse{74}{74}}
\put(1200,561){\ellipse{74}{74}}
\path(450,861)(300,561)
\path(450,861)(300,561)
\blacken\path(326.833,681.748)(300.000,561.000)(380.498,654.915)(326.833,681.748)
\dashline{60.000}(150,261)(750,861)(1200,561)
\dottedline{60}(450,561)(750,861)(975,711)(450,561)
\put(300,936){\makebox(0,0)[lb]{\smash{{{\SetFigFont{8}{9.6}{\rmdefault}{\mddefault}{\updefault}$b_2$}}}}}
\put(1125,936){\makebox(0,0)[lb]{\smash{{{\SetFigFont{8}{9.6}{\rmdefault}{\mddefault}{\updefault}$b_3$}}}}}
\put(675,1011){\makebox(0,0)[lb]{\smash{{{\SetFigFont{8}{9.6}{\rmdefault}{\mddefault}{\updefault}$p_2$}}}}}
\put(1125,711){\makebox(0,0)[lb]{\smash{{{\SetFigFont{8}{9.6}{\rmdefault}{\mddefault}{\updefault}$\ell_{23}$}}}}}
\put(1350,486){\makebox(0,0)[lb]{\smash{{{\SetFigFont{8}{9.6}{\rmdefault}{\mddefault}{\updefault}$p_3$}}}}}
\put(1425,111){\makebox(0,0)[lb]{\smash{{{\SetFigFont{8}{9.6}{\rmdefault}{\mddefault}{\updefault}$b_4$}}}}}
\put(0,36){\makebox(0,0)[lb]{\smash{{{\SetFigFont{8}{9.6}{\rmdefault}{\mddefault}{\updefault}$b_1=p_2$}}}}}
\put(375,336){\makebox(0,0)[lb]{\smash{{{\SetFigFont{8}{9.6}{\rmdefault}{\mddefault}{\updefault}$\ell_{12}$}}}}}
\put(675,561){\makebox(0,0)[lb]{\smash{{{\SetFigFont{8}{9.6}{\rmdefault}{\mddefault}{\updefault}$f_{123}$}}}}}
\end{picture}
}\end{center}
\caption{A general point of $X_{\cb}$ (where $\cb$ is given in Figure~\ref{lubbock1}) \lremind{lubbock2}}
\label{lubbock2}
\end{figure}

\begin{figure}[ht]
\begin{center}\setlength{\unitlength}{0.00083333in}
\begingroup\makeatletter\ifx\SetFigFont\undefined%
\gdef\SetFigFont#1#2#3#4#5{%
  \reset@font\fontsize{#1}{#2pt}%
  \fontfamily{#3}\fontseries{#4}\fontshape{#5}%
  \selectfont}%
\fi\endgroup%
{\renewcommand{\dashlinestretch}{30}
\begin{picture}(4995,2085)(0,-10)
\put(262.500,1198.500){\arc{855.132}{4.9786}{5.6221}}
\put(2250.000,2098.500){\arc{1968.661}{0.7045}{1.8805}}
\put(4050.000,1686.000){\arc{1209.339}{6.1588}{8.3731}}
\put(487.500,1348.500){\arc{977.880}{5.2791}{6.2064}}
\texture{0 115111 51000000 444444 44000000 151515 15000000 444444 
	44000000 511151 11000000 444444 44000000 151515 15000000 444444 
	44000000 115111 51000000 444444 44000000 151515 15000000 444444 
	44000000 511151 11000000 444444 44000000 151515 15000000 444444 }
\put(150,1161){\shade\ellipse{150}{150}}
\put(150,1161){\ellipse{150}{150}}
\put(150,1161){\blacken\ellipse{74}{74}}
\put(150,1161){\ellipse{74}{74}}
\put(1050,1761){\blacken\ellipse{74}{74}}
\put(1050,1761){\ellipse{74}{74}}
\put(1350,1161){\blacken\ellipse{74}{74}}
\put(1350,1161){\ellipse{74}{74}}
\put(1200,1461){\shade\ellipse{74}{74}}
\put(1200,1461){\ellipse{74}{74}}
\put(1950,1161){\shade\ellipse{150}{150}}
\put(1950,1161){\ellipse{150}{150}}
\put(1950,1161){\blacken\ellipse{74}{74}}
\put(1950,1161){\ellipse{74}{74}}
\put(2850,1761){\blacken\ellipse{74}{74}}
\put(2850,1761){\ellipse{74}{74}}
\put(3150,1161){\blacken\ellipse{74}{74}}
\put(3150,1161){\ellipse{74}{74}}
\put(4950,1161){\blacken\ellipse{74}{74}}
\put(4950,1161){\ellipse{74}{74}}
\put(3750,1161){\shade\ellipse{150}{150}}
\put(3750,1161){\ellipse{150}{150}}
\put(3750,1161){\blacken\ellipse{74}{74}}
\put(3750,1161){\ellipse{74}{74}}
\put(4650,1761){\shade\ellipse{150}{150}}
\put(4650,1761){\ellipse{150}{150}}
\put(4650,1761){\blacken\ellipse{74}{74}}
\put(4650,1761){\ellipse{74}{74}}
\put(4200,1461){\shade\ellipse{74}{74}}
\put(4200,1461){\ellipse{74}{74}}
\put(2400,1461){\shade\ellipse{74}{74}}
\put(2400,1461){\ellipse{74}{74}}
\put(3000,1461){\shade\ellipse{74}{74}}
\put(3000,1461){\ellipse{74}{74}}
\dashline{60.000}(150,1161)(1200,1461)
\path(1950,1161)(2250,1761)
\path(3750,1161)(4050,1761)
\texture{88555555 55000000 555555 55000000 555555 55000000 555555 55000000 
	555555 55000000 555555 55000000 555555 55000000 555555 55000000 
	555555 55000000 555555 55000000 555555 55000000 555555 55000000 
	555555 55000000 555555 55000000 555555 55000000 555555 55000000 }
\dashline{60.000}(4650,1761)(3750,1161)
\path(3300,186)(2700,486)
\path(3300,186)(2700,486)
\blacken\path(2820.748,459.167)(2700.000,486.000)(2793.915,405.502)(2820.748,459.167)
\path(3750,186)(4200,486)
\path(3750,186)(4200,486)
\blacken\path(4116.795,394.474)(4200.000,486.000)(4083.513,444.397)(4116.795,394.474)
\dashline{60.000}(1950,1161)(2850,1761)
\dashline{60.000}(3000,1461)(2400,1461)
\dashline{60.000}(900,1911)(150,1161)
\dottedline{45}(150,1161)(600,1461)
\path(150,1161)(450,1761)
\put(0,936){\makebox(0,0)[lb]{\smash{{{\SetFigFont{8}{9.6}{\rmdefault}{\mddefault}{\updefault}$p_1=p_2$}}}}}
\put(825,1536){\makebox(0,0)[lb]{\smash{{{\SetFigFont{8}{9.6}{\rmdefault}{\mddefault}{\updefault}$f_{123}$}}}}}
\put(1350,1461){\makebox(0,0)[lb]{\smash{{{\SetFigFont{8}{9.6}{\rmdefault}{\mddefault}{\updefault}$p_3$}}}}}
\put(675,1161){\makebox(0,0)[lb]{\smash{{{\SetFigFont{8}{9.6}{\rmdefault}{\mddefault}{\updefault}$\ell_{23}$}}}}}
\put(1800,936){\makebox(0,0)[lb]{\smash{{{\SetFigFont{8}{9.6}{\rmdefault}{\mddefault}{\updefault}$p_1$}}}}}
\put(3675,936){\makebox(0,0)[lb]{\smash{{{\SetFigFont{8}{9.6}{\rmdefault}{\mddefault}{\updefault}$p_1$}}}}}
\put(4800,1836){\makebox(0,0)[lb]{\smash{{{\SetFigFont{8}{9.6}{\rmdefault}{\mddefault}{\updefault}$p_3$}}}}}
\put(3150,1386){\makebox(0,0)[lb]{\smash{{{\SetFigFont{8}{9.6}{\rmdefault}{\mddefault}{\updefault}$p_3$}}}}}
\put(600,636){\makebox(0,0)[lb]{\smash{{{\SetFigFont{12}{14.4}{\rmdefault}{\mddefault}{\updefault}$D_{\{L\}}$}}}}}
\put(2400,636){\makebox(0,0)[lb]{\smash{{{\SetFigFont{12}{14.4}{\rmdefault}{\mddefault}{\updefault}$D_{\{\}}$}}}}}
\put(4050,636){\makebox(0,0)[lb]{\smash{{{\SetFigFont{12}{14.4}{\rmdefault}{\mddefault}{\updefault}$D_{\{R\}}$}}}}}
\put(3000,36){\makebox(0,0)[lb]{\smash{{{\SetFigFont{8}{9.6}{\rmdefault}{\mddefault}{\updefault}geometrically irrelevant}}}}}
\put(675,1986){\makebox(0,0)[lb]{\smash{{{\SetFigFont{8}{9.6}{\rmdefault}{\mddefault}{\updefault}$\ell_{12}$}}}}}
\put(300,1761){\makebox(0,0)[lb]{\smash{{{\SetFigFont{8}{9.6}{\rmdefault}{\mddefault}{\updefault}$\ell$}}}}}
\put(2400,1011){\makebox(0,0)[lb]{\smash{{{\SetFigFont{8}{9.6}{\rmdefault}{\mddefault}{\updefault}$f_{123}$}}}}}
\put(2475,1686){\makebox(0,0)[lb]{\smash{{{\SetFigFont{8}{9.6}{\rmdefault}{\mddefault}{\updefault}$\ell_{12}$}}}}}
\put(2250,1536){\makebox(0,0)[lb]{\smash{{{\SetFigFont{8}{9.6}{\rmdefault}{\mddefault}{\updefault}$p_2$}}}}}
\put(2775,1536){\makebox(0,0)[lb]{\smash{{{\SetFigFont{8}{9.6}{\rmdefault}{\mddefault}{\updefault}$\ell_{23}$}}}}}
\put(4200,1011){\makebox(0,0)[lb]{\smash{{{\SetFigFont{8}{9.6}{\rmdefault}{\mddefault}{\updefault}$f_{123}$}}}}}
\put(4050,1611){\makebox(0,0)[lb]{\smash{{{\SetFigFont{8}{9.6}{\rmdefault}{\mddefault}{\updefault}$\ell_{12}=\ell_{23}$}}}}}
\put(4200,1311){\makebox(0,0)[lb]{\smash{{{\SetFigFont{8}{9.6}{\rmdefault}{\mddefault}{\updefault}$p_2$}}}}}
\end{picture}
}\end{center}
\caption{Degenerations of Figure~\ref{lubbock2} (planes
are denoted by arcs; a plane contains a line if the arc meets
the line)
\lremind{lubbock3}}
\label{lubbock3}
\end{figure}
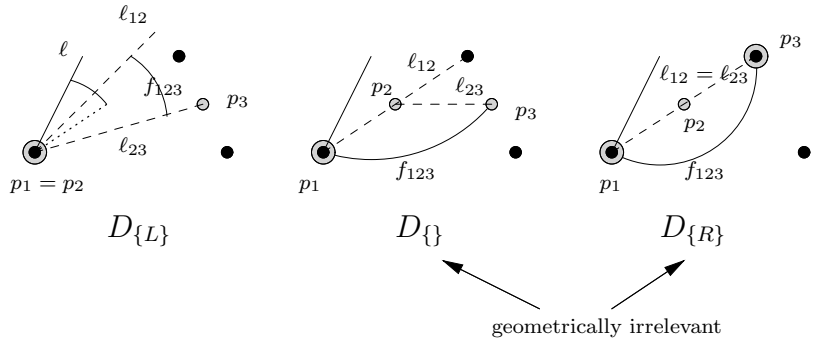

The first corresponds to the subset $S=\{L \}$, and is 
geometrically  {\em relevant}.  For a fixed point of $X_{\bn}$,
this locus visibly has dimension 2, one for the choice of $l_{12}$
through $p_1=p_2=b_1=b_2$ in the plane spanned by $\ell$ and $b_3$, and one for the
choice of $p_3$ on the line $b_3 b_4$.

The second corresponds to the subset $S=\{  \}$, and is
geometrically  {\em irrelevant}.   For a fixed point of $X_{\bn}$,
this locus visibly has dimension 2, one for the choice of $p_2$ on the
line $b_1b_3$, and one for the choice of 
$p_3$ on the line $b_3 b_4$.  However, the image in $G(3,4)= \mathbb{G}(2,3)$ has dimension
1, as the plane $f_{123}$ depends only on the choice of $p_3$.
This observation is a special case of Proposition~\ref{redmond}.

The third corresponds to the subset $\{ R \}$, and is also
geometrically {\em irrelevant}.  For a fixed point of $X_{\bn}$,
this locus visibly has dimension 2, one for the choice of $p_2$ on the
line $b_1b_3$, and one for the choice of the plane $f_{123}$
containing the line $b_1b_3$.  Again, the image in $G(3,4)$ has
dimension 1.  This observation is a special case of
Proposition~\ref{tacoma}.

Note that there is no component corresponding to $S=\{ L, R \}$, as as
the locus $Q_{\ell = 2}$ on $\PF(\poset_{\circ})$ is empty (as $\dim
\mathbf{a} = \dim \mathbf{a'} = 1$, so $\ell := \dim V_{\mathbf{a}}
\cap V_{\mathbf{a'}} - \dim V_{\inf(\mathbf{a},\mathbf{a'})} \leq 1-0
= 1$).

\tpoint{Proposition} {\em If $S \neq \{ \}, \{ \text{top left good quad.}\}$, then $D_S$ is geometrically  irrelevant.} \label{tacoma}
\lremind{tacoma}

\bpf
Note that the projection $$
\Cl_{\PF(\poset_{\circ}) \times (X_{\bullet} \cup X_{\bn})} X_{\cb}
\rightarrow
\Cl_{G(k,n) \times (X_{\bullet} \cup X_{\bn})} X_{\cb}$$
 involves forgetting
all elements of $\poset_{\circ}$ but the maximum, and that 
its construction 
(given in Sect.~\ref{provincetown}) depends only
on the vector spaces on the northeast and southeast borders
of $\poset_{\circ}$.  
Fix a general
point of $D^Z_S$ (corresponding to a configuration of
subspaces  $\{ V_{\mathbf{c}} \subset K^n : \mathbf{c} \in \poset_\circ \}$).
It will suffice to
produce a one-parameter family through 
this point 
that preserves the spaces 
in the northeast and southwest borders.  
All possibilities for $S$ fall into at least one of the
following three cases.

{\em Case 1:  $S$ contains a left good quadrilateral that is not
the top left good quadrilateral.}  
Name the elements of $\poset_{\circ}$ as in Figure~\ref{dallas1}
(showing the quadrilateral in question).  
Then a one-parameter family corresponds to fixing
all $V_{\mathbf{m}}$ ($\mathbf{m} \in \poset_{\circ}$) except $V_{\mathbf{d}}$, 
and letting $V_{\mathbf{d}}$ move freely in $\proj( V_{\mathbf{b}}/V_{\mathbf{e}})
\cong 
\proj^1$.

(Note that this would fail if the left good quadrilateral
were at the top.  In this case, $p_{\cB}(\mathbf{d}) = \x$, but
the $V_{\mathbf{d}}$ corresponding to a general element of $\proj(V_{\mathbf{b}}/V_{\mathbf{e}})$
would lie in $V_{\y} = V_{\yp}$, and not necessarily in $V_{\x}$.)

\begin{figure}[ht]
\begin{center}\setlength{\unitlength}{0.00083333in}
\begingroup\makeatletter\ifx\SetFigFont\undefined%
\gdef\SetFigFont#1#2#3#4#5{%
  \reset@font\fontsize{#1}{#2pt}%
  \fontfamily{#3}\fontseries{#4}\fontshape{#5}%
  \selectfont}%
\fi\endgroup%
{\renewcommand{\dashlinestretch}{30}
\begin{picture}(2850,1836)(0,-10)
\put(900,1662){\ellipse{150}{150}}
\put(1350,912){\ellipse{150}{150}}
\put(1500,312){\ellipse{150}{150}}
\path(900,1587)(1350,987)
\path(1350,837)(1500,387)
\path(1575,312)(2550,312)
\path(1425,912)(2550,912)(2550,312)
\blacken\path(1950,1212)(1988,1212)(1988,12)
	(1950,12)(1950,1212)
\path(1950,1212)(1988,1212)(1988,12)
	(1950,12)(1950,1212)
\path(1650,462)(2475,837)
\path(1650,462)(2475,837)
\path(1650,387)(2475,762)
\path(1650,387)(2475,762)
\put(2625,912){\makebox(0,0)[lb]{\smash{{{\SetFigFont{8}{9.6}{\rmdefault}{\mddefault}{\updefault}$\mathbf{c}$}}}}}
\put(0,537){\makebox(0,0)[lb]{\smash{{{\SetFigFont{8}{9.6}{\rmdefault}{\mddefault}{\updefault}$p_{\circ \Box}(\mathbf{b})=p_{\circ \Box}(\mathbf{d})=\mathbf{y}$}}}}}
\put(1200,1587){\makebox(0,0)[lb]{\smash{{{\SetFigFont{8}{9.6}{\rmdefault}{\mddefault}{\updefault}$p_{\circ \Box}(\mathbf{e})=\mathbf{x}\text{ or }\mathbf{y}$}}}}}
\put(2850,762){\makebox(0,0)[lb]{\smash{{{\SetFigFont{8}{9.6}{\rmdefault}{\mddefault}{\updefault}$p_{\circ \Box}(\mathbf{c})=\mathbf{z}$}}}}}
\put(1350,162){\makebox(0,0)[lb]{\smash{{{\SetFigFont{8}{9.6}{\rmdefault}{\mddefault}{\updefault}$\mathbf{b}$}}}}}
\put(1125,912){\makebox(0,0)[lb]{\smash{{{\SetFigFont{8}{9.6}{\rmdefault}{\mddefault}{\updefault}$\mathbf{d}$}}}}}
\put(675,1737){\makebox(0,0)[lb]{\smash{{{\SetFigFont{8}{9.6}{\rmdefault}{\mddefault}{\updefault}$\mathbf{e}$}}}}}
\end{picture}
}\end{center}
\caption{Case 1 of Proposition~\ref{tacoma}\lremind{dallas1}}
\label{dallas1}
\end{figure}

{\em Case 2: $S$ contains a right good quadrilateral $q$, and no
  quadrilaterals in the column directly to the left of $q$, and at
  least as high as $q$.}  Name the elements of $\poset_{\circ}$ as in
Figure~\ref{dallas2}. 
 The desired one-parameter family 
corresponds to fixing all the vector spaces $V_{\mathbf{m}}$ ($\mathbf{m} \in \poset_{\circ}$)
except for
$V_{\mathbf{d}}$, $V_{\mathbf{g_1}}$, \dots, $V_{\mathbf{g_m}}$, and choosing $V_{\mathbf{d}}$ from an open
set of $\proj ( V_{\mathbf{b}} / V_{\mathbf{e}})  \cong \proj^1$, such that $V_{\mathbf{g_i}} := V_{\mathbf{d}}
\cap V_{\mathbf{h_i}}$ has dimension $\dim \mathbf{g_i}$ ($1 \leq i \leq m$).  Note that
$V_{\mathbf{f_i}}$ is automatically contained in $V_{\mathbf{g_i}}$ as $V_{\mathbf{g_i}} = V_{\mathbf{d}}
\cap V_{\mathbf{h_i}}$ contains $V_{\mathbf{e}} \cap V_{\mathbf{h_i}}=V_{\mathbf{f_i}}$ by construction.

\begin{figure}[ht]
\begin{center}\setlength{\unitlength}{0.00083333in}
\begingroup\makeatletter\ifx\SetFigFont\undefined%
\gdef\SetFigFont#1#2#3#4#5{%
  \reset@font\fontsize{#1}{#2pt}%
  \fontfamily{#3}\fontseries{#4}\fontshape{#5}%
  \selectfont}%
\fi\endgroup%
{\renewcommand{\dashlinestretch}{30}
\begin{picture}(2862,2155)(0,-10)
\put(450,1381){\ellipse{76}{76}}
\put(1050,1681){\ellipse{76}{76}}
\path(450,1343)(450,481)
\path(1088,1681)(1650,1681)(1650,181)(1050,181)
\path(488,1381)(1650,1381)
\path(150,1981)(450,1418)
\path(150,1981)(1013,1681)
\path(1088,256)(1613,443)
\path(1088,218)(1613,406)
\path(450,1081)(1650,1081)
\path(450,781)(1650,781)
\path(1050,1643)(1050,181)
\path(450,481)(1650,481)
\blacken\path(750,1531)(2850,1531)(2850,1493)
	(750,1493)(750,1531)
\path(750,1531)(2850,1531)(2850,1493)
	(750,1493)(750,1531)
\put(1463,256){\makebox(0,0)[lb]{\smash{{{\SetFigFont{8}{9.6}{\rmdefault}{\mddefault}{\updefault}$q$}}}}}
\put(825,1418){\makebox(0,0)[lb]{\smash{{{\SetFigFont{8}{9.6}{\rmdefault}{\mddefault}{\updefault}$\mathbf{g_{m-1}}$}}}}}
\put(225,1081){\makebox(0,0)[lb]{\smash{{{\SetFigFont{8}{9.6}{\rmdefault}{\mddefault}{\updefault}$\vdots$}}}}}
\put(150,781){\makebox(0,0)[lb]{\smash{{{\SetFigFont{8}{9.6}{\rmdefault}{\mddefault}{\updefault}$\mathbf{f_1}$}}}}}
\put(150,406){\makebox(0,0)[lb]{\smash{{{\SetFigFont{8}{9.6}{\rmdefault}{\mddefault}{\updefault}$\mathbf{e}$}}}}}
\put(900,31){\makebox(0,0)[lb]{\smash{{{\SetFigFont{8}{9.6}{\rmdefault}{\mddefault}{\updefault}$\mathbf{b}$}}}}}
\put(900,331){\makebox(0,0)[lb]{\smash{{{\SetFigFont{8}{9.6}{\rmdefault}{\mddefault}{\updefault}$\mathbf{d}$}}}}}
\put(975,1831){\makebox(0,0)[lb]{\smash{{{\SetFigFont{8}{9.6}{\rmdefault}{\mddefault}{\updefault}$\mathbf{g_m}$}}}}}
\put(0,2056){\makebox(0,0)[lb]{\smash{{{\SetFigFont{8}{9.6}{\rmdefault}{\mddefault}{\updefault}$\mathbf{f_m}$}}}}}
\put(900,856){\makebox(0,0)[lb]{\smash{{{\SetFigFont{8}{9.6}{\rmdefault}{\mddefault}{\updefault}$\mathbf{g_1}$}}}}}
\put(1725,1756){\makebox(0,0)[lb]{\smash{{{\SetFigFont{8}{9.6}{\rmdefault}{\mddefault}{\updefault}$\mathbf{h_m}$}}}}}
\put(1725,1381){\makebox(0,0)[lb]{\smash{{{\SetFigFont{8}{9.6}{\rmdefault}{\mddefault}{\updefault}$\mathbf{h_{m-1}}$}}}}}
\put(1725,781){\makebox(0,0)[lb]{\smash{{{\SetFigFont{8}{9.6}{\rmdefault}{\mddefault}{\updefault}$\mathbf{h_1}$}}}}}
\put(1725,481){\makebox(0,0)[lb]{\smash{{{\SetFigFont{8}{9.6}{\rmdefault}{\mddefault}{\updefault}$\mathbf{c}$}}}}}
\put(75,1381){\makebox(0,0)[lb]{\smash{{{\SetFigFont{8}{9.6}{\rmdefault}{\mddefault}{\updefault}$\mathbf{f_{m-1}}$}}}}}
\put(1800,1081){\makebox(0,0)[lb]{\smash{{{\SetFigFont{8}{9.6}{\rmdefault}{\mddefault}{\updefault}$\vdots$}}}}}
\put(2700,1681){\makebox(0,0)[lb]{\smash{{{\SetFigFont{12}{14.4}{\rmdefault}{\mddefault}{\updefault}$\mathbf{y'}$}}}}}
\put(2700,181){\makebox(0,0)[lb]{\smash{{{\SetFigFont{12}{14.4}{\rmdefault}{\mddefault}{\updefault}$\mathbf{z}$}}}}}
\end{picture}
}\end{center}
\caption{Case 2 of Proposition~\ref{tacoma}\lremind{dallas2}}
\label{dallas2}
\end{figure}

{\em Case 3: $S$ contains a right good quadrilateral $q$, and another
good 
  quadrilateral in the column directly to the left of $q$, at
  least as high as $q$.}  Name the elements of $\poset_{\circ}$ as in
Figure~\ref{dallas3}.   Then the argument from the Case 2 applies
verbatim.
\epf

\begin{figure}[ht]
\begin{center}\setlength{\unitlength}{0.00083333in}
\begingroup\makeatletter\ifx\SetFigFont\undefined%
\gdef\SetFigFont#1#2#3#4#5{%
  \reset@font\fontsize{#1}{#2pt}%
  \fontfamily{#3}\fontseries{#4}\fontshape{#5}%
  \selectfont}%
\fi\endgroup%
{\renewcommand{\dashlinestretch}{30}
\begin{picture}(1575,1633)(0,-10)
\path(225,1606)(225,406)
\path(225,706)(1425,706)
\path(863,181)(1388,368)
\path(863,143)(1388,331)
\path(225,1006)(1425,1006)
\path(825,106)(1425,106)(1425,1306)
\path(225,1606)(825,1606)
\path(263,1343)(788,1531)
\path(263,1381)(788,1568)
\path(225,406)(1425,406)
\path(825,1606)(825,106)
\path(225,1306)(1425,1306)
\put(1238,181){\makebox(0,0)[lb]{\smash{{{\SetFigFont{8}{9.6}{\rmdefault}{\mddefault}{\updefault}$q$}}}}}
\put(0,1306){\makebox(0,0)[lb]{\smash{{{\SetFigFont{8}{9.6}{\rmdefault}{\mddefault}{\updefault}$\mathbf{f_{m}}$}}}}}
\put(675,781){\makebox(0,0)[lb]{\smash{{{\SetFigFont{8}{9.6}{\rmdefault}{\mddefault}{\updefault}$\mathbf{g_1}$}}}}}
\put(75,331){\makebox(0,0)[lb]{\smash{{{\SetFigFont{8}{9.6}{\rmdefault}{\mddefault}{\updefault}$\mathbf{e}$}}}}}
\put(75,706){\makebox(0,0)[lb]{\smash{{{\SetFigFont{8}{9.6}{\rmdefault}{\mddefault}{\updefault}$\mathbf{f_1}$}}}}}
\put(75,1006){\makebox(0,0)[lb]{\smash{{{\SetFigFont{8}{9.6}{\rmdefault}{\mddefault}{\updefault}$\vdots$}}}}}
\put(675,1381){\makebox(0,0)[lb]{\smash{{{\SetFigFont{8}{9.6}{\rmdefault}{\mddefault}{\updefault}$\mathbf{g_{m}}$}}}}}
\put(675,31){\makebox(0,0)[lb]{\smash{{{\SetFigFont{8}{9.6}{\rmdefault}{\mddefault}{\updefault}$\mathbf{b}$}}}}}
\put(1500,406){\makebox(0,0)[lb]{\smash{{{\SetFigFont{8}{9.6}{\rmdefault}{\mddefault}{\updefault}$\mathbf{c}$}}}}}
\put(1500,706){\makebox(0,0)[lb]{\smash{{{\SetFigFont{8}{9.6}{\rmdefault}{\mddefault}{\updefault}$\mathbf{h_1}$}}}}}
\put(1575,1006){\makebox(0,0)[lb]{\smash{{{\SetFigFont{8}{9.6}{\rmdefault}{\mddefault}{\updefault}$\vdots$}}}}}
\put(1500,1306){\makebox(0,0)[lb]{\smash{{{\SetFigFont{8}{9.6}{\rmdefault}{\mddefault}{\updefault}$\mathbf{h_{m}}$}}}}}
\put(675,331){\makebox(0,0)[lb]{\smash{{{\SetFigFont{8}{9.6}{\rmdefault}{\mddefault}{\updefault}$\mathbf{d}$}}}}}
\end{picture}
}\end{center}
\caption{Case 3 of Proposition~\ref{tacoma}\lremind{dallas3}}
\label{dallas3}
\end{figure}

\tpoint{Proposition} {\em If $S = \{ \}$ and (i) the white checker in
  the critical row is in the descending checker's square, or (ii) the
  white checker in the critical diagonal is in the rising checker's
  square, then $D_S$ is geometrically irrelevant.} \label{redmond}
\lremind{redmond}

\bpf  
We deal with case (ii); case (i) is essentially identical and hence
omitted.

As in the proof of Proposition~\ref{tacoma}, we 
note that 
the construction 
of Section~\ref{provincetown} depended only
on the vector spaces on the northeast and southeast borders
of $\poset_{\circ}$, but {\em not} the vector space
on the northeast border in the critical row. 
Hence it will suffice to
produce a one-parameter family through a general
point of $D^Z_S$  that preserves the spaces 
in the northeast and southwest borders, except 
the element of the critical row. 

\begin{figure}[ht]
\begin{center}\setlength{\unitlength}{0.00083333in}
\begingroup\makeatletter\ifx\SetFigFont\undefined%
\gdef\SetFigFont#1#2#3#4#5{%
  \reset@font\fontsize{#1}{#2pt}%
  \fontfamily{#3}\fontseries{#4}\fontshape{#5}%
  \selectfont}%
\fi\endgroup%
{\renewcommand{\dashlinestretch}{30}
\begin{picture}(3249,1207)(0,-10)
\put(312,1033){\ellipse{76}{76}}
\put(987,658){\ellipse{76}{76}}
\put(387,358){\ellipse{76}{76}}
\path(425,358)(2787,358)
\path(1025,658)(2787,658)
\path(1587,958)(2787,958)
\path(2787,958)(2787,358)
\path(2187,958)(2187,358)
\path(1587,958)(1587,358)
\path(987,620)(987,358)
\path(312,995)(387,395)
\path(350,1033)(950,658)
\blacken\path(687,808)(706,808)(706,58)
	(687,58)(687,808)
\path(687,808)(706,808)(706,58)
	(687,58)(687,808)
\blacken\path(12,508)(3237,508)(3237,489)
	(12,489)(12,508)
\path(12,508)(3237,508)(3237,489)
	(12,489)(12,508)
\blacken\path(687,808)(3237,808)(3237,789)
	(687,789)(687,808)
\path(687,808)(3237,808)(3237,789)
	(687,789)(687,808)
\put(1550,995){\makebox(0,0)[lb]{\smash{{{\SetFigFont{8}{9.6}{\rmdefault}{\mddefault}{\updefault}$\mathbf{f_1}$}}}}}
\put(2150,995){\makebox(0,0)[lb]{\smash{{{\SetFigFont{5}{6.0}{\rmdefault}{\mddefault}{\updefault}$\cdots$}}}}}
\put(2750,995){\makebox(0,0)[lb]{\smash{{{\SetFigFont{8}{9.6}{\rmdefault}{\mddefault}{\updefault}$\mathbf{f_m}$}}}}}
\put(2825,695){\makebox(0,0)[lb]{\smash{{{\SetFigFont{8}{9.6}{\rmdefault}{\mddefault}{\updefault}$\mathbf{g_m}$}}}}}
\put(1625,695){\makebox(0,0)[lb]{\smash{{{\SetFigFont{8}{9.6}{\rmdefault}{\mddefault}{\updefault}$\mathbf{g_1}$}}}}}
\put(12,808){\makebox(0,0)[lb]{\smash{{{\SetFigFont{12}{14.4}{\rmdefault}{\mddefault}{\updefault}$\mathbf{x}$}}}}}
\put(12,58){\makebox(0,0)[lb]{\smash{{{\SetFigFont{12}{14.4}{\rmdefault}{\mddefault}{\updefault}$\mathbf{y}$}}}}}
\put(3087,583){\makebox(0,0)[lb]{\smash{{{\SetFigFont{12}{14.4}{\rmdefault}{\mddefault}{\updefault}$\mathbf{y'}$}}}}}
\put(3087,58){\makebox(0,0)[lb]{\smash{{{\SetFigFont{12}{14.4}{\rmdefault}{\mddefault}{\updefault}$\mathbf{z}$}}}}}
\put(1550,245){\makebox(0,0)[lb]{\smash{{{\SetFigFont{8}{9.6}{\rmdefault}{\mddefault}{\updefault}$\mathbf{h_1}$}}}}}
\put(950,245){\makebox(0,0)[lb]{\smash{{{\SetFigFont{8}{9.6}{\rmdefault}{\mddefault}{\updefault}$\mathbf{b}$}}}}}
\put(2150,245){\makebox(0,0)[lb]{\smash{{{\SetFigFont{8}{9.6}{\rmdefault}{\mddefault}{\updefault}$\cdots$}}}}}
\put(2750,245){\makebox(0,0)[lb]{\smash{{{\SetFigFont{8}{9.6}{\rmdefault}{\mddefault}{\updefault}$\mathbf{h_m}$}}}}}
\put(275,245){\makebox(0,0)[lb]{\smash{{{\SetFigFont{8}{9.6}{\rmdefault}{\mddefault}{\updefault}$\mathbf{c}$}}}}}
\put(387,1108){\makebox(0,0)[lb]{\smash{{{\SetFigFont{8}{9.6}{\rmdefault}{\mddefault}{\updefault}$\mathbf{e}$}}}}}
\put(1017,680){\makebox(0,0)[lb]{\smash{{{\SetFigFont{8}{9.6}{\rmdefault}{\mddefault}{\updefault}$\mathbf{d}$}}}}}
\end{picture}
}\end{center}
\caption{\lremind{dallas4}}
\label{dallas4}
\end{figure}

Name the elements of $\poset_{\circ}$ as in
Figure~\ref{dallas4}. 
As the point of $D^Z_S$  maps to the
dense open stratum of $\PF(\poset_{\circ})$, we have $V_{\mathbf{h_i}} = V_{\mathbf{c}}+V_{\mathbf{g_i}}$;
as $V_{\mathbf{c}} \subset V_{\y} = V_{\yp}$ and $V_{\mathbf{g_i}} \subset V_{\yp}$, 
we have $V_{\mathbf{h_i}} \subset
V_{\y}$ as well.
 The desired one-parameter family 
corresponds to fixing all the vector spaces $V_{\mathbf{l}}$ ($\mathbf{l} \in \poset_{\circ}$)
except for
$V_{\mathbf{d}}$, $V_{\mathbf{g_1}}$, \dots, $V_{\mathbf{g_m}}$, and choosing $V_{\mathbf{d}}$ from an open
set of $\proj ( V_{\mathbf{b}} / V_{\mathbf{e}})  \cong \proj^1$, such that $V_{\mathbf{g_i}} := V_{\mathbf{d}}
+ V_{\mathbf{f_i}}$ has dimension $\dim \mathbf{g_i}$ ($1 \leq i \leq m$). 
\epf

\bpoint{Multiplicity 1} 
\label{amherst} \lremind{amherst} 
We have shown that at most one or two of
the components of $D'_X$ are geometrically relevant.  We now show that
these one or two components  appear with multiplicity 1 in $\Cl_{\PF(\poset_{\circ}) \times
  \left( \Xbb \right)} X_{\cb}$.    By the top row of (\ref{newton}) or the middle row of (\ref{miami}), it
    suffices to prove that the analogous divisors on $Z$ appear with
    multiplicity 1.

\tpoint{Proposition}
{\em
\begin{enumerate} \item[(a)]  The divisor $D^Z_{ \{ \}}$ is contained in the
pullback of $D_\Box$ with multiplicity 1.
\item[(b)] 
The divisor $D^Z_{ \{\text{top left good quad.} \}}$ is contained in the
pullback of $D_\Box$ with multiplicity 1 if
there is
no blocker (i.e.  there is a top left good quadrilateral). 
\end{enumerate}}

\bpf
(a)  It suffices to verify the result
in the preimage of the dense open stratum of $\PF(\poset_{\circ})$;
it then suffices to verify the result for a fixed point of $\PF(\poset_{\circ})$.
Here
$$
V_{\mathbf{a}} \cap V_{\mathbf{a'}} = V_{\mathbf{a}} \cap V_{{\mathbf{a''}}} = V_{\inf(\mathbf{a},\mathbf{a'})}.
$$
The Proposition then reduces to a straightforward (hence omitted) local calculation
about the locus in $\PF(\pBox)$ such that $V_{\x}$ contains $V_{{\mathbf{a''}}}$, $V_{\y}$
contains $V_{\mathbf{a}}$, and $V_{\yp}$ contains $V_{\mathbf{a'}}$ (where $V_{\mathbf{a}}$, $V_{\mathbf{a'}}$, $V_{{\mathbf{a''}}}$
are fixed).

(b) It suffices to verify the result in the preimage of the union of
the dense open stratum of $\PF(\poset_{\circ})$ and the (divisorial) open stratum of
$\PF(\poset_{\circ})$ corresponding to $\{\text{top left good
    quad.} \}$. This union admits a $\proj^1$-fibration to the
divisorial stratum of $\PF(\poset_{\circ})$, and it suffices to verify
the result for a fixed fiber as follows.  Name the elements of
$\poset_{\circ}$ as in Figure~\ref{dallas4}.  (Note that $\mathbf{a'}=\mathbf{g_m}$,
${\mathbf{a''}} = \mathbf{f_m}$, $\mathbf{e}=\inf(\mathbf{a},\mathbf{a'})$.)  Hold fixed all the spaces of
$\poset_{\circ}$ except those in the critical row $\{ V_{\mathbf{d}}, V_{\mathbf{g_1}},
\dots, V_{\mathbf{g_m}} \}$.  The space $V_{\mathbf{d}}$ corresponding to the white
checker varies in $\proj( V_{\mathbf{b}}/V_{\mathbf{e}})$, and $V_{\mathbf{g_i}}$ is defined to be
$V_{\mathbf{d}} + V_{\mathbf{f_i}}$.  The Proposition then reduces to a straightforward (hence omitted) 
local calculation about the space parametrizing  $(V_{\x}, V_{\y}, V_{\yp}, V_{\mathbf{d}})$ such that $V_{\mathbf{d}}$ varies in the pencil,  $V_{\x}$
contains the fixed space $V_{\mathbf{f_m}}$, $V_{\y}$ contains the fixed space $V_{\mathbf{a}}$, and $V_{\yp}$ contains $V_{\mathbf{g_m}} = V_{\mathbf{d}} + V_{\mathbf{f_m}}$.
\epf

Finally, we show that the map from $D_S$ to $\Cl_{G(k,n) \times \left( \Xbb \right)} X_{\cb}$ is birational to $\overline{X}_{\cbn}$ or
$\overline{X}_{\cbs}$; this will conclude the proof of the
\LR.

\tpoint{Proposition} \label{birational}\lremind{birational}
{\em  (a) $D_{ \{ \} }$ is birational to $\overline{X}_{\cbn}$.
 (b) $D_{ \{ \text{top left good quad.}\} }$ is birational to $\overline{X}_{\cbs}$.}

\bpf
(a) The open subset of $D_{ \{ \} }$ contained in the preimage of the dense open stratum
of $\PF(\poset_{\circ})$ is isomorphic to $X_{\cbn}$.
{\em Sketch of reason:}  The subspaces parametrized by the dense open 
stratum are determined by the subspaces corresponding to the white checkers.
(For $\mathbf{c} \in \poset_\circ$, $V_{\mathbf{c}} = \operatorname{span}( V_{\mathbf{w}})$, where
$\mathbf{w}$ runs over the white checkers dominated by $\mathbf{c}$.)
As in the construction of Section~\ref{provincetown}, choose
(i) $S_i$ to contain those $V_{\mathbf{w}}$ where $\mathbf{w}$ is
a white checker in column at most $i$, and to contain no $V_{\mathbf{w}}$
where $\mathbf{w}$ is a white checker in column greater than $i$, and (ii)
$T_i$ analogously.  This gives a subset of the ``divisor  at $\infty$''
of $\Cl_{\PF(\poset_\circ) \times \left( \Xbb \right) } X_{\circ \bB}$
whose image in
$\Cl_{\PF(\poset_{\circ}) \times \left( \Xbb \right) } X_{\cb}$
is $X_{\cbn}$.

(b)
The open subset of $D_{ \{ \text{top left good quad.} \}}$ contained in the preimage of the divisorial
open stratum of $\PF(\poset_{\circ})$ corresponding to $\{ \text{top left good quad.} \}$ is isomorphic to 
$X_{\cbs}$.
{\em Sketch of reason:}  The subspaces parametrized by the divisorial open 
stratum are determined by the subspaces corresponding to
\begin{enumerate}
\item[(i)] the white checkers of $\circ$, except the white
checkers at the opposite corners of the top left good quadrilateral
(call them $\mathbf{w_1} = (r_1,c_1)$, 
$\mathbf{w_2} = (r_2,c_2)$, with $r_1>r_2$, $c_1<c_2$),
\item[(ii)] the vector space $V_{\mathbf{w_1}} = V_{\mathbf{w_2}}$, and
\item[(iii)] $V_{\sup(\mathbf{w_1},\mathbf{w_2})}$ (corresponding
to the lower right corner of the top left good quadrilateral).
\end{enumerate}
We can interpret these subspaces as follows.  Identify $\mathbf{w_1}$
and $\mathbf{w_2}$ of $\poset_{\circ}$ to obtain the poset $\poset_{\cs}$.
The spaces described above correspond to those elements of $\poset_{\cs}$
with exactly  one edge to a lower-dimensional element (i.e. those $\mathbf{w}$
such that there is a unique covering relation $\mathbf{w'} \prec \mathbf{w}$).
They also correspond to the white checkers of $\cs$.
As in (a), the image in 
$\Cl_{\PF(\poset_{\circ}) \times \left( \Xbb \right) } X_{\cb}$
is (naturally isomorphic to) $X_{\cbs}$.  (The white checker of $\cs$ corresponding
$V_{\mathbf{w_1}}=V_{\mathbf{w_2}}$ must lie in the intersection
of the spaces where $V_{\mathbf{w_1}}$ and $V_{\mathbf{w_2}}$ were required
to lie, in the definition of $X_{\cb}$. 
Hence this checker lies in $(r_2, c_1)$.  The white checker of $\cbs$
corresponding to $V_{\sup(\mathbf{w_1},\mathbf{w_2})}$ is required
to lie in $(r_1,c_2)$, as it was in $X_{\cb}$.)
\epf

\section{Appendix:  The bijection between checker games and puzzles (with
A. Knutson)}
\label{checkerpuzzle}
\lremind{checkerpuzzle}
Fix $k$ and $n$. 
We fill in a puzzle with given inputs, one row of triangles at a time,
from left to right.   Row $m$ consists of those triangles between 
the $m^{\text{th}}$ edges from the top on the sides of the triangle.

The placement of  vertical rhombi may cause parts of
subsequent rows to be filled; call these {\em teeth}.  The
$m^{\text{th}}$ row ($1 < m \leq n$) corresponds to the part of the
checker game where the black checker in the $m^{\text{th}}$ column is
descending.  The possible choices for filling in puzzle pieces
correspond to the possible choices of next moves in the checker game;
this will give the bijection.

We now describe an injection from checker games to puzzles; to each
checker game we will associate a puzzle.  As both count
Littlewood-Richardson coefficients, this injection must be a
bijection.  Alternatively, to show that this is a bijection, one can
instead show that there are no puzzles not accounted for here.  For
example, one can show easily that there are no puzzles if the
checker game predicts there
shouldn't be (i.e. if the sets are $a_1 < \cdots < a_k$ and $b_1 <
\cdots < b_k$, and if $a_i + b_{k+1-i} \leq n$ for some $i$), by
focusing on a certain parallelogram-shaped region of the puzzle.
More generally, one should be able to show combinatorially that
if a partially-filled-in puzzle doesn't correspond to a valid checker
game (in progress), then there is no way to complete it.

\bpoint{Bijection of starting positions}  Fill in the top row of the
puzzle in the only way possible.  (As remarked earlier, the
translation to checkers will give an immediate criterion for there to
be no puzzles.)

\bpoint{The translation part-way through the checker game}
At each stage, the partially complete puzzle will look like Figure~\ref{pun}.
Any of $a$, $b$, and $c$ may be zero.   In the checker game, $a$, $b$, and
$c$ correspond to the the numbers shown in Figure~\ref{pub}.
The rows of the white checkers in the game are given by 
the edges of Figure~\ref{pus} --- a ``1'' indicates that there
is a white checker in that row.  The columns are given by the
edges of Figure~\ref{push}.  As the white checkers are mid-sort,
it turns out that this specifies their position completely.  See Figure~\ref{pummel}
for a more explicit description.

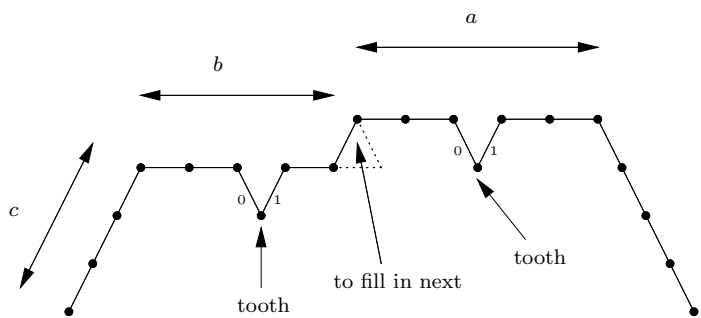
\begin{figure}[ht]
\begin{center}
\setlength{\unitlength}{0.00083333in}
\begingroup\makeatletter\ifx\SetFigFont\undefined%
\gdef\SetFigFont#1#2#3#4#5{%
  \reset@font\fontsize{#1}{#2pt}%
  \fontfamily{#3}\fontseries{#4}\fontshape{#5}%
  \selectfont}%
\fi\endgroup%
{\renewcommand{\dashlinestretch}{30}
\begin{picture}(4308,1932)(0,-10)
\put(375,33){\blacken\ellipse{50}{50}}
\put(375,33){\ellipse{50}{50}}
\put(525,333){\blacken\ellipse{50}{50}}
\put(525,333){\ellipse{50}{50}}
\put(675,633){\blacken\ellipse{50}{50}}
\put(675,633){\ellipse{50}{50}}
\put(825,933){\blacken\ellipse{50}{50}}
\put(825,933){\ellipse{50}{50}}
\put(1125,933){\blacken\ellipse{50}{50}}
\put(1125,933){\ellipse{50}{50}}
\put(1425,933){\blacken\ellipse{50}{50}}
\put(1425,933){\ellipse{50}{50}}
\put(1575,633){\blacken\ellipse{50}{50}}
\put(1575,633){\ellipse{50}{50}}
\put(1725,933){\blacken\ellipse{50}{50}}
\put(1725,933){\ellipse{50}{50}}
\put(2025,933){\blacken\ellipse{50}{50}}
\put(2025,933){\ellipse{50}{50}}
\put(2175,1233){\blacken\ellipse{50}{50}}
\put(2175,1233){\ellipse{50}{50}}
\put(2475,1233){\blacken\ellipse{50}{50}}
\put(2475,1233){\ellipse{50}{50}}
\put(2775,1233){\blacken\ellipse{50}{50}}
\put(2775,1233){\ellipse{50}{50}}
\put(2925,933){\blacken\ellipse{50}{50}}
\put(2925,933){\ellipse{50}{50}}
\put(3075,1233){\blacken\ellipse{50}{50}}
\put(3075,1233){\ellipse{50}{50}}
\put(3375,1233){\blacken\ellipse{50}{50}}
\put(3375,1233){\ellipse{50}{50}}
\put(3675,1233){\blacken\ellipse{50}{50}}
\put(3675,1233){\ellipse{50}{50}}
\put(3825,933){\blacken\ellipse{50}{50}}
\put(3825,933){\ellipse{50}{50}}
\put(3975,633){\blacken\ellipse{50}{50}}
\put(3975,633){\ellipse{50}{50}}
\put(4125,333){\blacken\ellipse{50}{50}}
\put(4125,333){\ellipse{50}{50}}
\put(4275,33){\blacken\ellipse{50}{50}}
\put(4275,33){\ellipse{50}{50}}
\path(375,33)(825,933)
\path(375,33)(825,933)
\path(825,933)(1425,933)
\path(825,933)(1425,933)
\path(1425,933)(1575,633)
\path(1425,933)(1575,633)
\path(1575,633)(1725,933)
\path(1575,633)(1725,933)
\path(1725,933)(2025,933)
\path(1725,933)(2025,933)
\path(2025,933)(2175,1233)
\path(2025,933)(2175,1233)
\path(2175,1233)(2775,1233)
\path(2175,1233)(2775,1233)
\path(2775,1233)(2925,933)
\path(2775,1233)(2925,933)
\path(2925,933)(3075,1233)
\path(2925,933)(3075,1233)
\path(3075,1233)(3675,1233)
\path(3075,1233)(3675,1233)
\path(3675,1233)(4275,33)
\path(3675,1233)(4275,33)
\blacken\path(2295.000,1713.000)(2175.000,1683.000)(2295.000,1653.000)(2295.000,1713.000)
\path(2175,1683)(3675,1683)
\path(2175,1683)(3675,1683)
\blacken\path(3555.000,1653.000)(3675.000,1683.000)(3555.000,1713.000)(3555.000,1653.000)
\path(2325,333)(2175,1083)
\path(2325,333)(2175,1083)
\blacken\path(2227.951,971.214)(2175.000,1083.000)(2169.117,959.447)(2227.951,971.214)
\blacken\path(498.167,962.252)(525.000,1083.000)(444.502,989.085)(498.167,962.252)
\path(525,1083)(75,183)
\path(525,1083)(75,183)
\blacken\path(101.833,303.748)(75.000,183.000)(155.498,276.915)(101.833,303.748)
\blacken\path(945.000,1413.000)(825.000,1383.000)(945.000,1353.000)(945.000,1413.000)
\path(825,1383)(2025,1383)
\path(825,1383)(2025,1383)
\blacken\path(1905.000,1353.000)(2025.000,1383.000)(1905.000,1413.000)(1905.000,1353.000)
\path(1575,183)(1575,558)
\blacken\path(1605.000,438.000)(1575.000,558.000)(1545.000,438.000)(1605.000,438.000)
\path(3225,483)(2925,858)
\blacken\path(3023.389,783.037)(2925.000,858.000)(2976.537,745.555)(3023.389,783.037)
\dottedline{45}(2025,933)(2325,933)(2175,1233)
\put(0,633){\makebox(0,0)[lb]{\smash{{{\SetFigFont{8}{9.6}{\rmdefault}{\mddefault}{\updefault}$c$}}}}}
\put(1275,1533){\makebox(0,0)[lb]{\smash{{{\SetFigFont{8}{9.6}{\rmdefault}{\mddefault}{\updefault}$b$}}}}}
\put(2850,1833){\makebox(0,0)[lb]{\smash{{{\SetFigFont{8}{9.6}{\rmdefault}{\mddefault}{\updefault}$a$}}}}}
\put(1425,708){\makebox(0,0)[lb]{\smash{{{\SetFigFont{5}{6.0}{\rmdefault}{\mddefault}{\updefault}$0$}}}}}
\put(1650,708){\makebox(0,0)[lb]{\smash{{{\SetFigFont{5}{6.0}{\rmdefault}{\mddefault}{\updefault}$1$}}}}}
\put(3000,1008){\makebox(0,0)[lb]{\smash{{{\SetFigFont{5}{6.0}{\rmdefault}{\mddefault}{\updefault}$1$}}}}}
\put(2775,1008){\makebox(0,0)[lb]{\smash{{{\SetFigFont{5}{6.0}{\rmdefault}{\mddefault}{\updefault}$0$}}}}}
\put(1425,33){\makebox(0,0)[lb]{\smash{{{\SetFigFont{8}{9.6}{\rmdefault}{\mddefault}{\updefault}tooth}}}}}
\put(2025,183){\makebox(0,0)[lb]{\smash{{{\SetFigFont{8}{9.6}{\rmdefault}{\mddefault}{\updefault}to fill in next}}}}}
\put(3150,333){\makebox(0,0)[lb]{\smash{{{\SetFigFont{8}{9.6}{\rmdefault}{\mddefault}{\updefault}tooth}}}}}
\end{picture}
}
\end{center}
\caption{The puzzle in the process of being filled \lremind{pun}}
\label{pun}
\end{figure}

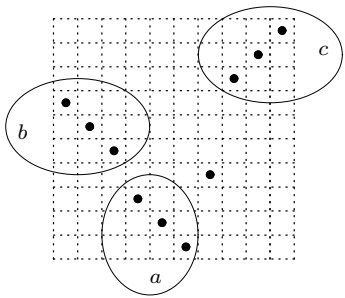
\begin{figure}[ht]
\begin{center}
\setlength{\unitlength}{0.00083333in}
\begingroup\makeatletter\ifx\SetFigFont\undefined%
\gdef\SetFigFont#1#2#3#4#5{%
  \reset@font\fontsize{#1}{#2pt}%
  \fontfamily{#3}\fontseries{#4}\fontshape{#5}%
  \selectfont}%
\fi\endgroup%
{\renewcommand{\dashlinestretch}{30}
\begin{picture}(2116,1830)(0,-10)
\put(1733,1658){\blacken\ellipse{50}{50}}
\put(1733,1658){\ellipse{50}{50}}
\put(1583,1508){\blacken\ellipse{50}{50}}
\put(1583,1508){\ellipse{50}{50}}
\put(1433,1358){\blacken\ellipse{50}{50}}
\put(1433,1358){\ellipse{50}{50}}
\put(1283,758){\blacken\ellipse{50}{50}}
\put(1283,758){\ellipse{50}{50}}
\put(1133,308){\blacken\ellipse{50}{50}}
\put(1133,308){\ellipse{50}{50}}
\put(983,458){\blacken\ellipse{50}{50}}
\put(983,458){\ellipse{50}{50}}
\put(833,608){\blacken\ellipse{50}{50}}
\put(833,608){\ellipse{50}{50}}
\put(683,908){\blacken\ellipse{50}{50}}
\put(683,908){\ellipse{50}{50}}
\put(533,1058){\blacken\ellipse{50}{50}}
\put(533,1058){\ellipse{50}{50}}
\put(383,1208){\blacken\ellipse{50}{50}}
\put(383,1208){\ellipse{50}{50}}
\put(908,383){\ellipse{600}{750}}
\put(458,1058){\ellipse{900}{600}}
\put(1658,1508){\ellipse{900}{600}}
\dottedline{45}(308,1733)(1808,1733)(1808,233)
	(308,233)(308,1733)
\dottedline{45}(1658,1733)(1658,233)
\dottedline{45}(1658,1733)(1658,233)
\dottedline{45}(1508,1733)(1508,233)
\dottedline{45}(1358,1733)(1358,233)
\dottedline{45}(1208,1733)(1208,233)
\dottedline{45}(1058,1733)(1058,233)
\dottedline{45}(908,1733)(908,233)
\dottedline{45}(758,1733)(758,233)
\dottedline{45}(608,1733)(608,233)
\dottedline{45}(458,1733)(458,233)
\dottedline{45}(308,1583)(1808,1583)
\dottedline{45}(308,1433)(1808,1433)
\dottedline{45}(308,1283)(1808,1283)
\dottedline{45}(308,1133)(1808,1133)
\dottedline{45}(308,983)(1808,983)
\dottedline{45}(308,833)(1808,833)
\dottedline{45}(308,683)(1808,683)
\dottedline{45}(308,533)(1808,533)
\dottedline{45}(308,383)(1808,383)
\put(83,983){\makebox(0,0)[lb]{\smash{{{\SetFigFont{8}{9.6}{\rmdefault}{\mddefault}{\updefault}$b$}}}}}
\put(908,83){\makebox(0,0)[lb]{\smash{{{\SetFigFont{8}{9.6}{\rmdefault}{\mddefault}{\updefault}$a$}}}}}
\put(1958,1508){\makebox(0,0)[lb]{\smash{{{\SetFigFont{8}{9.6}{\rmdefault}{\mddefault}{\updefault}$c$}}}}}
\end{picture}
}
\end{center}
\caption{The corresponding point in the checker game \lremind{pub}}
\label{pub}
\end{figure}

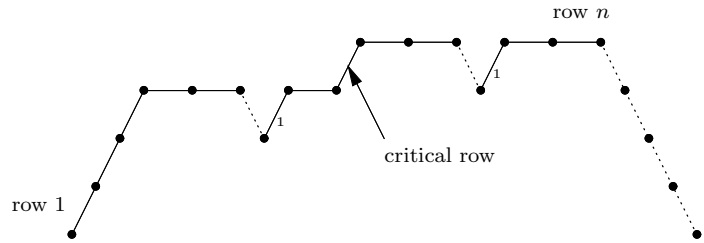
\begin{figure}[ht]
\begin{center}
\setlength{\unitlength}{0.00083333in}
\begingroup\makeatletter\ifx\SetFigFont\undefined%
\gdef\SetFigFont#1#2#3#4#5{%
  \reset@font\fontsize{#1}{#2pt}%
  \fontfamily{#3}\fontseries{#4}\fontshape{#5}%
  \selectfont}%
\fi\endgroup%
{\renewcommand{\dashlinestretch}{30}
\begin{picture}(4308,1482)(0,-10)
\put(375,33){\blacken\ellipse{50}{50}}
\put(375,33){\ellipse{50}{50}}
\put(525,333){\blacken\ellipse{50}{50}}
\put(525,333){\ellipse{50}{50}}
\put(675,633){\blacken\ellipse{50}{50}}
\put(675,633){\ellipse{50}{50}}
\put(825,933){\blacken\ellipse{50}{50}}
\put(825,933){\ellipse{50}{50}}
\put(1125,933){\blacken\ellipse{50}{50}}
\put(1125,933){\ellipse{50}{50}}
\put(1425,933){\blacken\ellipse{50}{50}}
\put(1425,933){\ellipse{50}{50}}
\put(1575,633){\blacken\ellipse{50}{50}}
\put(1575,633){\ellipse{50}{50}}
\put(1725,933){\blacken\ellipse{50}{50}}
\put(1725,933){\ellipse{50}{50}}
\put(2025,933){\blacken\ellipse{50}{50}}
\put(2025,933){\ellipse{50}{50}}
\put(2175,1233){\blacken\ellipse{50}{50}}
\put(2175,1233){\ellipse{50}{50}}
\put(2475,1233){\blacken\ellipse{50}{50}}
\put(2475,1233){\ellipse{50}{50}}
\put(2775,1233){\blacken\ellipse{50}{50}}
\put(2775,1233){\ellipse{50}{50}}
\put(2925,933){\blacken\ellipse{50}{50}}
\put(2925,933){\ellipse{50}{50}}
\put(3075,1233){\blacken\ellipse{50}{50}}
\put(3075,1233){\ellipse{50}{50}}
\put(3375,1233){\blacken\ellipse{50}{50}}
\put(3375,1233){\ellipse{50}{50}}
\put(3675,1233){\blacken\ellipse{50}{50}}
\put(3675,1233){\ellipse{50}{50}}
\put(3825,933){\blacken\ellipse{50}{50}}
\put(3825,933){\ellipse{50}{50}}
\put(3975,633){\blacken\ellipse{50}{50}}
\put(3975,633){\ellipse{50}{50}}
\put(4125,333){\blacken\ellipse{50}{50}}
\put(4125,333){\ellipse{50}{50}}
\put(4275,33){\blacken\ellipse{50}{50}}
\put(4275,33){\ellipse{50}{50}}
\path(375,33)(825,933)
\path(375,33)(825,933)
\path(825,933)(1425,933)
\path(825,933)(1425,933)
\dottedline{45}(1425,933)(1575,633)
\path(1575,633)(1725,933)
\path(1575,633)(1725,933)
\path(1725,933)(2025,933)
\path(1725,933)(2025,933)
\path(2025,933)(2175,1233)
\path(2025,933)(2175,1233)
\path(2175,1233)(2775,1233)
\path(2175,1233)(2775,1233)
\path(2925,933)(3075,1233)
\path(2925,933)(3075,1233)
\path(3075,1233)(3675,1233)
\path(3075,1233)(3675,1233)
\dottedline{45}(3675,1233)(4275,33)
\dottedline{45}(2775,1233)(2925,933)
\path(2325,633)(2100,1083)
\path(2325,633)(2100,1083)
\blacken\path(2180.498,989.085)(2100.000,1083.000)(2126.833,962.252)(2180.498,989.085)
\put(2325,483){\makebox(0,0)[lb]{\smash{{{\SetFigFont{8}{9.6}{\rmdefault}{\mddefault}{\updefault}critical row}}}}}
\put(1650,708){\makebox(0,0)[lb]{\smash{{{\SetFigFont{5}{6.0}{\rmdefault}{\mddefault}{\updefault}$1$}}}}}
\put(3000,1008){\makebox(0,0)[lb]{\smash{{{\SetFigFont{5}{6.0}{\rmdefault}{\mddefault}{\updefault}$1$}}}}}
\put(3375,1383){\makebox(0,0)[lb]{\smash{{{\SetFigFont{8}{9.6}{\rmdefault}{\mddefault}{\updefault}row $n$}}}}}
\put(0,183){\makebox(0,0)[lb]{\smash{{{\SetFigFont{8}{9.6}{\rmdefault}{\mddefault}{\updefault}row $1$}}}}}
\end{picture}
}
\end{center}
\caption{The rows of the white checkers \lremind{pus}}
\label{pus}
\end{figure}

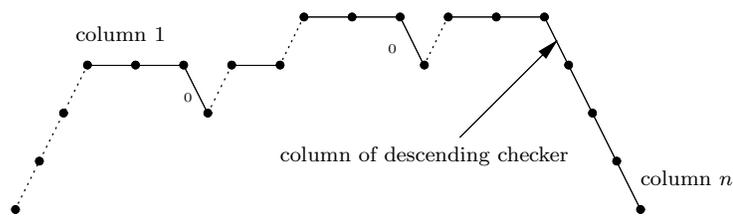
\begin{figure}[ht]
\begin{center}
\setlength{\unitlength}{0.00083333in}
\begingroup\makeatletter\ifx\SetFigFont\undefined%
\gdef\SetFigFont#1#2#3#4#5{%
  \reset@font\fontsize{#1}{#2pt}%
  \fontfamily{#3}\fontseries{#4}\fontshape{#5}%
  \selectfont}%
\fi\endgroup%
{\renewcommand{\dashlinestretch}{30}
\begin{picture}(3966,1281)(0,-10)
\put(33,33){\blacken\ellipse{50}{50}}
\put(33,33){\ellipse{50}{50}}
\put(183,333){\blacken\ellipse{50}{50}}
\put(183,333){\ellipse{50}{50}}
\put(333,633){\blacken\ellipse{50}{50}}
\put(333,633){\ellipse{50}{50}}
\put(483,933){\blacken\ellipse{50}{50}}
\put(483,933){\ellipse{50}{50}}
\put(783,933){\blacken\ellipse{50}{50}}
\put(783,933){\ellipse{50}{50}}
\put(1083,933){\blacken\ellipse{50}{50}}
\put(1083,933){\ellipse{50}{50}}
\put(1233,633){\blacken\ellipse{50}{50}}
\put(1233,633){\ellipse{50}{50}}
\put(1383,933){\blacken\ellipse{50}{50}}
\put(1383,933){\ellipse{50}{50}}
\put(1683,933){\blacken\ellipse{50}{50}}
\put(1683,933){\ellipse{50}{50}}
\put(1833,1233){\blacken\ellipse{50}{50}}
\put(1833,1233){\ellipse{50}{50}}
\put(2133,1233){\blacken\ellipse{50}{50}}
\put(2133,1233){\ellipse{50}{50}}
\put(2433,1233){\blacken\ellipse{50}{50}}
\put(2433,1233){\ellipse{50}{50}}
\put(2583,933){\blacken\ellipse{50}{50}}
\put(2583,933){\ellipse{50}{50}}
\put(2733,1233){\blacken\ellipse{50}{50}}
\put(2733,1233){\ellipse{50}{50}}
\put(3033,1233){\blacken\ellipse{50}{50}}
\put(3033,1233){\ellipse{50}{50}}
\put(3333,1233){\blacken\ellipse{50}{50}}
\put(3333,1233){\ellipse{50}{50}}
\put(3483,933){\blacken\ellipse{50}{50}}
\put(3483,933){\ellipse{50}{50}}
\put(3633,633){\blacken\ellipse{50}{50}}
\put(3633,633){\ellipse{50}{50}}
\put(3783,333){\blacken\ellipse{50}{50}}
\put(3783,333){\ellipse{50}{50}}
\put(3933,33){\blacken\ellipse{50}{50}}
\put(3933,33){\ellipse{50}{50}}
\dottedline{45}(33,33)(483,933)
\path(483,933)(1083,933)
\path(483,933)(1083,933)
\path(1083,933)(1233,633)
\path(1083,933)(1233,633)
\dottedline{45}(1233,633)(1383,933)
\path(1383,933)(1683,933)
\path(1383,933)(1683,933)
\dottedline{45}(1683,933)(1833,1233)
\path(1833,1233)(2433,1233)
\path(1833,1233)(2433,1233)
\path(2433,1233)(2583,933)
\path(2433,1233)(2583,933)
\dottedline{45}(2583,933)(2733,1233)
\path(2733,1233)(3333,1233)
\path(2733,1233)(3333,1233)
\path(3333,1233)(3933,33)
\path(3333,1233)(3933,33)
\path(2808,483)(3408,1083)
\path(2808,483)(3408,1083)
\blacken\path(3344.360,976.934)(3408.000,1083.000)(3301.934,1019.360)(3344.360,976.934)
\put(408,1083){\makebox(0,0)[lb]{\smash{{{\SetFigFont{8}{9.6}{\rmdefault}{\mddefault}{\updefault}column $1$}}}}}
\put(1083,708){\makebox(0,0)[lb]{\smash{{{\SetFigFont{5}{6.0}{\rmdefault}{\mddefault}{\updefault}$0$}}}}}
\put(2358,1008){\makebox(0,0)[lb]{\smash{{{\SetFigFont{5}{6.0}{\rmdefault}{\mddefault}{\updefault}$0$}}}}}
\put(3933,183){\makebox(0,0)[lb]{\smash{{{\SetFigFont{8}{9.6}{\rmdefault}{\mddefault}{\updefault}column $n$}}}}}
\put(1683,333){\makebox(0,0)[lb]{\smash{{{\SetFigFont{8}{9.6}{\rmdefault}{\mddefault}{\updefault}column of descending checker}}}}}
\end{picture}
}
\end{center}
\caption{The columns of the white checkers \lremind{push}}
\label{push}
\end{figure}

\begin{figure}[ht]
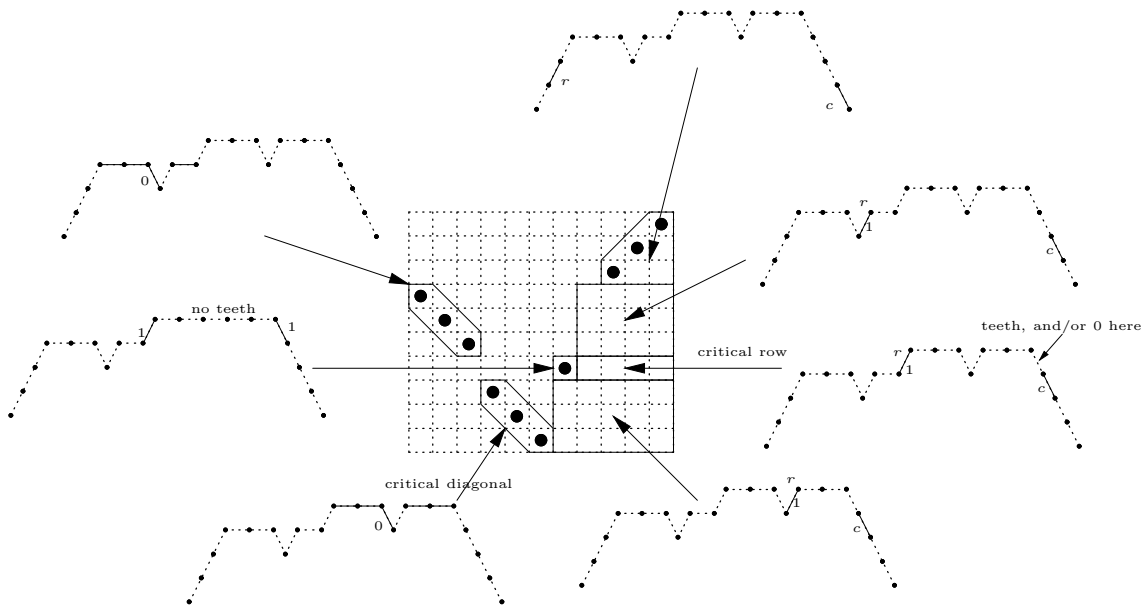

\begin{center}
\include{pummel}
\end{center}
\caption{How to locate the white checker given the partially completed puzzle ($r$ is the row and $c$ is the column; see Figures~\ref{pub} and~\ref{push}
to interpret them as numbers) \lremind{pummel}}
\label{pummel}
\end{figure}

We now go through the various cases of how to fill in the next part of
the puzzle, and verify that  they correspond to the possible next moves
of the checker games.  Each case is depicted in Figure~\ref{puff},
along with the portion of Table~\ref{keywest} 
that it corresponds to (in checkers).  The reader
should verify that all possible puzzle piece placements, and
all possible checker moves, are accounted for in the bijection.

{\em Case 1.}  
There is no white checker in the critical row, or in the next row.
Then make one move in the checker game.

{\em Case 2.}  There is no white checker in the critical row,
and there is a white checker in the next row, not on the rising black
checker.

{\em Case 3.}  There is no white checker in the critical row, and
there is a white checker on the rising black checker.

{\em Case 4.}  There is a white checker in the descending checker's square.
In this case, we finish the row of the puzzle, and make a series
of checker moves to move the descending checker to the bottom row.

{\em Case 5.}  There is a white checker in the critical
row but not in the descending checker's square, and there are no white checkers in any lower row.
We finish the row of the puzzle, and make a series
of checker moves to move the descending checker to the bottom row.

{\em Case 6.}  There is a white checker in the critical row, and
there is another white checker in a lower row, but in a higher
row than any white checkers on the critical diagonal (e.g.
a blacker if there is a white checker on the critical diagonal).
We finish the part of the row of the puzzle up to
the corresponding tooth, and make a series of checker moves to move the descending
checker to the blocker's row.

{\em Case 7.}  There is a white checker in the critical row but not
in the descending checker's square, and there is a white checker in the
rising checker's square.  Then we place two puzzle pieces and make one 
checker move, as shown.

{\em Case 8.}  There is a white checker in the critical row, but not
on the descending checker; there is a white checker in the
critical diagonal, but not on the rising checker; and there is no blocker.
Then there are two cases.  If the white checkers ``stay'', then
then we make one checker move, and place two pieces.  If the white checkers ``swap'', then we fill in the part of the puzzle until the ``1'' in the 
region marked $a$ in Figure~\ref{pun}, and make a series of checker moves
to move the descending checker to the row of the lower white checker
in question.\secretnote{In puff, ``Table 1'' was manually 
changed to ``Table keywest''.}

\begin{figure}[ht]
\begin{center}
\include{puff}
\end{center}
\caption{How to place the next piece in the puzzle?\lremind{puff}}
\label{puff}
\end{figure}

} % end of parskip; it started just before the introduction

\end{document}